\documentclass[reqno]{amsart}  
\theoremstyle{plain}
\newtheorem{theorem}{Theorem}[section]

\newtheorem{remark}[theorem]{Remark}
\newtheorem{proposition}[theorem]{Proposition}

\newtheorem{corollary}[theorem]{Corollary}
\newtheorem{example}[theorem]{Example}

\numberwithin{equation}{section}
\allowdisplaybreaks[1]

\theoremstyle{definition}

\newtheorem{definition}[theorem]{Definition}

\theoremstyle{remark}

\newcommand{\bbC}{{\mathbb C}}
\newcommand{\bA}{{\mathbf A}}
\newcommand{\ba}{{\mathbf a}}
\newcommand{\blam}{{\boldsymbol \lambda}}
\newcommand{\lam}{\lambda}
\newcommand{\bzeta}{{\boldsymbol \zeta}}

\newcommand{\bn}{{\mathbf n}}
\newcommand{\bm}{{\mathbf m}}

\newcommand{\cF}{{\mathcal F}}
\newcommand{\cH}{{\mathcal H}}
\newcommand{\cL}{{\mathcal L}}
\newcommand{\cM}{{\mathcal M}}

\newcommand{\cS}{{\mathcal S}}

\newcommand{\cU}{{\mathcal U}}
\newcommand{\cX}{{\mathcal X}}
\newcommand{\cY}{{\mathcal Y}}
\newcommand{\cO}{{\mathcal O}}

\newcommand{\bbZ}{{\mathbb Z}}
\newcommand{\C}{{\mathbb C}}
\newcommand{\B}{{\mathbb B}}

\begin{document}

\title[Backward-shift-invariant subspaces]{Multivariable
backward-shift-invariant subspaces and observability operators}
\author[J. A. Ball]{Joseph A. Ball}
\address{Department of Mathematics,
Virginia Tech,
Blacksburg, VA 24061-0123, USA}
\email{ball@math.vt.edu}
\author[V. Bolotnikov]{Vladimir Bolotnikov}
\address{Department of Mathematics,
The College of William and Mary,
Williamsburg VA 23187-8795, USA}
\email{vladi@math.wm.edu}
\author[Q. Fang]{Quanlei Fang}
\address{Department of Mathematics,
Virginia Tech,
Blacksburg, VA 24061-0123, USA}
\email{qlfang@math.vt.edu}

\begin{abstract}
         It is well known that subspaces of the Hardy space over the unit
         disk which are invariant under the backward shift occur as the
         image of an observability operator associated with a
         discrete-time linear system with stable state-dynamics, as well as
         the functional-model space for a Hilbert space contraction
         operator.  We discuss two multivariable extensions of this
         structure, where the classical Hardy space is replaced by (1) the
         Fock space of formal power series in a collection of $d$
         noncommuting indeterminates with norm-square-summable vector
         coefficients, and (2) the reproducing kernel Hilbert space (often
         now called the Arveson space) over the unit ball in ${\mathbb
         C}^{d}$ with reproducing kernel $k(\lam, \zeta) = 1/(1 - \langle
         \lam, \zeta \rangle)$ ($\lam, \zeta \in {\mathbb C}^{d}$ with $\|
         \lam \|, \| \zeta \| < 1$).  In the first case, the associated
         linear system is of noncommutative Fornasini-Marchesini type with
         evolution along a free semigroup with $d$ generators, while in the
         second case the linear system is a standard (commutative)
         Fornasini-Marchesini-type system with evolution along the integer
         lattice ${\mathbb Z}^{d}$.  An abelianization map (or
         symmetrization of the Fock space) links the first case with the
         second.  The second case has special features depending on whether
         the operator-tuple defining the state dynamics is commutative or
         not.  The paper focuses on multidimensional state-output linear systems
         and the associated observability operators; followup papers
         \cite{BBF2, BBF3} use the results here to extend the analysis to
         represent observability-operator ranges as reproducing kernel
         Hilbert spaces with reproducing kernels constructed from the
         transfer function of a conservative multidimensional
         (noncommutative or commutative) input-state-output linear system.

\end{abstract}


\maketitle

\tableofcontents

\section{Introduction}  \label{S:Intro}
\setcounter{equation}{0}
For $\cU$ and $\cY$ any pair of Hilbert spaces, we use the notation
$\cL(\cU, \cY)$ to denote the space of bounded, linear operators
from $\cU$ to $\cY$.  For $\cX$ a single Hilbert space, we shorten
the notation $\cL(\cX, \cX)$ to $\cL(\cX)$.
Let $\cX$, $\cU$ and $\cY$ be Hilbert spaces, let
$A \in \cL(\cX)$, $B \in \cL(\cU, \cX)$, $C \in \cL(\cX, \cY)$
and $D \in \cL(\cU, \cY)$ be bounded linear operators, and let us consider
the associated discrete-time linear time-invariant system
\begin{equation}
\left\{ \begin{array}{rcl}
x(n+1)&= & Ax(n)+Bu(n)\\
y(n)&= & Cx(n)+Du(n)\end{array} \right.
\label{1.1}
\end{equation}
with $x(n)$ taking values in the {\em state space} $\cX$, $u(n)$ taking
values in the {\em input-space} $\cU$ and $y(n)$ taking values in the
{\em output-space} $\cY$.
If we let the system evolve on the
nonnegative  integers $n \in {\mathbb Z}_{+}$, then the whole
trajectory $\{u(n), x(n), y(n)\}_{n \in {\mathbb Z}_{+}}$ is
determined from the input signal $\{u(n)\}_{n \in {\mathbb Z}_{+}}$
and the initial state $x(0)$ according to the formulas
\begin{align}
x(n) & = A^{n}x(0) + \sum_{k=0}^{n-1} A^{n-1-k} B u(k), \notag \\
y(n) & = C A^{n} x(0) + \sum_{k=0}^{n-1} C A^{n-1-k} B u(k) + D u(n)
\notag \\
         & = [{\mathcal O}_{C,A} x(0)]_{n} + \sum_{k=0}^{n-1} C A^{n-1-k} B
u(k) + D u(n)
\label{1.timedomainIO}
\end{align}
where ${\mathcal O}_{C,A}$ denotes the so-called {\em observability
operator}
$$
        {\mathcal O}_{C,A} \colon x \mapsto \left\{CA^nx\right\}_{n\in\bbZ_+}.
$$
If we introduce the $Z$-transform
\begin{equation}  \label{Ztrans}
\{f(n)\}_{n \in {\mathbb Z}_{+}} \mapsto \widehat f(\lambda) =
           \sum_{n=0}^{\infty} f(n) \lambda^{n},
\end{equation}
the $Z$-transformed version of the system-trajectory formulas
\eqref{1.timedomainIO} become
\begin{align}
           \widehat x(\lambda) & = (I - \lambda A)^{-1} x(0) + \lambda
           (I - \lambda A)^{-1} B \widehat
           u(\lambda),   \notag \\
           \widehat y(\lambda) & = C(I - \lambda A)^{-1} x(0) + [D + \lambda
           C(I - \lambda A)^{-1}B]
           \widehat u(\lambda) \notag \\
           & = \widehat {\mathcal O}_{C,A} x(0) + T_{\Sigma}(z) \widehat
           u(\lambda)
    \label{1.freqdomainIO}
\end{align}
where
\begin{equation}
\widehat{\cO}_{C,A} \colon \; x \mapsto
\sum_{n=0}^{\infty}(CA^n x)\, \lambda^n=C(I-\lambda A)^{-1} x
\label{1.3}
\end{equation}
is the $Z$-transformed version of the observability operator and where
$$
T_{\Sigma}(\lambda) = D + \lambda C (I - \lambda A)^{-1} B
$$
is the {\em transfer function} of the system $\Sigma$ given by \eqref{1.1}.
In particular, if the input signal $\{u(n)\}_{n \in {\mathbb Z}_{+}}$
is taken to be zero, the resulting output $\{y(n)\}_{n \in {\mathbb
Z}_{+}}$ is given by $y = {\mathcal O}_{C,A} x(0)$.
In case ${\mathcal O}_{C,A}$
is bounded as an operator from $\cX$ into $\ell^2_{\cY}:=\ell^{2}\otimes
\cY$ (here $\ell^{2}$ is the space of square-summable complex
sequences indexed by the nonnegative integers ${\mathbb Z}_{+}$, we say that
the pair $(C,A)$ is {\em output-stable}.
It is convenient to represent $\cO_{C,A}$ in the output-stable case in
the matrix form
$$
\cO_{C,A}=\operatorname{col}_{n\in\bbZ_+} \left[CA^{n}\right]
\colon \; {\mathcal X} \to\ell_{\cY}^{2}.
$$
        Since the {\em $Z$-transform} \eqref{Ztrans}
maps $\ell_{\cY}^{2}$ unitarily onto $H^2_{\cY}:=H^2(\mathbb
D)\otimes\cY$,
where
$H^{2}$, the image of $\ell^{2}$ under the $Z$-transform, is the
space of analytic functions on the unit disk with
modulus-square-summable sequence of Taylor coefficients:
$$
    H^{2} = \{ f(\lam) = \sum_{n=0}^{\infty} f_{n} \lam^{n} \colon
    \sum_{n=0}^{\infty} |f_{n}|^{2} < \infty\},
$$
the output stability of $(C,A)$ is equivalent to the
$Z$-transformed version of the observability operator \eqref{1.3}
being bounded as an operator from $\cX$ into $H_{\cY}^2$,
It is readily
seen that $\widehat{\cO}_{C,A} x =\widehat{\cO_{C,A} x}$.

If $(C,A)$  is output-stable, then the {\em observability gramian}
$$
{\mathcal G}_{C, A}:=(\cO_{C, A})^{*}\cO_{C,
A}=(\widehat{\cO}_{C,A})^*\widehat{\cO}_{C,A}
$$
is bounded on $\cX$ and can be represented via the series
\begin{equation}
{\mathcal G}_{C, A}= \sum_{n=0}^\infty A^{*n}C^*CA^n
\label{1.2}
\end{equation}
converging in the strong operator topology. The following result
gives a summary of well-known connections between output stability,
observability gramians and solutions of associated Stein equations
and inequalities.
\begin{theorem}
\label{T:1.1}
Let $(C,A)$ be a pair of operators with $C \colon
{\mathcal X} \to \cY$ and $A \colon {\mathcal X} \to {\mathcal X}$. Then:
\begin{enumerate}
           \item The pair $(C,A)$ is output-stable if and only if the Stein
           inequality
           \begin{equation}
           H-A^*HA \ge C^*C
           \label{Stein-ineq}
           \end{equation}
           has a positive semidefinite solution $H\in \cL(\cX)$.

           \item If $(C,A)$ is output-stable, then the observability gramian
           ${\mathcal G}_{C, A}$ satisfies the Stein equality
           \begin{equation}  \label{1.5}
	H - A^{*}HA = C^{*}C
           \end{equation}
            and is
           the minimal positive semidefinite solution of the Stein inequality
           \eqref{Stein-ineq}.

           \item There is a unique positive semidefinite solution of the
           Stein equality \eqref{1.5}
           if $A$ is {\em strongly stable}, i.e., powers $A^{n}$
           of $A$ tend to zero in the strong operator topology of
           ${\mathcal L}({\mathcal X})$.  If $A$ is a contraction operator,
           then the positive semidefinite solution of the Stein equation
           \eqref{1.5} is unique if and only if $A$ is strongly stable.
           \end{enumerate}
           \end{theorem}
A pair $(C,A)$ is called {\em observable} if the operator $\cO_{C,A}$
(equivalently, $\widehat{\cO}_{C,A}$, ${\mathcal G}_{C, \bA}$) is
injective. This
property means that a state space vector $x\in\cX$ is uniquely recovered
from the output string $\{y(n)\}_{n=0}^\infty$ generated by running the
system \eqref{1.1} with the zero input string and the initial condition
$x(0)=x$. A pair $(C,A)$ is called {\em exactly observable} if
$\cO_{C,A}$ (equivalently, ${\mathcal G}_{C, A}$) is bounded and bounded
from below.

Associated with an output-stable pair $(C,A)$ is the range of the
observability operator
$$
\operatorname{Ran} \widehat{\mathcal O}_{C,A} =
\{ C (I - zA)^{-1}x \colon x \in {\mathcal X}\}.
$$
The following theorem summarizes the connection between such ranges
and back\-ward-shift-invariant subspaces of $H^{2}_{\cY}$.

\begin{theorem} \label{T:1.2}
        Suppose that $(C,A)$ is an output-stable pair.  Then:
        \begin{enumerate}
\item
The linear manifold $\operatorname{Ran} \widehat {\mathcal O}_{C,A}$ is
invariant under the backward shift operator
\begin{equation}
\label{backwardshift}
S^{*} \colon f(\lambda) \to \frac{ f(\lambda) - f(0)}{\lambda}.
\end{equation}
\item  Let $H\ge 0$ be a solution of the Stein
inequality \eqref{Stein-ineq}
and let ${\mathcal X}'$ be the completion of ${\mathcal X}$ with
inner product $\| [x]\|_{{\mathcal X}'}^{2} = \langle H x, x
\rangle_{{\mathcal X}}$ (where $[x]$ denotes the equivalence class
modulo $\operatorname{Ker} H$ generated by $x$).
Then $A$ and $C$ extend to define bounded operators
$A' \colon {\mathcal X}' \to {\mathcal X}'$ and
$C' \colon {\mathcal X}' \to \cY$ and the observability operator $\widehat
{\mathcal O}_{C,A}$ extends to define a contraction operator
$\widehat {\mathcal O}_{C', A'}$
from ${\mathcal X}'$ into $H^{2}_{\cY}$. Moreover,
$\widehat {\mathcal O}_{C',A'} \colon {\mathcal X}' \to
H^{2}_{\cY}$ is an isometry if and only if $H$ satisfies the Stein
equation \eqref{1.5} and $A'$ is strongly
stable, i.e.,
$$
\langle H A^{n}x, A^{n}x \rangle \to 0\quad\mbox{for all} \; \;
x \in {\mathcal X}.
$$
\item If the linear manifold ${\mathcal M}:=
\operatorname{Ran} \widehat {\mathcal O}_{C,A}$ is given the lifted norm
$$
         \| \widehat {\mathcal O}_{C,A} x \|_{{\mathcal M}}^{2} =
        \inf_{y \in {\mathcal X} \colon {\mathcal O}_{C,A}y =
        {\mathcal O}_{C, A}x} \langle H y, y
         \rangle_{{\mathcal X}},
$$
then
\begin{enumerate}
\item ${\mathcal M}$ can be completed to ${\mathcal M}' =
\operatorname{Ran}{\mathcal O}_{C',A'}$ with contractive inclusion
in $H^{2}_{\cY}$:
$$
         \| f \|^{2}_{H^{2}_{\cY}} \le \| f \|^{2}_{{\mathcal M}'} \text{ for
         all } f \in {\mathcal M}'.
$$
Furthermore, ${\mathcal M}'$ is isometrically equal to the reproducing kernel
Hil\-bert space with reproducing kernel $K_{C,A;H}$ given by
\begin{equation}  \label{posker}
         K_{C,A;H}(\lambda, \zeta ) = C(I - \lambda A)^{-1}H (I -
         \overline{\zeta} A^{*})^{-1} C^{*}.
\end{equation}
\item The following difference-quotient inequality is satisfied
\begin{equation}
\label{1.dif-quot-ineq}
\| S^{*}f \|^{2}_{{\mathcal M}} \le \|f\|_{{\mathcal M}}^{2} -
\|f(0)\|^{2}_{\cY} \quad\mbox{for all} \; \; f\in{\mathcal M}
\end{equation}
and moreover, if the Stein equality \eqref{1.5} holds, then
\eqref{1.dif-quot-ineq} holds with equality.
\end{enumerate}
\item Conversely, if ${\mathcal M}$ is a Hilbert space contractively
included in $H^{2}_{\cY}$ which is invariant under $S^{*}$ and for which
the difference-quotient inequality \eqref{1.dif-quot-ineq} holds,
then there is a contractive pair $(C,A)$  (i.e., \eqref{Stein-ineq}
holds with $H = I_{{\mathcal X}}$) such that
${\mathcal M} = {\mathcal H}(K_{C,A;I}) = \operatorname{Ran} {\mathcal
O}_{C,A}$ isometrically. In case \eqref{1.dif-quot-ineq} holds with
equality, then $(C,A)$ can be taken to be isometric.
\end{enumerate}
\end{theorem}
Results of the type in Theorem \ref{T:1.1} are the basis for the
Lyapunov-function approach to stability analysis in system theory;
there are far-reaching generalizations to nonlinear and time-varying
settings which are far afield from our main interests here.
The goal of characterizing subspaces of $H^{2}$ of the form
$\cH(K_{C,A;H})$ (especially in a finite-dimensional context) was a key
feature in the approach to Nevanlinna-Pick interpolation developed by H.
Dym (see \cite{Dym}).

In this paper we present the analogues of Theorems \ref{T:1.1} and
\ref{T:1.2} for the two related multivariable settings:
(1) the case where the Hardy space $H^{2}$ on the unit disk is
replaced by the Fock space $H^{2}_{\cY}({\mathcal F}_{d})$, and (2)
the case where $H^{2}$ is replaced by the vector-valued Arveson
reproducing kernel Hilbert space ${\mathcal H}_{\cY}(k_{d})$.

To define the Fock space, we let ${\mathcal F}_{d}$ denote the free
semigroup on the set $\{1, \dots, d\}$ of the first $d$ natural numbers
and then let $H^{2}_{\cY}({\mathcal F}_{d})$ consist of the space of all
formal power series $\sum_{v \in \cF_{d}} f_{v} z^{v}$ in $d$
noncommuting indeterminates $z = (z_{1}, \dots, z_{d})$ with
coefficients $f_{v}$ in a coefficient Hilbert space $\cY$ which are
square-summable in norm: $\sum_{v \in \cF_{d}} \|f_{v}\|^{2}_{\cY}<
\infty$.  Here we write $z^{v} = z_{i_{N}} z_{i_{N-1}} \cdots z_{i_{1}}$ if
$v = i_{N}i_{N-1} \dots i_{1} \in \cF_{d}$.  The shift operator
$S  \colon f(\lambda)\mapsto \lambda f(\lambda)$ acting on the Hardy space
$H^{2}$ is replaced by the noncommuting $d$-tuple ${\mathbf S} = (S_{1},
\dots, S_{d})$ on $H^{2}_{\cY}(\cF_{d})$ given by
\begin{equation}  \label{NCshifts}
S_{j} \colon f(z) \mapsto f(z) z_{j} \; \; \text{ for } \; \; j = 1,
\dots,
d. \end{equation}
The system \eqref{1.1} is replaced by a noncommutative
multidimensional input-state-output system of the form
\begin{equation}
\left\{ \begin{array}{ccc}
x(1 v) &= & A_{1} x(v)+ B_1 u(v) \\
\vdots & \vdots  & \vdots  \\
x(d v) & = & A_d x(v)+ B_d u(v)  \\
y(v) &  = & Cx(v)+Du(v).\end{array}\right.
\label{2.3}
\end{equation}
Here the system evolves along the free semigroup $\cF_{d}$, and, for
each $v \in \cF_{d}$, the state vector $x(v)$, input signal $u(v)$
and output signal $y(v)$ take values in the {\em state space} $\cX$,
{\em input space} $\cU$ and {\em output space} $\cY$, and the {\em
system matrix} $U$ has the form
\begin{equation}  \label{sysmatrix}
   U = \begin{bmatrix} A & B \\ C & D \end{bmatrix} = \begin{bmatrix}
   A_{1} & B_{1} \\ \vdots & \vdots \\ A_{d} & B_{d} \\ C & D
\end{bmatrix} \colon \begin{bmatrix} \cX \\ \cU \end{bmatrix} \to
\begin{bmatrix} \cX \\ \vdots \\ \cX  \\ \cY \end{bmatrix}.
\end{equation}
Such systems were introduced in \cite{Cuntz-scat} and with further elaboration
in \cite{BGM1} and \cite{BGM2}; following \cite{BGM1} we call this 
type of system a
{\em noncommutative Fornasini-Marchesini linear system}.
The observability operator associated with an output map
$C \colon {\mathcal X} \to \cY$ and a
$d$-tuple $\bA =
(A_{1}, \dots, A_{d})$ of not necessarily commuting operators on a
Hilbert space ${\mathcal X}$, expressed in ``frequency-domain'' coordinates,
takes the form
$$ \widehat {\mathcal O}_{C, {\mathbf A}} \colon x \mapsto C (I - z_{1}A_{1} -
\cdots - z_{d} A_{d})^{-1} x.
$$
For the particular case where $\bA$ is a row contraction and
$$
C =(I-A_{1}^{*}A_{1} - \cdots - A_{d}^{*}A_{d})^{1/2}
$$
with $\cY$ taken
to be equal to the closure of the range of $C$,
this operator appears already in work of Popescu \cite{popjfa} under the
term ``Poisson kernel'' and as the adjoint of the key operator $L$
used in many constructions in the paper of Arveson \cite{Arveson}.
Reproducing kernel Hilbert spaces consisting
of formal power series were developed in a systematic way in
\cite{NFRKHS}.  Such spaces already appear (although not quite in our
notation) in the Sz.-Nagy-Foia\c s model theory for row contractions
developed by Popescu (see \cite{PopescuNF0, PopescuNF1, PopescuNF2}). We shall
see that Theorems \ref{T:1.1} and \ref{T:1.2} extend in a natural way
to this setting, where the observability gramian \eqref{1.2} in the
statement of Theorem \ref{T:1.1} is
replaced with the multivariable observability gramian
\begin{equation}  \label{NCobsgram}
       {\mathcal G}_{C,\bA} = \widehat{\mathcal O}^{*}_{C,\bA}
       \widehat{\mathcal O}_{C, \bA} = \sum_{v \in \cF_{d}} A^{v*}C^{*}C A^{v}
\end{equation}
(here we set $A^{v} = A_{i_{N}} \cdots A_{i_{1}}$ if $v = i_{N}
\cdots i_{1} \in {\mathcal F}_{d}$),
where the backward shift $S^{*}$
\eqref{backwardshift} in the statement of Theorem \ref{T:1.2} is
replaced by the $d$-tuple ${\mathbf S}^{*} = (S_{1}^{*}, \dots,
S_{d}^{*})$ of adjoints of the shift operators $S_{j}$ in
\eqref{NCshifts}, and where the positive kernel \eqref{posker}
becomes the kernel
\begin{equation}  \label{NCkernel}
          K_{C, \bA, H}(z,w) = C (I - z_{1} A_{1} - \cdots - z_{d}
          A_{d})^{-1} H (I - w_{1} A_{1}^{*} - \cdots - w_{d} A_{d}^{*})^{-1}
          C^{*}
\end{equation}
in two sets $z = (z_{1}, \dots, z_{d})$ and $w = (w_{1}, \dots,
w_{d})$ of noncommuting indeterminates
(see Theorems \ref{T:2-1.1} and \ref{T:2-1.2} below).

         In the second Arveson-space setting, the Hardy space $H^{2}$ over
         the unit disk is replaced by the so-called Arveson space, the
         reproducing kernel Hilbert space $\cH(k_{d})$ over the unit ball
         ${\mathbb B}^{d}$ in complex $d$-dimensional space ${\mathbb
         C}^{d}$ based on the reproducing kernel function
         $$
           k_{d}(\blam,\bzeta) = \frac{1}{1 - \langle \blam, \bzeta
\rangle_{{\mathbb C}^{d}}}
           \text{ for } \blam, \bzeta \in {\mathbb B}^{d}
         $$
         and the classical Hardy-space shift $f(\lambda) \mapsto \lambda \cdot
         f(\lambda)$ is
         replaced by the $d$-tuple of Arveson shift operators
          ${\mathbf M}_{\blam} = (M_{\lambda_{1}}, \dots, M_{\lambda_{d}})$
          where
         \begin{equation}  \label{Arveson-shift}
          M_{\lambda_{j}} \colon f(\blam) \mapsto \lambda_{j} f(\blam)
          \text{ for } f \in \cH(k_{d})
         \end{equation}
         (see \cite{Dr, ArvesonIII}).  In this case the underlying
	system evolves along the integer lattice ${\mathbb Z}^{d}_{+}
	= (\bn = (n_{1}, \dots, n_{d}) \colon n_{j} \in {\mathbb
	Z}_{+}\}$ has the form of what we call a  {\em (commutative)
	Fornasini-Marchesini system}
	\begin{equation}\left\{
	\begin{array}{rcl}
	x(\bn) & = & A_{1} x(\bn -e_{1}) + \cdots + A_{d} x(\bn - e_{d}) \\
	& & \qquad + B_{1}u(\bn - e_{1}) + \cdots + B_{d} u(\bn -
	e_{d}) \\
	y(\bn) & = & C x(\bn) + D u(\bn).
	\end{array} \right.
   \label{2.14}
\end{equation}
    Here and in what follows, $e_j$ denotes the element in
    ${\mathbb Z}^{d}_{+}$ having the $j$-th partial index equal to one
    and all other partial indices equal to zero:
   \begin{equation}
   e_j=(0,\ldots, 0,1,0, \ldots ,0)\in {\mathbb Z}^{d}_{+}.
   \label{2.16}
  \end{equation}
   Thus the system matrix $U$ has the same form \eqref{sysmatrix} as
   for the noncommutative setting but the domain for all the signals
   and the system evolution is the integer lattice ${\mathbb
   Z}^{d}_{+}$ rather than the free semigroup $\cF_{d}$ and the
   associated ``frequency-domain'' objects are functions or formal
   power series in the commuting variables $\blam = (\lam_{1}, \dots,
   \lam_{d})$ rather than in the noncommuting indeterminates $z =
   (z_{1}, \dots, z_{d})$.
   In Section \ref{S:Arveson}, we show how the Arveson space
         ${\mathcal H}(k_{d}) \otimes \cY$ and
         this Fornasini-Marchesini linear system can be derived as an
         abelianization (sometimes also called {\em symmetrization}) of the
         noncommutative Fock space $H^{2}_{\cY}(\cF_{d})$ and of
         the noncommutative Fornasini-Marchesini linear system, respectively;
         while it is well known that the Arveson space is a 
symmetrization of the
         Fock space and that the multiplier algebra on the Arveson space is
         the image under a completely positive map acting on the
         noncommutative multiplier algebra on the Fock space (see
         \cite{ArvesonIII, Arias-Popescu, DP, Davidson-Bordeaux} and
         \cite{Popescu-var1, Popescu-var2} for a recent, more general
         systematic framework),
         our extension of these
         ideas to the underlying system theory appears to be new.
         The observability operator, as in the noncommutative setting, is
         associated with a so-called output pair $(C, \bA)$ but now has the
         form
         $$ \widehat{\mathcal O}^{\ba}_{C, \bA} \colon x \mapsto C (I -
\lambda_{1}
         A_{1} - \cdots - \lambda_{d} A_{d})^{-1} x
         $$
         where the variables $\lambda_{1}, \dots, \lambda_{d}$ commute and
         the {\em abelianized} observability gramian ${\mathcal G}^{\ba}_{C,
         \bA} = (\widehat {\mathcal O}^{\ba}_{C, \bA})^{*}
         \widehat{\mathcal O}^{\ba}_{C, \bA}$
         has an infinite-series representation more complicated than the
         second expression in \eqref{NCobsgram} (see equation \eqref{2.23}
         below).
         In case the operator $d$-tuple  $\bA = (A_{1}, \dots,
         A_{d})$ is {\em commutative} (so $A_{i}A_{j} = A_{j} A_{i}$ for
         all $1 \le i,j \le d$), ${\mathcal G}^{\ba}_{C, \bA} = {\mathcal
         G}_{C,\bA}$ (see Proposition \ref{P:reverseStein} below),
          and Theorem \ref{T:1.2} has a natural analogue for
         this setting,
         with the abelianized multivariable observability gramian
         ${\mathcal G}^{\ba}_{C,\bA} = {\mathcal G}_{C, \bA}$
         playing the role of the observability gramian in Theorem \ref{T:1.1},
         with the operator $d$-tuple
         ${\mathbf M}_{\blam}^{*} = (M_{\lambda_{1}}^{*}, \dots,
         M_{\lambda_{d}}^{*})$ in place of the backward shift $S^{*}$
         \eqref{backwardshift} in Theorem \ref{T:1.2}, and with kernel
         \eqref{posker} now taken to be the multivariable positive kernel
         \begin{equation}  \label{MVkernel}
             K_{C, \bA; H}^{\ba}(\blam, \bzeta) =
             C (I - \lam_{1} A_{1} - \cdots - \lam_{d} A_{d})^{-1} H (I -
             \overline{\zeta_{1}} A_{1}^{*} - \cdots - \overline{\zeta_{d}}
             A_{d}^{*})^{-1} C^{*}
         \end{equation}
         (see Theorems \ref{T:3-1.2}, \ref{T:3-1.2c} and
\ref{T:3c-converse} below).
         In the general case where the operator $d$-tuple $\bA = (A_{1}, \dots,
         A_{d})$ is not assumed to be commutative, there is no
         characterization of the abelianized observability gramian as a
         minimal solution of a generalized Stein equation analogous to the
         classical case given in Theorem \ref{T:1.1}, but there still is a
         somewhat more implicit analogue of Theorem \ref{T:1.2}, where the
         backward shift $S^{*}$ \eqref{backwardshift} in Theorem
         \ref{T:1.2} is replaced by a solution of the so-called {\em
         Gleason problem} (see Theorems \ref{T:3-1.2nc} and
         \ref{T:3nc-converse} below).  The Gleason problem originates in
         the work of Henkin and Gleason (see \cite{gleason, Henkin}) and
         has been studied in the context of the Arveson space (with various
         formulas for the solution) in \cite{AlpayDubi} with an
         application to realization questions in \cite{ADR}.  Our
         analogue of Theorem \ref{T:1.2} for the Arveson space for the case
         of commutative $d$-tuple $\bA$ has already been given in \cite{BR1}
         (with a more general power-series setting worked out in
         \cite{BR2}) for the finite-dimensional case.

         We also give various numerical
         examples (constructed with the aid of the software program
         {\tt MATHEMATICA})
         to illustrate how $\widehat {\mathcal O}_{C, \bA}$
         and $\widehat {\mathcal O}^{\ba}_{C, \bA}$ can have divergent
         properties when ${\mathbf A}$ is not commutative (see
         Examples \ref{E:reverseStein}, \ref{E:a-stable}
         and \ref{E:a-obs} below).

	Backward-shift-invariant subspaces for the classical setting
	have been used for some time in the operator-theory
	literature as the model space for a more general (abstractly
	defined) Hilbert-space contraction operator (see \cite{deBR,
	NF}); connections of this work with linear system theory were
	only realized later (see e.g.~ \cite{Helton73, Helton74}).
	Our results develop the structure of such model spaces
	for the case of operator-tuples and therefore are of interest
	from the point of view of multivariable operator theory.  We
	find it satisfying that these model spaces in turn tie in
	with the theory of multidimensional linear systems in much
	the same way (but with some surprises) as in the classical
	case.

         As applications of the ideas, we obtain new system-theoretic
         derivations of the Beurling-Lax
         representation theorem for shift invariant subspaces in both the
         noncommutative and commutative settings; the result for the
         noncommutative setting is due originally to Popescu
         \cite{PopescuBL} and for the commutative setting to Arveson
         \cite{ArvesonIII} and McCullough-Trent \cite{MCT}. We also indicate
         connections with dilation theory and the von Neumann inequality for
         these settings (see \cite{Popescu-vN, popjfa, Dr, ArvesonIII}).

         Closely related to the kernels $K_{C,\bA}(z,w)$ and
         $K^{\ba}_{C,\bA}(\blam, \bzeta)$ (given by \eqref{NCkernel} and
         \eqref{MVkernel} with $H$ normalized to be the identity operator) are
         kernels
         of de Branges-Rovnyak type (see \cite{deBR} for the classical
         case)
         $$
         K_{S}(z,w) = k_{\text{Sz}}(z,w) - S(z) k_{\text{Sz}}(z,w)
         S(w)^{*}, \qquad K_{S}^{\ba}(\blam,\bzeta) = \frac{ I - S(\blam)
         S(\bzeta)^{*}}{1 - \langle \blam, \bzeta \rangle}
         $$
         (where $z = (z_{1}, \dots, z_{d})$ and $w = (w_{1}, \dots, w_{d})$
         are two sets of noncommuting indeterminates with
         $k_{\text{Sz}}(z,w) = \sum_{v \in \cF_{d}} z^{v} w^{v^{\top}}$
         equal to the {\em noncommutative Szeg\"o kernel} while $\blam =
         (\lam_{1}, \dots, \lam_{d})$ and $\bzeta = (\zeta_{1}, \dots,
         \zeta_{d})$ are two sets of commuting variables)
         for respective reproducing kernel Hilbert spaces
       $\cH(K_{S})$, $\cH(K_{S}^{\ba})$ in the respective noncommutative
and commutative settings.
     In this situation (where $K_{S}$ and $K_{S}^{\ba}$ are positive
     kernels in noncommuting and commuting variable, respectively),
     the respective power series
	$$
	S(z) = \sum_{v \in \cF_{d}} S_{v} z^{v}, \qquad
	S(\lam) = \sum_{\bn \in {\mathbb Z}^{d}_{+}} S_{\bn}
	\blam^{\bn},
	$$
         are {\em contractive multipliers}, i.e., the respective
         multiplication operators
         $$ M_{S} \colon f(z) \mapsto S(z) \cdot f(z), \qquad
         M_{S} \colon f(\blam) \mapsto S(\blam) \cdot f(\blam)
         $$
         are bounded from $H^{2}_{\cU}(\cF_{d})$ into
         $H^{2}_{\cY}(\cF_{d})$ and from $\cH_{\cU}(k_{d})$ into
         $\cH_{\cY}(k_{d})$ respectively with norm at most 1.
           A particular issue is the construction
         of operators
         $$ B_{j} \colon \cU \to \cX \quad\text{for} \; \; j = 1, \dots, d
         \quad\text{and} \quad D \colon \cU \to \cY
         $$
         for some input space $\cU$ so that
         \begin{align*}
           S(z) &= D + C (I - z_{1}A_{1} - \cdots -
           z_{d}A_{d})^{-1}(z_{1}B_{1} + \cdots + z_{d} B_{d}), \\
           S(\blam) & = D + C (I - \lam_{1} A_{1} - \cdots -
           \lam_{d}A_{d})^{-1} (\lam_{1}B_{1} + \cdots + \lam_{d}B_{d})
         \end{align*}
         satisfy
         $$
         K_{C,\bA}(z,w) = K_{S}(z,w), \qquad K_{C, \bA}^{\ba}(\blam,
         \bzeta) = K_{S}^{\ba}(\blam, \bzeta).
         $$
          With the resolution of this issue, then the results
         here lead directly to representations of backward-shift-invariant
         subspaces as reproducing kernel Hilbert spaces of the form
         $\cH(K_{S})$ and $\cH(K^{\ba}_{S})$ for a Schur multiplier $S$ in
         both the noncommutative and commutative settings as well as
         linear-fractional realizations for Beurling-Lax representers of
         shift-invariant subspaces for
         both the noncommutative (see \cite{PopescuBL}) and commutative
         (see \cite{MCT}) settings.  We work out these issues for the
         commutative setting and for the noncommutative setting
         in \cite{BBF2} and \cite{BBF3} respectively.

         The paper is organized as follows. After the present Introduction,
         Section \ref{S:Fock} focuses on the noncommutative Fock space
         setting while Section \ref{S:Arveson} focuses on the Arveson-space
         setting.  Section \ref{S:Fock} is divided into Section
         \ref{S:NC-Stein} dealing with the connections between solutions of
         generalized Stein equations and strong stability of the state
         dynamics for noncommutative Fornasini-Marchesini systems and
         Section \ref{S:NC-Obs} dealing with characterizing ranges of
         observability operators as  backward-shift-invariant subspaces of
         the Fock space with a certain reproducing-kernel-Hilbert-space
         structure.  The first subsection (Section \ref{S:C-Stein}) of
         Section \ref{S:Arveson} deals with the less tractable issues
         parallel to the material in Section \ref{S:NC-Stein} of
         generalized Stein equations and stability for commutative
         Fornasini-Marchesini
         systems and also presents the abelianization map giving the
         connection between noncommutative and commutative
         Fornasini-Marchesini systems.  The second subsection (Section
         \ref{S:C-Obs}) of Section \ref{S:Arveson}, parallel to Section
         \ref{S:NC-Obs}, discusses characterizations of
         observability-operator ranges for the case of a commutative
         Fornasini-Marchesini state-output system.  The results are the most
         satisfying in case the operator-tuple $\bA$ giving the state
         dynamics is commutative---these are collected in Subsection
         \ref{S:C-Obs-comA}.  The more implicit results for the case of
         noncommutative $\bA$ are given in Subsection \ref{S:C-Obs-noncomA}.

\section{The Fock-space setting}  \label{S:Fock}

\subsection{Output stability and Stein equations: the noncommutative case}
\label{S:NC-Stein}

For $d$ a positive integer, let ${\cF}_{d}$ be the free semigroup
${\cF}_{d}$ generated by the set of $d$ letters $\{1, \dots, d\}$.
Elements of ${\mathcal F}_{d}$ are words of the form $i_{N}
\cdots i_{1}$ where $i_{\ell} \in \{1, \dots, d\}$ for each $\ell = 1,
\dots, N$ with multiplication given by concatenation.  We also use
$\emptyset$ to denote the empty word; this serves as the unit element for
${\mathcal F}_{d}$.  For $v = i_{N} i_{N-1} \cdots i_{1} \in {\mathcal
F}_{d}$, we let $|v|$ denote the number $N$ of letters in $v$ and we let
$v^{\top} : = i_{1}  \cdots i_{N-1} i_{N}$ denote the {\em transpose}
of $v$. We let $z = (z_{1}, \dots,z_{d})$ to be a
collection of $d$ formal noncommuting variables and let $\cY\langle\langle
z\rangle\rangle$ denote the set of formal noncommutative series
$\sum_{v\in\cF_d} f_vz^v$ where $f_v\in\cY$ and where
\begin{equation}
\label{NCfunccalc}
z^{v} =z_{i_{N}}z_{i_{N-1}} \cdots
z_{i_{1}}\quad\mbox{if}\quad v= i_{N}i_{N-1} \cdots i_{1}.
\end{equation}
The Fock space $\ell_\cY^{2}(\cF_d)$ is defined as
\begin{equation}
\ell_{\cY}^{2}(\cF_d): = \left\{ \{f_v\}_{v\in\cF_d}
\colon \, \sum_{v\in\cF_d}\|f_v\|_{\cY}^{2} < \infty \right\}.
\label{focknorm}
\end{equation}
If we let $\chi_{v}$ be the characteristic function of the word $v$, so
$$
\chi_{v} = \{ \chi_{v}(v')\}_{v' \in \cF_{d}} \quad\text{where}\quad
          \chi_{v}(v') = \begin{cases} 1 &\text{if } v' = v, \\
             0 & \text{otherwise,}
             \end{cases}
          $$
          and we let ${\mathcal B}_{\cY}$ be an orthonormal basis for $\cY$,
          then $\{ \chi_{v} y_{i} \colon v \in \cF_{d}, y_{i} \in {\mathcal
          B}_{\cY}\}$ is an orthonormal basis for $\ell^{2}_{\cY}(\cF_{d})$.
          The space $\ell^{2}_{\cY}(\cF_{d})$ can be identified as the
          tensor product $\ell^{2}(\cF_d)\otimes\cY$
and is mapped unitarily onto the space
\begin{equation}
H_{\cY}^{2}(\cF_{d}) = \left\{\sum_{v \in \cF_{d}} f_{v}z^{v}
\in \cY\langle\langle z\rangle\rangle\colon \;
\sum_{v \in \cF_{d}} \|f_{v}\|_{\cY}^{2} <\infty \right\}
\label{focknorm1}
\end{equation}
by the noncommutative $Z$-transform
\begin{equation}  \label{NC-Z}
\{f_{v} \}_{v \in \cF_{d}} \mapsto f^{\wedge}(z) = \sum_{v \in
\cF_{d}} f_{v} z^{v}
\end{equation}
with the monomials $z^{v}$ playing the role of the basis vectors
$\chi_{v}$.

The noncommutative multidimensional analogue of the system \eqref{1.1} is the
system with evolution along the free semigroup $\cF_{d}$ given by
\eqref{2.3}.
Upon running the system \eqref{2.3} with the zero input string
$u(v)=0$ for $v\in\cF_d$ and a fixed
initial condition $x(\emptyset)=x\in\cX$ we get
\begin{equation}
y(v)=C\bA^vx + \sum_{v'', v' \in \cF_{d}, j \in \{1,
\dots, d \} \colon v'' j v' =v} C A^{v''}B_{j} u(v').
\label{2.4}
\end{equation}
Here we extend the noncommutative functional calculus \eqref{NCfunccalc}
from noncommuting indeterminates $z = (z_{1}, \dots, z_{d})$ to a
$d$-tuple of operators ${\mathbf A} = (A_{1}, \dots, A_{d})$;  we use
the notation
$$
         {\mathbf A}^{v} = A_{i_{N}}A_{i_{N-1}} \cdots A_{i_{1}}\quad
\text{if}\quad
         v = i_{N} i_{N-1} \cdots i_{1} \in \cF_{d}
$$
where the multiplication is now operator composition.
Application of the formal noncommutative $Z$-transform \eqref{NC-Z} then gives
\begin{equation}  \label{NC-IO}
\widehat y(z) = C (I - Z(z) A)^{-1} x(\emptyset) +
T_{\Sigma}(z) \widehat u(z)
\end{equation}
where the formal power series $T_{\Sigma}(z)$ (by definition equal to the {\em
transfer function} of the system \eqref{2.3}) is given by
$$
T_{\Sigma}(z) = D + C (I - Z(z) A)^{-1} Z(z) B
$$
where we have set
\begin{equation}
Z(z) = \begin{bmatrix} z_{1} & \cdots & z_{d}
\end{bmatrix}\otimes I_{\cX}, \quad
A=\begin{bmatrix} A_{1} \\ \vdots \\ A_{d} \end{bmatrix},
\quad B = \begin{bmatrix} B_{1} \\ \vdots \\ B_{d} \end{bmatrix}.
\label{2.8}
\end{equation}
For details see \cite{Cuntz-scat} or, for a more general setting of
structured noncommutative multidimensional systems, see \cite{BGM1}.

In analogy to the classical case, the system \eqref{2.3} is called {\em
output-stable} (and in this case we will say that the pair $(C,\bA)$ is
output-stable) if the output string $\{y(v)\}_{v\in\cF_d}$, defined as in
\eqref{2.4} but with the input string $\{u(v)\}_{v \in \cF_{d}}$ assumed
to be equal to 0, belongs to $\ell_\cY^{2}(\cF_d)$ for every $x\in\cX$ and the
observability operator
\begin{equation}
\cO_{C,\bA}\colon \; x\mapsto \left\{C\bA^vx\right\}_{v\in\cF_d}
\label{2.5}
\end{equation}
is bounded as an operator from $\cX$ into $\ell_\cY^{2}(\cF_d)$.
The $Z$-transformed version of $\cO_{C,\bA}$ is
$$
\widehat{\cO}_{C,\bA}\colon \; x\mapsto \sum_{v\in\cF_d}(C\bA^vx) \,
z^v\in\cY\langle\langle z\rangle\rangle
$$
and the following realization formula for $\widehat{\cO}_{C,\bA}$ is
immediate:
$$
\widehat{\cO}_{C,\bA}x=C(I-Z(z) A)^{-1} x.
$$
If $(C,\bA)$ is output-stable, then $\widehat{\cO}_{C,\bA}$ maps $\cX$
into $H_{\cY}^{2}(\cF_{d})$ and is bounded. In this case it makes sense to
introduce the {\em observability gramian}
\begin{equation}
{\mathcal G}_{C, \bA}:=({\cO}_{C, \bA})^{*}{\cO}_{C, \bA}
=(\widehat{\cO}_{C, \bA})^{*}\widehat{\cO}_{C, \bA}
\label{2.10}
\end{equation}
and its representation in terms of strongly converging series
\begin{equation}
{\mathcal G}_{C, \bA}=\sum_{v\in\cF_d}\bA^{*v^\top}C^*C\bA^{v}
\label{2.11}
\end{equation}
follows immediately by definition \eqref{2.5} of $\cO_{C,\bA}$ and the
formula \eqref{focknorm} for the norm in $\ell_\cY^{2}(\cF_d)$. The second
equality in \eqref{2.10} follows by definition of $\widehat{\cO}_{C,\bA}$
and the formula \eqref{focknorm1} for the norm in $H_\cY^{2}(\cF_d)$.
\begin{definition}
A pair $(C,\bA)$ is called {\em observable} if ${\mathcal G}_{C, \bA}$ is
positive definite and {\em  exactly observable} if
${\mathcal G}_{C, \bA}$ is strictly positive definite.
We say that the $d$-tuple of operators ${\mathbf A} = (A_{1}, \dots,
A_{d})$ is {\em strongly stable} if
\begin{equation}  \label{3.3}
           \lim_{N \to \infty} \sum_{v \in \cF_{d} \colon |v| = N} \|{\mathbf
           A}^{v} x\|^{2} \to 0 \quad \text{for all} \; \; x \in {\mathcal X}.
          \end{equation}
\label{D:2.1}
\end{definition}

We mention that the term {\em pure} rather than {\em strongly stable}
has been used in this context (see \cite{Arveson}), but we prefer the
present terminology since {\em pure} so as to avoid confusion with
the use of the term {\em pure} in the context of contractive
operator-valued functions (see \cite{NF}).

In analogy with the classical case one can introduce the unobservable
subspace
\begin{equation}
{\rm Ker} \, {\mathcal G}_{C,\bA}={\rm Ker} \, \cO_{C,\bA}=
\bigcap_{v\in\cF_d}{\rm Ker} \, C\bA^v.
\label{2.12}
\end{equation}
Thus, observability of $(C,\bA)$ means that ${\rm Ker} \, {\mathcal
G}_{C,\bA}$ is the zero subspace or that
\begin{equation}
C\bA^vx=0 \; \; (\forall v\in\cF_d) \; \; \Longrightarrow \; \; x=0.
\label{2.13}
\end{equation}
The following is the noncommutative Fock-space counterpart to Theorem
\ref{T:1.1}.

\begin{theorem}  \label{T:2-1.1}
	Let ${\mathbf A} = (A_{1}, \dots,
	A_{d}) \in \cL(\cX)^{d}$ and let $C \in \cL(\cX, \cY)$.
	Then the pair $(C, {\mathbf A})$
	is output-stable if and only if the
	(generalized) Stein inequality
\begin{equation} \label{3.4a}
	  H - A_{1}^{*}H A_{1} - \cdots - A_{d}^{*} H A_{d} \ge C^{*}C
\end{equation}
          has a positive semidefinite solution $H \in {\mathcal L}({\mathcal
X})$. In this case,
\begin{enumerate}
\item The observability gramian ${\mathcal G}_{C,\bA}$ satisfies the
          generalized Stein equation
          \begin{equation} \label{3.4}
	 H - A_{1}^{*}H A_{1} - \cdots - A_{d}^{*} H A_{d} = C^{*}C
	  \end{equation}
         and is the minimal
          positive semidefinite solution of the generalized Stein
          inequality \eqref{3.4a}.
\item The  positive semidefinite solution of  the Stein equation
          \eqref{3.4} is unique if  ${\mathbf A}$ is {\em strongly
          stable}, i.e., \eqref{3.3} holds.
          Moreover, in case ${\mathbf A}$ is contractive in the sense that
          \begin{equation}  \label{contractive}
	 A_{1}^{*}A_{1} + \cdots + A_{d}^{*}A_{d} \le I_{\cX},
          \end{equation}
          then the solution of the Stein equation \eqref{3.4} is unique
	if and only if ${\mathbf A}$ is strongly stable.
\end{enumerate}
         \end{theorem}

     \begin{proof} Suppose first that $(C, {\mathbf A})$
         is output-stable.  Then for each $x \in {\mathcal X}$,
         $$
         \{ C {\mathbf A}^{v} \}_{v \in \cF_{d}} \in
         \ell^{2}_{\cY}(\cF_{d}), \text{ i.e., }
         \sum_{v \in \cF_{d}} \| C {\mathbf A}^{v} x \|^{2}_{\cY} < \infty.
         $$
         This has the consequence that the infinite series
         $$
           \sum_{N = 0}^{\infty} \sum_{v \in \cF_{d} \colon |v| = N} {\mathbf
           A}^{v*} C^{*} C {\mathbf A}^{v}
         $$
         converges in the strong operator topology to an operator $H \in
         {\mathcal L}({\mathcal X})$ (in fact, $H = {\mathcal G}_{C, {\mathbf
         A}}$ is the observability gramian).  From this infinite-series
         representation for ${\mathcal G}_{C, {\mathbf A}}$ it is easily
         verified that ${\mathcal G}_{C, {\mathbf A}}$ is
         positive semidefinite and satisfies the Stein equation
         \eqref{3.4} and hence also the Stein inequality \eqref{3.4a}.

         Conversely, suppose that the Stein inequality \eqref{3.4}
         has a positive semidefinite solution $H$.  We first claim that
         \begin{equation}
         H \ge  \sum_{v \in \cF_{d} \colon |v| \le N}\bA^{*v^{\top}}C^*C\bA^{v}
         + \sum_{v\in\cF_{d} \colon |v|=N+1}\bA^{*v^{\top}}H\bA^{v}
         \label{3.5}
         \end{equation}
         for each $N \in {\mathbb Z}_{+}$.
         For $N=0$, \eqref
         {3.5} collapses to \eqref{3.4a} which is
         given. Inductively assume that
         $$
         H \ge \sum_{v\in \cF_{d}\colon |v| < N} \bA^{*v^{\top}}C^{*}C \bA^{v}
         + \sum_{v\in \cF_{d}\colon |v| = N} \bA^{*v^{\top}}H \bA^{v}.
         $$
         Use the Stein inequality \eqref{3.4a} to replace  $H$ on the
         right side by its lower bound
         $C^{*}C + A_{1}^{*}HA_{1} +\cdots + A_{d}^{*}H A_{d}$
         to get from this
         $$
             H \ge
             \sum_{v\in \cF_{d} \colon |v| <N} \bA^{*v^{\top}}C^{*}C \bA^{v}
             + \sum_{v\in \cF_{d} \colon |v| = N+1} \bA^{*v^{\top}} H \bA^{v}
             + \sum_{v\in \cF_{d} \colon |v| = N} \bA^{*v^{\top}}C^{*}C \bA^{v}
         $$
         which then simplifies to \eqref{3.5} as wanted.

        We rewrite \eqref{3.5} in the form
        \begin{equation}
        \sum_{v \in \cF_{d} \colon |v| < N}\bA^{*v^{\top}}C^*C
        \bA^{v}  \le H-\sum_{v\in\cF_{d} \colon
|v|=N}\bA^{*v^{\top}}H\bA^{v}\le H.
        \label{3.6}
        \end{equation}
        By letting $N\to\infty$ in \eqref{3.6} we conclude that the left
        hand side sum converges (weakly and therefore, since all the terms are
        positive semidefinite, strongly) to a bounded positive semidefinite
        operator. By \eqref{2.11},
        $$
        \lim_{N\to\infty}\sum_{v \in \cF_{d} \colon |v| < N}\bA^{*v^{\top}}C^*C
        \bA^{v}=\sum_{v \in \cF_{d}}\bA^{*v^{\top}}C^*C\bA^{v}=
        {\mathcal G}_{C, \bA}
        $$
        and passing to the limit in \eqref{3.6} as $N\to \infty$ gives
        ${\mathcal G}_{C, \bA}\le H$.  In particular  the
        operator ${\mathcal G}_{C,\bA}$ is bounded (since $H$ is) and
        therefore the pair $(C,\bA)$ is output-stable.
        \end{proof}

         \begin{proof}[Proof of (1):] As observed in the proof of the
	  first part of the theorem,
         from the infinite-series
         representation \eqref{2.11} it follows that ${\mathcal
         G}_{C,{\mathbf A}}$ satisfies the Stein equation \eqref{3.4}.
         If $H$ is any solution of the Stein inequality, the computation
         leading to \eqref{3.6} shows that $H$ satisfies \eqref{3.6}. By
         taking the limit as $N \to \infty$ we conclude that ${\mathcal
         G}_{C, {\mathbf A}} \le H$ as asserted.
         \end{proof}

         \begin{proof}[Proof of (2):]  Suppose that ${\mathbf A}$ is
         strongly stable and that $H$ solves the Stein equation \eqref{3.4}.
         Then the proof of \eqref{3.5} shows that in this case \eqref{3.5}
         holds with equality:
         \begin{equation}
           H =  \sum_{v \in \cF_{d} \colon |v| \le N}\bA^{*v^{\top}}C^*C\bA^{v}
           + \sum_{v\in\cF_{d} \colon |v|=N+1}\bA^{*v^{\top}}H\bA^{v}
           \label{3.5=}
           \end{equation}
          for each $N=0,1,2,\dots$.  Taking the limit as $N \to \infty$ and
          using the stability assumption \eqref{3.3} we conclude that
          $H = {\mathcal G}_{C,{\bf A}}$.

          For the converse direction we assume in addition that ${\mathbf
          A}$ is contractive (i.e., \eqref{contractive} holds).
          We prove the contrapositive:  {\em if ${\mathbf
          A}$ does not satisfy the stability condition \eqref{3.3}, then
          the solution of the Stein equation \eqref{3.4} is not unique.}
          Assume therefore that ${\mathbf A}$ is not stable.  By the
          assumption \eqref{contractive}, the sequence of operators
          $$
            \Delta_{N} = \sum_{v \in \cF_{d} \colon |v| = N} {\mathbf
A}^{*v^{\top}}
            {\mathbf A}^{v},\ N=1,2, \dots
          $$
          is decreasing and bounded below and therefore has a strong limit
          $\Delta$.  Since ${\mathbf A}$ is assumed not to be stable, this
          limit $\Delta$ is not zero.  However it is easily verified that
          \begin{equation}  \label{iterate}
          A_{1}^{*}\Delta_{N} A_{1} + \cdots + A_{d}^{*} \Delta_{N}
          A_{d} = \Delta_{N+1}.
          \end{equation}
          Taking limits in \eqref{iterate} gives that $\Delta = \lim_{N \to
          \infty} \Delta_{N}$ satisfies the homogeneous Stein equation
          $$ \Delta - A_{1}^{*} \Delta A_{1} - \cdots - A_{d}^{*} \Delta
          A_{d} = 0.
          $$
          We conclude that the solution of the Stein equation \eqref{3.4}
          cannot be unique.
          \end{proof}

          Particular cases of output pairs $(C, {\mathbf A})$ are the cases
          where $(C, {\mathbf A})$ is {\em contractive} (i.e., the Stein
          inequality \eqref{3.4a} holds with $H = I_{{\mathcal X}}$) and
where $(C,
          {\mathbf A})$ is {\em isometric} (i.e., the Stein equality \eqref{3.4}
          holds with $H = I_{{\mathcal X}}$).  For these cases some additional
          observations can be made along the lines of Theorem \ref{T:2-1.1}.

          \begin{proposition}  \label{P:2-1.1}
              \begin{enumerate}
	   \item Suppose that $(C, {\mathbf A})$ is a
              contractive pair.  Then $(C, {\mathbf A})$ is
              output-stable with ${\mathcal G}_{C, {\mathbf A}} \le
              I_{{\mathcal X}}$ and the observability gramian ${\mathcal G}_{C,
              {\mathbf A}}$ is the unique positive semidefinite solution of
              the Stein equation \eqref{3.4} if and only if ${\mathbf A}$ is
              strongly stable.

              \item Suppose that $(C, {\mathbf A})$ is an isometric pair.  Then
              $(C, {\mathbf A})$ is  output-stable. Moreover $H = I$
              is the unique solution of the Stein equation \eqref{3.4} if and
              only if ${\mathbf A}$ is strongly stable.  In this case
              ${\mathcal O}_{C, {\mathbf A}}$ is isometric and hence also $(C,
              {\mathbf A})$ is exactly observable.
              \end{enumerate}
         \end{proposition}

         \begin{proof}
         Statement (1) immediately follows from statements in Theorem
         \ref{T:2-1.1} combined with the observation that $(C, {\mathbf A})$
         being a contractive pair implies that ${\mathbf A}$ is contractive.

         The first two assertions in statement (2) follow in a similar way.
         As for the last assertion,
         for the case where $(C, \bA)$ is isometric, $I_{\cX}$ is a
         solution of the Stein equation \eqref{3.4}; for the situation
         where ${\mathbf A}$ is
         strongly stable, uniqueness implies that the observability
         gramian ${\mathcal G}_{C, \bA} = I_{\cX}$, i.e., that $\cO_{C,
         \bA}$ is isometric.  Then also $(C, \bA)$ is exactly observable by
         definition.
         \end{proof}

\begin{remark}  \label{R:obs-stable}  {\rm The converse of the last part
         of Proposition \ref{P:2-1.1} does not hold even for the case
         $d=1$.  More precisely}, there exists an isometric pair of
         operators
         $(C, A)$ such that $(C,A)$ is observable but $A$ is not strongly
         stable.

         {\rm An example necessarily requires that $\operatorname{dim}
         {\mathcal X} = \infty$.  In the terminology of Sz.-Nagy-Foia\c s
         \cite{NF}. it suffices to produce a {\em completely
         non-isometric} ({\em c.n.i.}) contraction operator $A$ on a
         nontrivial Hilbert space ${\mathcal X}$ (so $\operatorname{dim}
         {\mathcal X} > 0$ and
         there is no nonzero-invariant subspace ${\mathcal M}$ for $A$ such
         that $A|_{{\mathcal M}}$ is an isometry)
         in the class $C_{1 \cdot}$ (so $A^{n} x \to 0$ in ${\mathcal X}$ for
some $x \in {\mathcal X}$ implies that $x = 0$).  Indeed, if $A$ is such
         an operator, set $C = (I - A^{*}A)^{1/2}$ considered as an operator
         from ${\mathcal X}$ into $\cY: = \overline{\operatorname{Ran}} (I -
         A^{*}A)^{1/2}$ (the closure of the range of $(I - A^{*}A)^{1/2}$).
         Such an $A$ is not strongly stable by the definition of the class
         $C_{1 \cdot}$, the definition of $C$
         makes the pair $(C,A)$ isometric, and the condition that $A$ is
         c.n.i.~implies that $(C,A)$ is observable.

To construct such an operator $A$, let $\theta$ be a Schur-class outer
function such that $\log (1-|\theta|^{2})$ is not integrable
(with respect to arc-length Lebesgue measure) over the unit
circle ${\mathbb T}$.
Furthermore, let ${\mathbb K}(\theta)$ be the associated Sz.-Nagy-Foia\c s
model space
         $$
          {\mathbb K}(\theta) = \begin{bmatrix} H^{2}({\mathbb D}) \\
          L^{2}({\mathbb T}) \end{bmatrix} \ominus \begin{bmatrix}
          \theta(\lambda) \\ (1 - | \theta(\zeta)|^{2})^{1/2} \end{bmatrix}
          H^{2}({\mathbb D}) \; \text{ where }
          \lambda \in {\mathbb D} \text{ and } \zeta \in {\mathbb T},
         $$
and let $S(\theta)$ be the Sz.-Nagy-Foia\c s model operator
         $$
           T = \left. P_{{\mathbb K}(\theta)} \begin{bmatrix} M_{\lambda} & 0
\\ 0    & M_{\zeta}\end{bmatrix} \right|_{{\mathbb K}(\theta)},
         $$
where $M_{\lambda}$ and $M_{\zeta}$ are the operators of multiplication
by $\lambda$ and by $\zeta$, respectively. Now we let $A:=S(\theta)^{*}$
and note that $A$ is in the class $C_{1\cdot}$ by Proposition 3.5 in
\cite{NF} (since $\theta$ is outer) and $A$ is c.n.i. by Theorem 5 in
\cite{BK} (since the non-log-integrability property of $1 - |\theta|^{2}$
implies that there is no $H^{\infty}$-function $a(z)$ for which
$|a(\zeta)|^{2} \le 1 - |\theta(\zeta)|^{2}$ for $\zeta \in {\mathbb
T}$). This completes the construction.}
\end{remark}

           Let us say that  {\em the pair $(C, \bA)$ is similar to the
	    pair $(\widetilde C, \widetilde \bA)$ }
          if there is an invertible operator $S$ on $\cX$ so that
	$$
	 \widetilde C = C S^{-1}, \qquad \widetilde A_{j} = S A_{j}S^{-1}
	 \quad\text{for} \; \; j = 1, \dots, d.
	$$
Then we have the following characterization of pairs
$(C, {\mathbf A})$ which are similar to a contractive or to an isometric pair.

	\begin{proposition}  \label{P:similar}
	    \begin{enumerate}
         \item The pair $(C, \bA)$ is similar to a contractive pair
             $(\widetilde C, \widetilde \bA)$ if and only if there exists a
             bounded, strictly positive-definite solution $H$ to the
             Stein inequality \eqref{3.4a}.

             \item The pair $(C, \bA)$ is
             similar to an isometric pair if and only if there exists a bounded,
             strictly positive-definite solution $H$ of the  Stein
             equation \eqref{3.4}.
         \end{enumerate}

         \end{proposition}

         \begin{proof}
           Suppose that $H$ is a strictly positive-definite solution
        of \eqref{3.4a}.   Factor $H$ as $H = S^{*}S$ with
         $S$ invertible and set
             \begin{equation}  \label{tilde-pair}
             \widetilde C = C S^{-1}, \qquad \widetilde A_{k} = S A_{k}
             S^{-1}\quad\text{for} \; \; k = 1, \dots, d.
             \end{equation}
             Multiplying \eqref{3.4a} on the left by $S^{*-1}$ and on the right
             by $S^{-1}$ then leads us to
             $$
	 I - \widetilde A_{1}^{*} \widetilde A_{1} - \cdots - \widetilde
	 A_{d}^{*} \widetilde A_{d} \ge \widetilde C^{*} \widetilde C,
             $$
             i.e., $(\widetilde C, \widetilde \bA)$ is a contractive
pair which is
             similar to the original pair $(C, \bA)$.  Conversely, if
$(\widetilde
             C, \widetilde \bA)$ given by \eqref{tilde-pair} is 
contractive, then
             $H = S^{*}S$ is bounded and positive-definite and
satisfies the Stein
             inequality \eqref{3.4a}. This verifies the first statement of
             the Proposition.  The second statement follows in a similar way.
         \end{proof}

As a consequence of the observations in Proposition \ref{P:similar},
Proposition \ref{P:2-1.1} can be formulated more generally as follows.
\begin{proposition}  \label{P:2-1.1'}
	  \begin{enumerate}
	      \item If the pair $(C, {\mathbf A})$ is such that the Stein
	      inequality \eqref{3.4a} has a strictly positive-definite
	      solution $H$, then $(C, {\mathbf A})$ is output-stable.
	      Moreover, the observability gramian
	      ${\mathcal G}_{C, {\mathbf A}}$ is the unique positive
	      semidefinite solution of the Stein equation \eqref{3.4} if and
	      only if ${\mathbf A}$ is strongly stable.
          \item If the pair $(C, {\mathbf A})$ is such that the Stein
          equation \eqref{3.4} has a strictly positive-definite solution $H$,
          then $(C, {\mathbf A})$ is output-stable and the observability gramian
        ${\mathcal G}_{C, {\mathbf A}}$ is the unique positive
        semidefinite solution of the Stein equation \eqref{3.4} if and
        only if ${\mathbf A}$ is strongly stable.  In this
        case $(C, {\mathbf A})$ is moreover exactly observable.
         \end{enumerate}
\end{proposition}

The last part of Proposition \ref{P:2-1.1'} has a converse.

\begin{proposition} \label{P:2-1.1'converse}
           Suppose that the pair $(C, {\mathbf A})$ is output-stable and
           exactly observable.  Then ${\mathbf A}$ is strongly stable, i.e.,
           \eqref{3.3} holds.
        \end{proposition}

        \begin{proof}
            If $(C, {\mathbf A})$ is output-stable and exactly observable,
            then the observability gramian ${\mathcal G}_{C, {\mathbf A}}$
            is a strictly positive-definite solution of the Stein equation
            \eqref{3.4}.  Hence \eqref{3.5=} holds with $H = {\mathcal G}_{C,
            {\mathbf A}}$:
           \begin{equation}  \label{3.5=gram}
        \langle {\mathcal G}_{C, {\mathbf A}} x, x \rangle =
        \sum_{v \in \cF_{d} \colon |v| \le N}
        \langle \bA^{*v^{\top}}C^*C\bA^{v} x, x \rangle
	   + \sum_{v\in\cF_{d} \colon |v|=N+1}\langle
	   {\mathcal G}_{C, {\mathbf A}}\bA^{v} x, \bA^{v}x \rangle.
\end{equation}
            From the infinite-series representation \eqref{2.11}
            for ${\mathcal G}_{C, {\mathbf A}}$, taking limits in
            \eqref{3.5=gram} gives
        \begin{equation}  \label{gram-lim}
        \lim_{N \to \infty}
        \sum_{v\in\cF_{d} \colon |v|=N+1}
        \langle {\mathcal G}_{C, {\mathbf A}}{\mathbf A}^{v}x, {\mathbf A}^{v} x
        \rangle  = 0.
\end{equation}
The strict positive-definiteness of ${\mathcal G}_{C, {\mathbf A}}$
tells us that there is an $\varepsilon > 0$ so that
\begin{equation}  \label{coercive}
        \varepsilon  \| x \|^{2} \le \langle {\mathcal G}_{C, {\mathbf A}} x,
        x \rangle\quad  \text{for all} \; \; x \in {\mathcal X}.
\end{equation}
In particular, from \eqref{coercive} with
${\mathbf A}^{v}x$ in place of $x$ combined with \eqref{gram-lim} we get
$$
        \varepsilon \sum_{v \in \cF_{d} \colon |v| = N+1} \| {\mathbf
A}^{v} x \|^{2}
          \le \sum_{v\in\cF_{d} \colon |v|=N+1}
        \langle {\mathcal G}_{C, {\mathbf A}}{\mathbf A}^{v} x, {\mathbf A}^{v}x
        \rangle \to 0
$$
for all $x \in {\mathcal X}$, and we conclude that ${\mathbf A}$ is
strongly stable as asserted.
        \end{proof}

        \subsection{Observability-operator range spaces and reproducing
        kernel Hilbert spaces: the noncommutative-variable case}
        \label{S:NC-Obs}

        To develop the noncommutative analogue of Theorem \ref{T:1.2}, we
        first introduce the right noncommutative shift operators
        $S^{R}_{1}, \dots,
        S^{R}_{d}$ on $H_{\cY}^{2}({\mathcal F}_{d})$ as follows:
        \begin{equation}
        S^{R}_j: \; \sum_{v\in\cF_d}f_vz^v\mapsto
        \sum_{v\in\cF_d}f_vz^vz_j=
        \sum_{v\in\cF_d}f_vz^{vj}\quad (j=1,\ldots,d).
        \label{4.5}
        \end{equation}
        It is readily seen that their adjoints (backward shifts) are given by
        \begin{equation}
        (S^{R}_j)^*: \; \sum_{v\in\cF_d}f_vz^v\mapsto
        \sum_{v\in\cF_d}f_{vj}z^v\quad (j=1,\ldots,d).
        \label{4.6}
        \end{equation}
        Their left counterparts $S^{L}_{1}, \dots, S^{L}_{d}$, also on
        $H^{2}_{\cY}(\cF_{d})$, are given by
   \begin{equation} \label{4.5L}
       S^{L}_{j} \colon \sum_{v \in \cF_{d}} f_{v} z^{v} \mapsto
       \sum_{v \in \cF_{d}} f_{v} z_{j} z^{v} = \sum_{v \in \cF_{d}}
       f_{v} z^{j v}
   \end{equation}
   with adjoints given by
   \begin{equation}  \label{4.6L}
       (S^{L}_{j})^{*} \colon \sum_{v \in \cF_{d}} f_{v }z^{v} \mapsto
       \sum_{v \in \cF_{d}} f_{j v} z^{v}.
   \end{equation}
     Let $\tau$ denote the unitary involution on
         $H^2_{\cY}(\cF_d)$ given by
         \begin{equation}
	 \tau \colon \; \sum_{v\in\cF_d}f_v z^v\to
\sum_{v\in\cF_d}f_{v^\top} z^v.
         \label{4.9a}
         \end{equation}
         In addition to the unitary property $\tau^{*} = \tau^{-1}$ of
         $\tau$, note also that $\tau$ intertwines the left shifts
	with the right shifts:
         \begin{equation}  \label{Sjtau}
             (S^{R}_{j})^{*} \tau = \tau (S^{L}_{j})^{*}, \quad
	   S_{j}^{R} \tau = \tau S_{j}^{L}\quad
	   \text{for}\quad j = 1, \dots, d.
         \end{equation}
        Then we have the following Fock-space analogue of Theorem \ref{T:1.2}.

\begin{theorem}
\label{T:2-1.2}
Suppose that $(C, {\mathbf A})$ is an output-stable pair.  Then:
\begin{enumerate}
\item The intertwining relation
\begin{equation}
\label{4.8a}
(S^{R}_{j})^{*} \widehat{\mathcal O}_{C, {\mathbf A}}x = \widehat{\mathcal
O}_{C, {\mathbf A}} A_{j} x\quad (x\in\cX)
\end{equation}
holds for every backward-shift operator $(S^{R}_j)^*$ defined in \eqref{4.6}
and hence $\operatorname{Ran}\widehat{\mathcal O}_{C, {\mathbf A}}$
is
$(S^{R}_j)^*$-invariant
for $j = 1, \dots, d$.
\item  Let $H\ge 0$ be a solution of the Stein inequality \eqref{3.4a}
             and let ${\mathcal X}'$ be the completion of ${\mathcal X}$ with
             $H$-inner product $\| [x]\|_{{\mathcal X}'}^{2} = \langle H x, x
             \rangle_{{\mathcal X}}$. Then $A_{j}$ and $C$ extend to define
              bounded operators $A'_{j} \colon {\mathcal X}' \to {\mathcal X}'$
              for $j = 1, \dots, d$ and $C' \colon {\mathcal X}' \to \cY$
              and the observability operator $\widehat
              {\mathcal O}_{C,{\mathbf A}}$ extends to define a contraction
              operator $\widehat{{\mathcal O}}_{C', {\mathbf A}'}$ from
${\mathcal X}'$
              into $H^{2}_{\cY}(\cF_{d})$. Moreover, $\widehat {\mathcal
O}_{C',{\mathbf A}'} \colon{\mathcal X}' \to H^{2}_{\cY}(\cF_{d})$
              is an isometry if and only if $H$ satisfies the Stein equation
              \eqref{3.4} and ${\mathbf A}' = (A'_{1}, \dots, A'_{d})$ is
              strongly stable, i.e.,
               \begin{equation}
               \label{H-stable}
               \sum_{v \in \cF_{d} \colon |v| = N} \langle H A^{v}x,
               A^{v}x \rangle \to 0 \quad\text{for all}\; \;  x \in 
{\mathcal X}.
               \end{equation}

\item If $H \ge 0$ is a solution of the Stein inequality \eqref{3.4a}
and the linear manifold ${\mathcal M}:=\operatorname{Ran} \widehat
            {\mathcal O}_{C,{\bf A}}$ is given  the lifted norm
            \begin{equation}  \label{lifted}
             \| \widehat {\mathcal O}_{C,{\bf A}} x \|_{{\mathcal M}}^{2} =
             \inf_{y \in {\mathcal X} \colon
             \widehat{\mathcal O}_{C,{\mathbf A}}y = \widehat{\mathcal O}_{C,
             {\mathbf A}}x} \langle H y, y \rangle_{{\mathcal X}},
             \end{equation}
             then
             \begin{enumerate}
	   \item
       ${\mathcal M}$ can be completed to ${\mathcal M}'=
\operatorname{Ran} \widehat{\mathcal O}_{C',{\bf A}'}$
(with $(C', \bA')$ as in \#2 above)
with
contractive inclusion in $H^{2}_{\cY}(\cF_{d})$:
$$
\| f \|^{2}_{H^{2}_{\cY}(\cF_{d})} \le
\| f \|^{2}_{{\mathcal M}'}\quad\text{for all} \quad f\in {\mathcal M}'.
$$
Furthermore, ${\mathcal M}'$ is isometrically equal to the formal
noncommutative reproducing kernel Hilbert space with reproducing kernel
$K_{C,{\bf A};H}$ given by \eqref{NCkernel}.

\item The following difference-quotient inequality is valid
\begin{equation}
\label{dif-quot-ineq}
\sum_{j=1}^{d}\| (S^{R}_{j})^{*}f \|^{2}_{{\mathcal H}(K_{C, {\mathbf A};H})}
\le \|f\|_{{\mathcal H}(K_{C,{\mathbf A};H})}^{2} -\|f_{\emptyset}\|^{2}_{\cY}
\end{equation}
for every $f\in {\mathcal M}' = {\mathcal H}(K_{C,{\mathbf A};H})$
with equality holding in \eqref{dif-quot-ineq} if and only if
\eqref{3.4} holds.
\end{enumerate}

\item Conversely, if ${\mathcal M}$ is a Hilbert space
included in $H^{2}_{\cY}(\cF_{d})$ which is invariant under
$(S^{R}_{j})^{*}$ for
$j = 1, \dots, d$ such that the difference-quotient inequality
\begin{equation}
\label{dif-quot-ineq1}
\sum_{j=1}^{d}\| (S^{R}_{j})^{*}f \|^{2}_{{\mathcal M}}
\le \|f\|_{{\mathcal M}}^{2} -\|f_{\emptyset}\|^{2}_{\cY}
\end{equation}
holds for every $f\in{\mathcal M}$,
then $\cM$ is contractively included in $H^{2}_{\cY}(\cF_{d})$ and
there exists a contractive pair $(C,{\mathbf A})$ (so $H = I$
positive definite
solution of the Stein inequality \eqref{3.4a}) such that
$$
{\mathcal M} = {\mathcal H}(K_{C,{\bf A};I}) = \operatorname{Ran}
\widehat{\mathcal O}_{C,{\bf A}}
$$
isometrically.  In case \eqref{dif-quot-ineq1} holds with
equality, then $(C,{\mathbf A})$ can be taken to be an isometric pair.
An explicit $(C, {\mathbf A})$ meeting these conditions is given as
         follows.  Take ${\mathcal X}$ to be the Hilbert space ${\mathcal
         X} = \tau ({\mathcal M})$ (where $\tau$ is the involution given by
         \eqref{4.9a}) with $\|\tau(f) \|_{{\mathcal X}} = \| f
         \|_{{\mathcal M}}$ and define $C
               \colon {\mathcal X} \to \cY$ and
               ${\mathbf A} = (A_{1}, \dots, A_{d})$ on ${\mathcal X}$  by
               \begin{equation} \label{modelAC}
                A_{j} = (S^{L}_{j})^{*}|_{{\mathcal X}} \text{ for } j = 1,
\dots, d; \qquad
                C =E|_{{\mathcal X}} \colon {\mathcal X} \to \cY, \qquad
                 H = I_{{\mathcal X}}
               \end{equation}
                where $E \colon H^{2}_{\cY}(\cF_{d}) \to \cY$ is given by
                \begin{equation}  \label{4.10a}
              E \colon \sum_{v \in \cF_{d}} f_{v} z^{v} \mapsto f_{\emptyset}.
              \end{equation}
          \end{enumerate}
\end{theorem}
        \begin{proof}[Proof of (1):]
            Applying $(S^{R}_j)^*$ to a typical element from $\operatorname{Ran}
            \widehat{\cO}_{C, \bA}$ (with notation as in \eqref{2.8}), we get
            \begin{eqnarray}
            (S^{R}_j)^*(C(I-Z(z)A)^{-1}x)&=&
            (S^{R}_j)^*(\sum_{v\in\cF_d}(C\bA^vx)z^v)\nonumber\\
            &=&\sum_{v\in\cF_d}(C\bA^{vj}x)z^v\nonumber\\
            &=&C(I-Z(z)A)^{-1}A_jx\in{\rm Ran} \,
            \widehat{\cO}_{C, \bA}.\label{4.7a}
            \end{eqnarray}
            The latter equality shows that ${\rm Ran} \,
            \widehat{\cO}_{C, \bA}$ is invariant under $(S^{R}_j)^*$
            for all  $j=1,\ldots,d$ (backward-shift-invariant) and
            \eqref{4.8a} follows.
            \end{proof}

        \begin{proof}[Proof of (2):]
            The Stein inequality \eqref{3.4a} amounts to the statement
            that $(C, {\mathbf A})$ is contractive  and well-defined on the
            dense subset
            $[{\mathcal X}]$ of ${\mathcal X}'$ (where $[x]$ is the
	 equivalence class containing $x$) and hence extends to a
            contractive pair $(C', {\mathbf A}')$ on all of ${\mathcal X}'$
            and moreover the inequality \eqref{3.5} holds for all $N = 1,2,
            \dots$.  From this we see that $\widehat {\mathcal
            O}_{C,{\mathbf A}}$ is contractive from ${\mathcal X}$ with the
            $H$-inner product to $H^{2}_{\cY}(\cF_{d})$,
            and hence also $\widehat {\mathcal O}_{C',
            {\mathbf A}'}$ is contractive from ${\mathcal X}'$ to
            $H^{2}_{\cY}(\cF_{d})$. The inequality \eqref{3.5} is actually a
            chain of inequalities $W_{N} \ge W_{N+1}$ for $N=1,2, \dots$,
            where
$$
W_{N} =\sum_{v \in \cF_{d} \colon |v| \le N-1}\bA^{*v^{\top}}C^*C\bA^{v}
             + \sum_{v\in\cF_{d} \colon |v|=N}\bA^{*v^{\top}}H\bA^{v}.
$$
       Note that
         $$
       \operatorname{s-lim}_{N \to \infty} W_{N} = {\mathcal
         G}_{C,{\mathbf A}} + \Delta_{H,{\mathbf A}}
$$
where
$$
\Delta_{H,{\mathbf A}} = \operatorname{s-lim}_{N \to
          \infty} \sum_{|v|=N} {\mathbf A}^{*v^{\top}} H {\mathbf A}.
$$
         In particular it follows from \eqref{contractive} that
         $$
         H \ge C^{*}C + \sum_{j=1}^{d} A_{j}^{*}H A_{j} \ge W_{N}\quad \text{for
         all} \; \; N \ge 2
         $$
         and hence, by taking the strong limit on the right hand side, we get
         \begin{equation} \label{ineq-chain}
             H \ge C^{*}C + \sum_{j=1}^{d} A_{j}^{*}H A_{j} \ge {\mathcal
         G}_{C, {\mathbf A}} + \Delta_{H,{\mathbf A}}.
         \end{equation}
         By definition, $\widehat{\mathcal O}_{C',{\mathbf A}'} \colon
{\mathcal X}'
         \to H^{2}_{\cY}(\cF_{d})$ being an isometry means
         that ${\mathcal G}_{C, {\mathbf A}} = H$ in which case
         \eqref{ineq-chain} becomes
         \begin{equation}  \label{ineq-chain'}
          H \ge C^{*}C + \sum_{j=1}^{d} A_{j}^{*} H A_{j} \ge H +
          \Delta_{H,{\mathbf A}}
         \end{equation}
         which in turn forces $\Delta_{H,{\mathbf A}} = 0$ and equalities
         throughout \eqref{ineq-chain'}.
         The condition $\Delta_{H, {\mathbf A}}=0$ just means that
         ${\mathbf A}'$ is strongly stable.
         From  equality holding in \eqref{ineq-chain'} we see that the Stein
         inequality \eqref{3.4a} holds with equality, i.e., the Stein
         equation \eqref{3.4} holds.  Conversely, by reversing the steps of
         the argument, we see that ${\mathbf A}'$ being strongly stable and
         the Stein equality holding leads to ${\mathcal G}_{C, {\mathbf A}} =
         H$, i.e., to $\widehat{\mathcal O}_{C', {\mathbf A}'}$ being 
an isometry
         from ${\mathcal X}'$ into $H^{2}_{\cY}(\cF_{d})$.
            \end{proof}

        \begin{proof}[Proof of (3a):]  Statement (3a) follows from general
        principles laid out in \cite{NFRKHS}.
          \end{proof}

          \begin{proof}[Proof of (3b):]
              For $f$ of the form $f(z) = C (I - Z(z) A)^{-1}x$, we have
          $$
	 \| f \|^{2}_{{\mathcal H}(K_{C, {\mathbf A}; H})} = \langle H x, x
              \rangle.
              $$
             We see that then
              $$
              \| S_{j}^{*} f \|^{2}_{{\mathcal H}(K_{C, {\mathbf A}; H})} =
              \langle H A_{j} x, A_{j} x \rangle \text{ for } \;  j = 1, \dots,
              d \text{ (from \eqref{4.7a})},  \qquad f_{\emptyset} = C x.
              $$
	 With these substitutions, we see that \eqref{dif-quot-ineq} is
	 equivalent to
	$$ \sum_{j=1}^{d} \langle H A_{j} x, A_{j} x \rangle_{\cX} \le
	\langle H x, x \rangle_{\cX} - \| C x \|^{2}_{\cY},
	$$
	or, in operator-theoretic form,
	\begin{equation} \label{Steinineq'}
	A_{1}^{*}H A_{1} + \cdots + A_{d}^{*}H A_{d} \le H - C^{*}C,
	\end{equation}
	with equality in \eqref{dif-quot-ineq} equivalent to equality in
	\eqref{Steinineq'}.  This completes the verification of part (3b) of
	Theorem \ref{T:2-1.2}.
	\end{proof}

         Before commencing the proof of part (4) of Theorem \ref{T:2-1.2},
         we collect some useful facts concerning $H^{2}_{\cY}(\cF_{d})$ itself.

             \begin{proposition}  \label{P:4.2}
             Let ${\mathbf S} = (S_1,\ldots,S_d)$ denote either the
	   right shift ${\mathbf S}^{R}$ or the left shift ${\mathbf
	   S}^{L}$
	   $$ {\mathbf S}^{R} = (S_{1}^{R}, \dots, S_{d}^{R}), \qquad
	   {\mathbf S}^{L} = (S^{L}_{1}, \dots, S^{L}_{d})
	   $$
	   defined as in \eqref{4.5} and \eqref{4.5L} and let
              the operator $E\colon \, H^2_{\cY}(\cF_d),\to \cY$ be 
defined as in
             \eqref{4.10a}.
             Then:
	\begin{enumerate}
	    \item
	    The operator-tuple ${\mathbf S}^{*} = (S_{1}^{*},
	    \dots, S_{d}^{*})$ is strongly stable, i.e.,
	\begin{equation}  \label{S*ss}
	    \lim_{N \to \infty} \sum_{v \in \cF_{d} \colon |v| = N}
	    \| {\mathbf S}^{* v} f \|^{2} = 0 \quad\text{for each}\quad f
\in
	    H^{2}_{\cY}(\cF_{d}).
	 \end{equation}

	 \item The operator
             $$
             \begin{bmatrix}S_1^*\\ \vdots \\ S_d^* \\ E\end{bmatrix}: \;
             H^2_{\cY}(\cF_d)\to (H^2_{\cY}(\cF_d))^d\oplus \cY
             $$
             is unitary, i.e.,
             \begin{equation}
              EE^*=I_{\cY},\quad ES_i=0,\quad S_j^*S_i=\delta_{ij}I\quad
             \text{for} \quad i,j=1,\ldots,d
             \label{4.11a}
	  \end{equation}
	  (where $\delta_{ij}$ stands for the Kronecker symbol), and
	  \begin{equation}
	      I-S_1S_1^*- \cdots -S_dS_d^*=E^*E.
             \label{4.12a}
             \end{equation}

	  \item
              $X = I$ is the unique solution of the Stein equation
           \begin{equation}  \label{NC-model-Stein}
	  X - S_{1}X S_{1}^{*} - \cdots - S_{d} X S_{d}^{*} = E^{*}E.
            \end{equation}

	\item
             For every $f\in H^2_{\cY}(\cF_d)$,
             \begin{equation}
             f(z)-f_{\emptyset}=\sum_{j=1}^d (S_jS_j^*f)(z)
	    = \begin{cases} {\displaystyle\sum_{j=1}^d}
             ((S^{R}_j)^*h)(z)\cdot z_j &\text{if } \; {\mathbf S} =
            {\mathbf S}^{R}, \\[3mm]
	    {\displaystyle\sum_{j=1}^{d}}
	   z_{j} \cdot ((S^{L}_j)^*h)(z) &\text{if } \; {\mathbf S} =
            {\mathbf S}^{L}.
	 \end{cases}
             \label{4.13a}
             \end{equation}

	\item
             The observability operator $\widehat{\cO}_{E, {\bf S}^{R *}}$ is
             equal to the operator $\tau$ defined in \eqref{4.9a} and
	hence is unitary.
         \end{enumerate}
             \end{proposition}

             \begin{proof}
	    If $f(z) = \sum_{v \in \cF_{d}} f_{v} z^{v}$, then
	 $$
	      ({\mathbf S}^{R * v'}f)(z) = \sum_{v \in \cF_{d}} f_{v
	      v^{\prime \top}} z^{v}, \qquad
	      ({\mathbf S}^{L* v'}f)(z) = \sum_{v \in \cF_{d}}
	      f_{v^{\prime \top} v} z^{v}
	 $$
	 and hence, in either the left or the right case, we have
	 \begin{equation}  \label{model-cO} E {\mathbf S}^{*v}f = f_{v^{\top}}.
	\end{equation}
	 Therefore,
	  $$
	\sum_{v \in \cF_{d} \colon |v| = N} \| {\mathbf S}^{*v} f
	\|^{2} = \sum_{v \in \cF_{d} \colon |v| \ge N} \| f_{v} \|^{2}
	\to 0 \quad\text{as}\quad  N \to \infty,
	$$
	and \eqref{S*ss} follows.
             Equalities \eqref{4.11a} and \eqref{4.12a} follow from \eqref{4.5},
             \eqref{4.6}, \eqref{4.5L}, \eqref{4.6L} \eqref{4.10a} and the fact
	   that $E^*$ is the inclusion map
             of $\cY$ into $H^2_{\cY}(\cF_d)$.
	Applying the operator identity
             \eqref{4.12a} to an $f\in H^2_{\cY}(\cF_d)$, we get \eqref{4.13a}.
             Finally, from \eqref{model-cO} we see that, for both the
	   left and the right case,
             $$
             \widehat{\cO}_{E, {\mathbf S}^{*}}f=\sum_{v\in\cF_d}(E{\bf
             S}^{*v}f)z^v=\sum_{v\in\cF_d}f_{v^\top}z^v=\tau f
             $$
	for all $f \in H^{2}_{\cY}(\cF_{d})$.
	That $X = I$ is the unique
	solution of the Stein equation \eqref{NC-model-Stein} is now
	a consequence of \eqref{S*ss} combined with the last part of
	Theorem \ref{T:2-1.1}.
         \end{proof}

        \begin{proof}[Proof of (4) in  Theorem \ref{T:2-1.2}:]
         Suppose that ${\mathcal M}$ is a Hilbert space contractively
         included in $H^{2}_{\cY}(\cF_{d})$ which is invariant under
         $(S^{R}_{j})^{*}$ for each $j = 1, \dots, d$ such that the
         difference-quotient inequality \eqref{dif-quot-ineq} holds.
         Set ${\mathcal X} = \tau({\mathcal M})$ with norm inherited
         from ${\mathcal M}$.  From the intertwining relations
         \eqref{Sjtau} we see that $\cX$ is invariant under the left
         backward shifts $(S^{L}_{1})^{*}, \dots, (S^{L}_{d})^{*})$.
         Define operators $A_{j} \colon {\mathcal X} \to
         {\mathcal X}$ for $j = 1, \dots, d$ and $C \colon {\mathcal X} \to \cY$
         by \eqref{modelAC}. From the difference-quotient inequality
         \eqref{dif-quot-ineq1} together with the definition of the
         ${\mathcal X}$-norm and the intertwining relations
         \eqref{Sjtau},  we have
         \begin{align*}
	   \sum_{j=1}^{d} \| (S^{L}_{j})^{*}(\tau f ) \|^{2}_{\cX} & =
	   \sum_{j=1}^{d} \| \tau (S^{R}_{J})^{*}f) \|^{2}_{\cX}
	   = \sum_{j=1}^{d} \| (S^{R}_{j})^{*}f \|^{2}_{\cM} \\
	   & \le \| f \|^{2}_{\cM} - \| f_{\emptyset}\|^{2}_{\cY}
	   = \| \tau(f) \|^{2}_{\cX} - \| (\tau
	   f)_{\emptyset}\|^{2}_{\cY}
          \end{align*}
	and hence  $H = I_{{\mathcal X}}$ satisfies the Stein inequality
         \eqref{3.4a}.  From Proposition \ref{P:4.2} we see that
         $$
         \widehat {\mathcal O}_{C, {\mathbf A}} =
         \widehat {\mathcal O}_{E, {\mathbf S}^{L*}}|_{{\mathcal X}} =
         \tau|_{{\mathcal X} = \tau ({\mathcal M})}
         $$
         and hence $\widehat {\mathcal O}_{C, {\mathbf A}} \tau|_{{\mathcal M}}
         = I_{{\mathcal M}}$. Therefore, for each $f \in {\mathcal M}$ we have
         $$
         \|f\|_{\cH(K_{C,\bA;I})}=\|\widehat{\cO}_{C, \bA}\tau
f\|_{\cH(K_{C,\bA})}=
         \|\tau f\|_{\cX}=\|f\|_{{\mathcal M}}
         $$
         and thus ${\mathcal M}=\cH(K_{C,\bA})$ isometrically.  It
	then follows from part (3a) of the theorem that in fact $\cM$
	is contractively included in $H^{}_{\cY}(\cF_{d})$.
         \end{proof}

         As explained by part (4) of Theorem \ref{T:2-1.2}, for purposes of
         study of contractively-included, backward-shift-invariant subspaces
         of $H^{2}_{\cY}(\cF_{d})$ which satisfy the
         difference-quotient-inequality \eqref{dif-quot-ineq}, without loss
         of generality we may suppose at the start that we are working with
         ${\mathcal X}'$ as the original state space ${\mathcal X}$ and with
         the solution $H$ of the Stein inequality \eqref{3.4a} to be
         normalized to $H = I_{{\mathcal X}}$.  Then certain simplifications
         occur in parts (1)-(4) of Theorem \ref{T:2-1.2} as explained in the
         next result.

         \begin{theorem}  \label{T:2-1.2'}
Suppose that $(C, {\mathbf A})$
         is a contractive pair with state space ${\mathcal X}$ and output
         space $\cY$.  Then:
         \begin{enumerate}
         \item $(C, {\mathbf A})$ is output-stable and the intertwining relation
         \eqref{4.8a} holds. Hence
         $\operatorname{Ran}  \widehat \cO_{C,{\mathbf A}}$ is invariant under
         the backward shifts $(S^{R}_{j})^{*}$ for $j = 1, \dots, d$.

         \item The observability operator $\widehat {\mathcal
         O}_{C,{\mathbf A}}$ is a
         contraction from ${\mathcal X}$ into $H^{2}_{\cY}(\cF_{d})$.
         Moreover $\widehat {\mathcal O}_{C, {\mathbf A}}$ is isometric if
         and only if $(C, {\mathbf A})$ is an isometric pair and ${\mathbf
         A}$ is strongly stable.

         \item If the linear manifold ${\mathcal M}:= \operatorname{Ran}
         \widehat {\mathcal O}_{C,{\mathbf A}}$ is given the lifted norm
         \begin{equation}  \label{lifted'}
           \| \widehat {\mathcal O}_{C, {\mathbf A}}x\|_{{\mathcal M}} =
           \| Q x \|_{{\mathcal X}}
         \end{equation}
         where $Q$ is the orthogonal projection of ${\mathcal X}$ onto
         $(\operatorname{Ker} \widehat{\mathcal O}_{C, {\mathbf
A}})^{\perp}$, then
         $\widehat{\mathcal O}_{C, {\mathbf A}}$ is a coisometry of
         ${\mathcal X}$ onto ${\mathcal M}$. Moreover,
         ${\mathcal M}$ is contained contractively in $H^{2}_{\cY}(\cF_{d})$
         and is isometrically equal to the formal noncommutative reproducing
         kernel Hilbert space ${\mathcal H}(K_{C,{\mathbf A}})$ with
         reproducing kernel $K_{C,{\mathbf A}}(z,w)$ given by
         $$
           K_{C, {\mathbf A}}(z,w) = C (I - Z(z) A)^{-1} (I -
           A^{*}Z(w)^{*})^{-1} C^{*}.
         $$

         \item If ${\mathcal O}_{C, {\mathbf A}}$ is given the lifted norm
         $\| \cdot \|_{{\mathcal H}(K_{C,{\mathbf A}})}$ as in
         \eqref{lifted'}, then the difference-quotient inequality
         \begin{equation}  \label{dif-quot-ineq'}
             \sum_{j=1}^{d} \| (S^{R}_{j})^{*} f\|^{2}_{{\mathcal
H}(K_{C,{\mathbf
             A}})} \le \| f\|^{2}_{{\mathcal H}(K_{C,{\mathbf A}})} - \|
             f_{\emptyset}\|^{2}_{\cY}
         \end{equation}
         holds for all $f \in \cH(K_{C,\bA})$.
	Moreover, \eqref{dif-quot-ineq'} holds with equality
         if and only the orthogonal projection $Q$ of ${\mathcal X}$ onto
         $(\operatorname{Ker} {\mathcal O}_{C, {\mathbf A}})^{\perp}$
         satisfies the Stein equation
         \begin{equation}  \label{Q-Stein}
           Q - \sum_{j=1}^{d}  A_{j}^{*} Q A_{j} = C^{*}C.
         \end{equation}
         In particular, if $(C, {\mathbf A})$ is observable, then
\eqref{dif-quot-ineq'} holds with equality if and only if $(C, {\mathbf
A})$ is an isometric pair.
         \end{enumerate}
         \end{theorem}

         \begin{proof}  Statements (1)--(3) and all but the last part of
         statement (4) are direct specializations to the case $H =
         I_{{\mathcal X}}$ of the corresponding results in Theorem 
\ref{T:2-1.2}.
         It remains only to analyze the conditions for equality in
         \eqref{dif-quot-ineq'}.

         From the intertwining relation \eqref{4.8a}, we see that equality
         in \eqref{dif-quot-ineq'} for a generic element $f = \widehat{\mathcal
         O}_{C,{\mathbf A}}x \in {\mathcal H}(K_{C,{\mathbf A}})$ means that
         $$
           \sum_{j=1}^{d} \| \widehat{\mathcal O}_{C, {\mathbf A}} A_{j}x
           \|^{2}_{{\mathcal H}(K_{C, {\mathbf A}})} = \| \widehat{\mathcal
          O}_{C, {\mathbf A}} x \|^{2}_{{\mathcal H}(K_{C,{\mathbf A}})} - \|
           C x\|^{2}_{\cY}
         $$
         for all $x \in {\mathcal X}$.  Using the definition \eqref{lifted'}
         of the ${\mathcal H}(K_{C,{\mathbf A}})$-norm, we rewrite this last
         equality as
         $$
           \sum_{j=1}^{d} \| Q A_{j} x \|^{2}_{{\mathcal X}} = \| Q x
           \|^{2}_{{\mathcal X}} - \| C x \|^{2}_{\cY}.
         $$
         This holding for all $x \in {\mathcal X}$ is finally equivalent to
         the Stein equation \eqref{Q-Stein}.    \end{proof}

         \begin{remark}  \label{R:2-1.2} {\em In Theorems \ref{T:2-1.2}
	   and \ref{T:2-1.2'} we could equally well have interchanged
	   the roles of left versus right.  For a given output pair
	   $(C, \bA)$, define the associated {\em left observability
	   operator} $\cO^{L}_{C, \bA} \colon \cX \to
	   H^{2}_{\cY}(\cF_{d})$ by
	   $$ \cO^{L}_{C, \bA}x = \sum_{v \in \cF_{d}} C
	   \bA^{v^{\top}} z^{v}.
	   $$
	   Then the linear manifold $\operatorname{Ran}\, \cO^{L}_{C,
	   \bA}$ is invariant under the left backward shifts
	   $(S^{L}_{1})^{*}, \dots, (S^{L}_{d})^{*})$ as verified by
	   the intertwining relation
	   $$
	  (S^{L}_{j})^{*} \cO^{L}_{C, \bA} = \cO^{L}_{C, \bA} A_{j}.
	  $$
	  We leave the precise statements and proofs to the
	  interested reader.}
\end{remark}

          The characterization \eqref{Q-Stein} of the difference-quotient
          inequality holding with equality for a space ${\mathcal H}(K_{C,
          {\mathbf A}})$ in Theorem \ref{T:2-1.2'} can be made more explicit
          as follows.
          \begin{proposition}  \label{P:2-1.2'}
              Suppose that $(C,{\mathbf A})$ is a contractive pair as in
              Theorem \ref{T:2-1.2'} and let $Q$ be the orthogonal
              projection onto $(\operatorname{Ker} \widehat{\mathcal O}_{C,
            {\mathbf A}})^{\perp}$.  Then $Q$ satisfies the Stein inequality
            \begin{equation}  \label{Q-Stein-ineq}
	 Q - A_{1}^{*}Q A_{1} - \cdots A_{d}^{*}Q A_{d} \ge C^{*}C
          \end{equation}
          and we have the inequalities
            \begin{equation} \label{Q-ineq}
	 {\mathcal G}_{C, {\mathbf A}} \le Q \le I_{{\mathcal X}}.
          \end{equation}
              If we write $C, A_{j},Q$ in $2 \times 2$-block
              matrix form with respect to the decomposition ${\mathcal X} =
              \operatorname{Ker} \widehat{\mathcal O}_{C, {\mathbf A}} \oplus
              (\operatorname{Ker} \widehat{\mathcal O}_{C, {\mathbf
A}})^{\perp}$ as
         \begin{equation}  \label{matrix-decom}
          C = \begin{bmatrix} 0 & C^{0} \end{bmatrix}, \qquad
            A_{j}^{0} = \begin{bmatrix} A_{j1} & A_{j2} \\ 0 & A_{j}^{0}
            \end{bmatrix}, \qquad
            Q = \begin{bmatrix} 0 & 0 \\ 0 & I \end{bmatrix}
        \end{equation}
        for $j = 1, \dots, d$,  then $Q$ satisfies the Stein equation
        \eqref{Q-Stein} if and only if the pair $(C^{0},
        {\mathbf A}_{j}^{0})$ is an isometric pair, in which case we also
        have that $A_{j2} = 0$ (so
        $(\operatorname{Ker} \widehat{\mathcal O}_{C, {\mathbf A}})^{\perp}$ is
        invariant for $A_{j}$) for $j = 1, \dots, d$.
        \end{proposition}
        \begin{proof}
            First note that $\operatorname{Ker} \widehat{\mathcal O}_{C,
            {\mathbf A}}$ is invariant for each $A_{j}$ and that
            $ \operatorname{Ker} \widehat{\mathcal
            O}_{C, {\mathbf A}} \subset  \operatorname{Ker} C$.
            Therefore the matrix decompositions of $C,A_{j},Q$ with respect
            to the decomposition ${\mathcal X} =
            \operatorname{Ker} \widehat{\mathcal O}_{C, {\mathbf A}} \oplus
            (\operatorname{Ker} \widehat{\mathcal O}_{C, {\mathbf
            A}})^{\perp}$ have the form as given in \eqref{matrix-decom}.  Next
            note that the contractive property of the pair $(C, {\mathbf A})$
            means that
            \begin{equation}  \label{work1}
             \begin{bmatrix} 0 & 0 \\ 0 & C^{0*}C^{0} \end{bmatrix}
	  + \sum_{j=1}^{d} \begin{bmatrix} A_{j1}^{*} A_{j1} &
	  A_{j1}^{*}A_{j2} \\ A_{j2}^{*}A_{j1} & A_{j2}^{*}A_{j2} +
	  A_{j}^{0*}A_{j}^{0*} \end{bmatrix} \le \begin{bmatrix} I & 0 \\
	  0 & I \end{bmatrix}.
          \end{equation}
          On the other hand, the Stein inequality \eqref{Q-Stein-ineq} works
          out to be
          \begin{equation}  \label{work2}
          \begin{bmatrix} 0 & 0 \\ 0 & C^{0*}C^{0} \end{bmatrix}
	     + \sum_{j=1}^{d} \begin{bmatrix} A_{j1}^{*} A_{j1} &
	     A_{j1}^{*}A_{j2} \\ A_{j2}^{*}A_{j1} &
	     A_{j}^{0*}A_{j}^{0*} \end{bmatrix} \le \begin{bmatrix} I & 0 \\
	     0 & I \end{bmatrix}.
         \end{equation}
         As the left hand side of \eqref{work2} is dominated by the left
         hand side of \eqref{work1}, it is clear that \eqref{work2} follows
         from \eqref{work1}, and hence \eqref{Q-Stein-ineq} holds as asserted.
         Since ${\mathcal G}_{C, {\mathbf A}}$ is the minimal
         positive semidefinite solution of the Stein inequality \eqref{3.4a}
         (by part (2) of Theorem \ref{T:2-1.1})
         and we now know that $Q$ is one such solution, it follows that
         ${\mathcal G}_{C, {\mathbf A}} \le Q$. As $Q$ is an orthogonal
         projection on ${\mathcal X}$, we also have $Q \le I_{{\mathcal X}}$
         and \eqref{Q-ineq} now follows.

          From the $(2,2)$ entry of \eqref{work1}, we read off
          \begin{equation}  \label{Stein-1}
              C^{0*}C^{0} + \sum_{j=1}^{d} A_{j2}^{*}A_{j2} + \sum_{j=1}^{d}
              A_{j}^{0*}A_{j}^{0} \le I_{(\operatorname{Ker} \widehat
              \cO_{C, \bA})^{\perp}}.
          \end{equation}
          In particular
          $$
          C^{0*}C^{0} + \sum_{j=1}^{d} A_{j}^{0*}A_{j}^{0} \le
          I_{(\operatorname{Ker} \widehat \cO_{C,\bA})^{\perp}},
          $$
          i.e., $Q$ satisfies \eqref{Q-Stein-ineq}.
          On the other hand, the validity of \eqref{Q-Stein} reduces to
         $$
         \begin{bmatrix} 0 & 0 \\ 0 & C^{0*}C^{0} \end{bmatrix} + \sum_{j=1}^{d}
         \begin{bmatrix} 0 & 0 \\ 0 & A_{j}^{0*} A_{j}^{0} \end{bmatrix} =
         \begin{bmatrix} 0 & 0 \\ 0 & I \end{bmatrix},
          $$
          or simply to
          \begin{equation}  \label{Stein-2}
              C^{0*} C^{0} + \sum_{j=1}^{d} A_{j}^{0*} A_{j}^{0} =
              I_{(\operatorname{Ker} \widehat \cO_{C,\bA})^{\perp}}.
          \end{equation}
          Thus, the validity of \eqref{Q-Stein} is equivalent to
          $(C^{0},{\mathbf A}^{0})$ being an isometric pair, in which case we
          also have that $A_{j2}=0$.
\end{proof}

Finally, we have the following uniqueness result.

\begin{theorem}  \label{T:unique}
         Suppose that $(C, \bA)$ and $(\widetilde C, \widetilde \bA)$ are two
         output-stable, observable pairs realizing the same positive kernel
       \begin{eqnarray}
       K_{C, \bA}(z,w)&:=& C(I - Z(z) A)^{-1}(I - A^{*}Z(w)^{*})^{-1}C^{*}
       \label{sameker}\\
       &=& \widetilde C(I - Z(z) \widetilde A)^{-1}(I - \widetilde
       A^{*}Z(w)^{*})^{-1}\widetilde C^{*}=: K_{\widetilde C,\widetilde
\bA}(z,w).
       \nonumber
       \end{eqnarray}
       Then $(C,\bA)$ and $(\widetilde C,\widetilde \bA)$ are
       {\em unitarily equivalent}, i.e.,
       there is a unitary operator $U \colon \cX \to \widetilde \cX$ such that
       $$
         C = \widetilde C U \quad\text{and}\quad A_{j} = U^{-1} \widetilde
A_{j}U \; \; \text{ for } \; j = 1,
         \dots, d.
       $$
       \end{theorem}

       \begin{proof}
           For any two words $\alpha, \beta \in \cF_{d}$, equating
           coefficients of $z^{\alpha} w^{\beta^{\top}}$ in \eqref{sameker}
           gives
         $$
         C \bA^{\alpha} \bA^{*\beta}C^{*} = C' \widetilde \bA^{ \alpha}
         \widetilde \bA^{* \beta} C^{\prime *}.
         $$
         Hence the operator $U$ defined by
       \begin{equation}  \label{intertwine}
         U \colon \bA^{* \beta} C^{*} y \mapsto \widetilde \bA^{ * \beta}
         \widetilde C^{*}
         y
       \end{equation}
       extends by linearity and continuity to define an isometry from
       $$
       {\mathcal D}_{U} = \operatorname{ \overline{span}} \{ \bA^{* \beta}
       C^{*}y \colon \beta \in \cF_{d},\ y \in \cY \}
       $$
       onto
       $$
       {\mathcal R}_{U} = \operatorname{\overline{span}} \{ \widetilde\bA^{ *
       \beta} \widetilde C^{*} y \colon \beta \in \cF_{d},\ y \in \cY\}.
       $$
       The observability assumption implies that ${\mathcal D}_{U} = \cX$
       and ${\mathcal R}_{U} = \widetilde {\mathcal X}$; hence $U \colon \cX
       \to \widetilde \cX$
       is unitary.  From  \eqref{intertwine} it is easily seen that
       $$
         U C^{*} = \widetilde C^{*}\quad \text{and}\quad U A_{j}^{*} =
\widetilde A_{j}^{*}U \quad\text{for} \; \;  j = 1, \dots, d.
       $$
       Since $U$ is unitary we then get
       $$
        \widetilde C U = C \quad \text{and}\quad\widetilde A_{j}U = U A_{j}
\quad\text{for} \; \; j = 1, \dots, d
       $$
       and we conclude that $(C, \bA)$ and $( \widetilde{C},
\widetilde{\bA})$ are unitarily equivalent as desired.

       \end{proof}

       \subsection{Applications of observability operators: the
       noncommutative setting}  \label{S:NC-appl}

       As an application we give a proof of the Beurling-Lax theorem for
       the Fock-space setting originally given by Popescu \cite{PopescuBL}.
       We shall in fact prove a more general version of the
       Beurling-Lax-Halmos theorem for contractively-included (rather than
       isometrically included) subspaces of $H^{2}_{\cY}(\cF_{d})$ due in
       the classical setting to de Branges (see \cite{deBR}).
       Our proof is similar to that in \cite{PopescuBL} but highlights more
       explicitly the role of an associated observability operator.
       For this purpose we say that a formal power series $\theta(z) =
       \sum_{v \in \cF_{d}} \theta_{v} z^{v} \in {\mathcal L}(\cU, \cY)\langle
       \langle z \rangle \rangle$ is a {\em contractive multiplier}, also
       written as $\theta$ is in the $d$-variable, noncommutative Schur-class
       $\cS_{nc,d}(\cU, \cY)$,  if the operator $M_{\theta}$ of multiplication
       by $\theta$
       $$
        M_{\theta} \colon f(z) \mapsto \theta(z) \cdot f(z)
       $$
       defines a bounded linear operator from $H^{2}_{\cU}(\cF_{d})$ to
       $H^{2}_{\cY}(\cF_{d})$ with operator norm at most $1$.  Such a formal
       power series $\theta(z)$ is said to be {\em inner} if moreover
       the operator
       $M_{\theta}$ from $H^{2}_{\cU}(\cF_{d})$ to $H^{2}_{\cY}(\cF_{d})$ is
       an isometry\footnote{We prefer to define {\em inner} to be
       isometric rather than partially isometric as in \cite{PopescuBL}.}.

       \begin{theorem} \label{T:NC-BL}
	 \begin{enumerate}
	\item
	A Hilbert space $\cM$ is such that
	 \begin{enumerate}
	     \item $\cM$ is contractively included in
	     $H^{2}_{\cY}(\cF_{d})$,
	     \item $\cM$ is invariant under the right shift operators
	     $S^{R}_{1}, \dots, S^{R}_{d}$:
	     $$ S^{R}_{j} \cM \subset \cM,
	     $$
        \item the $d$-tuple
        $$ {\mathbf S}^{R}_{\cM} = (S^{R}_{\cM,1}, \dots, S^{R}_{\cM,d})
\quad \text{where} \; \;
        S^{R}_{\cM,j}:= S^{R}_{j}|_{\cM} \; \text{ for } \; j = 1, \dots, d
        $$
        is a row contraction
        $$
         S^{R}_{\cM,1}(S^{R}_{\cM,1})^{*} + \cdots +
S^{R}_{\cM,d}(S^{R}_{\cM,d})^{*} \le
         I_{\cM},
        $$
        and
        \item $({\mathbf S}^{R}_{\cM})^{*}$ is strongly stable, i.e.
        $$ \lim_{n\to \infty} \sum_{v \in \cF_{d}\colon |v| = n} \|
({\mathbf S}^{R}_{\cM})^{*v}
         f \|^{2}_{\cM} \to 0 \quad\text{for all}\quad f \in \cM,
        $$
        \end{enumerate}
         if and only if there is a coefficient Hilbert space $\cU$ and a
         contractive multiplier $\theta \in \cS_{nc,d}(\cU, \cY)$ so that
         $$
	 \cM = \theta \cdot H^{2}_{\cU}(\cF_{d})
	$$
	with lifted norm
	\begin{equation}  \label{M-norm}
	\| \theta \cdot f \|_{\cM} = \| Q f \|_{H^{2}_{\cU}(\cF_{d})}
	\end{equation}
	where $Q$ is the orthogonal projection onto $(\operatorname{Ker}
         \, M_{\theta})^\perp$.

	\item The subspace $\cM$ in part (1) above is isometrically included in
	$H^{2}_{\cY}(\cF_{d})$ if and only if the associated contractive
	multiplier $\theta$ is inner.
	 \end{enumerate}
        \end{theorem}

        \begin{proof}  We first verify sufficiency in statement (1).
	   Suppose that $\cM$ has
	  the form $\cM = \theta \cdot H^{2}_{\cY}(\cF_{d})$ for a
	  contractive multiplier $\theta$ with $\cM$-norm given by
	  \eqref{M-norm}.  From the fact that $\|M_{\theta}\| \le 1$ it is
	  easily verified that $\|\theta \cdot f\|_{H^{2}_{\cY}(
	  \cF_{d})} \le \| \theta \cdot f \|_{\cM}$, i.e., (a) holds.
	  From the intertwining property $S^{R}_{j} M_{\theta} = M_{\theta}
	  S^{R}_{j}$ (note that $S^{R}_{j}$ is multiplication  by
$z_{j}$ on the
	  right while $M_{\theta}$ is multiplication by $\theta$ on the
	  left), property (b) follows.

	  If $Q$ is the orthogonal projection onto $(\operatorname{Ker}
	  M_{\theta})^{\perp} \subset H^{2}_{\cU}(\cF_{d})$, then the
	  intertwining property $S^{R}_{j} M_{\theta} = M_{\theta}
	  S^{R}_{j}$
	  implies that
	  \begin{equation}  \label{Q-Sj}
         Q S^{R}_{j} = Q S^{R}_{j}Q\quad \text{and}\quad (S^{R}_{j})^{*} Q
         = Q (S^{R}_{j})^{*}Q \quad\text{for}\quad j = 1, \dots, d.
         \end{equation}
         Thus
         \begin{align*}
	 \left\| \begin{bmatrix} S^{R}_{\cM,1} & \cdots & S^{R}_{\cM, d}
         \end{bmatrix} \begin{bmatrix} \theta \cdot f_{1} \\ \vdots \\
         \theta \cdot f_{d} \end{bmatrix} \right\|^{2}_{\cM} & =
         \left\| \theta \begin{bmatrix} S^{R}_{1} & \cdots & S^{R}_{d}
\end{bmatrix}
         \begin{bmatrix} f_{1} \\ \vdots \\ f_{d} \end{bmatrix}
	   \right\|^{2}_{\cM} \\
	  &= \left\| Q \begin{bmatrix} S^{R}_{1} & \cdots & S^{R}_{d}
         \end{bmatrix} \begin{bmatrix} Q f_{1} \\ \vdots \\ Q f_{d}
       \end{bmatrix} \right\|^{2}_{H^{2}_{\cU}(\cF_{d})}\\
       &\le \left\| \begin{bmatrix} Qf_{1} \\ \vdots \\ Qf_{d}\end{bmatrix}
       \right\|^{2}_{H^{2}_{\cU}(\cF_{d})^{d}} = \left\| \begin{bmatrix}
       \theta \cdot f_{1} \\ \vdots \\ \theta \cdot f_{d} \end{bmatrix}
       \right\|^{2}_{\cM^{d}}
       \end{align*}
       and property (c) follows.  Finally, a short computation shows that
       $$
        (S^{R}_{\cM,j})^{*} \colon \theta \cdot f \mapsto \theta \cdot
        S^{R*}_{j}
        Q f, \qquad ({\mathbf S}^{R}_{\cM})^{*v} \colon \theta \cdot f
        \mapsto\theta \cdot {\mathbf S}^{R *v} Q f
       $$
       and hence
       $$
       \sum_{v \in \cF_{d} \colon |v| = n} \| ({\mathbf S}^{R}_{\cM})^{*v}
\theta
       \cdot f \|^{2}_{\cM} =
       \sum_{v \in \cF_{d} \colon |v| = n} \| {\mathbf S}^{R*v} Q f
       \|^{2}_{H^{2}_{\cU}(\cF_{d})} \to 0
       $$
       as $n \to \infty$, and property (d) follows as well.  Moreover, if
$\theta$ is inner and
       $\cM = \theta \cdot H^{2}_{\cU}(\cF_{d})$ with the lifted norm
       \eqref{M-norm}, it is clear that $\cM$ is contained in
       $H^{2}_{\cY}(\cF_{d})$ isometrically.  This completes the proof of
       sufficiency in Theorem \ref{T:NC-BL}.

       Suppose now that the Hilbert space $\cM$ satisfies conditions
       (a), (b), (c), (d) in statement (1) of Theorem \ref{T:NC-BL}.
       Define a $d$-tuple of operators ${\mathbf A} = (A_{1}, \dots, A_{d})$
       on ${\mathcal M}$ by
       $$
        A_{j} = (S^{R}_{\cM,j})^{*} \quad\text{for} \quad j = 1, \dots, d,
       $$
       where we use hypothesis (b) to set
       $S^R_{\cM,j}:= S^R_{j}|_{\cM}$ for $j = 1, \dots, d$,
       and choose the coefficient Hilbert space $\cU$ so that
       $$
        \operatorname{dim} \cU = \operatorname{rank} (I - A_{1}^{*}A_{1} -
        \cdots - A_{d}^{*}A_{d}).
       $$
       By hypothesis (c) we may then choose the operator
       $C \colon \cM \to \cU$ so that
       $$
         C^{*}C = I - A_{1}^{*}A_{1} - \cdots A_{d}^{*}A_{d}.
       $$
       Then $(C, \bA)$ is an isometric pair and, by hypothesis (d),
       $\bA^{*}$ is strongly stable.
       Thus by part (2) of Proposition \ref{P:2-1.1} it follows that the
       observability operator
       $$
         \widehat \cO_{C, \bA} \colon f \mapsto ( I - Z(z) A)^{-1} f
       $$
       is an isometry from $\cM$ into $H^{2}_{\cU}(\cF_{d})$.  As observed
       for the general case in part (1) of Theorem \ref{T:2-1.2}, we have
       the intertwining condition
       $$ (S^{R}_{j})^{*} \widehat \cO_{C,\bA} =
       \widehat \cO_{C, \bA} (S^{R}_{\cM,j})^{*}.
       $$
       Taking adjoints then gives
       \begin{equation}  \label{intertwine1}
         (\widehat \cO_{C, \bA})^{*} S^{R}_{j} =
         S^R_{\cM,j} (\widehat \cO_{C, \bA})^{*}.
       \end{equation}
       Let us set
       $$
       \Theta = \iota \circ (\widehat \cO_{C, \bA})^{*} \colon
H^{2}_{\cU}(\cF_{d})
       \to H^{2}_{\cY}(\cF_{d})
       $$
       where $\iota \colon \cM \to H^{2}_{\cY}(\cF_{d})$ is the inclusion map.
       From hypothesis (a) that $\| \iota \| \le 1$, we see that $\| \Theta
       \| \le 1$.  From the intertwining relation \eqref{intertwine1}
       (together with hypothesis (b)) it
       follows that
       $$
         \Theta S^{R}_{j} = S^{R}_{j} \Theta
       $$
       and it follows (see e.g.~\cite{Popescu-multi}) that $\Theta$ is a
       multiplication operator, i.e., there is a contractive multiplier
       $\theta \in \cS_{nc, d}(\cU, \cY)$ so that $\Theta = M_{\theta}$.
       From the fact that $\widehat \cO_{C, \bA} \colon \cM \to
H^{2}_{\cU}(\cF_{d})$
       is an isometry, it follows that $\operatorname{Ran} (\widehat \cO_{C,
       \bA})^{*} = \cM$ and also that $\cM = \theta \cdot
       H^{2}_{\cU}(\cF_{d})$ with $\cM$-norm given by
       \eqref{M-norm}.  This completes the proof of necessity in
       statement (1) of Theorem \ref{T:NC-BL} for the
       general case.

       We now consider statement (2).
       In case $\cM$ is isometrically included in
       $H^{2}_{\cY}(\cF_{d})$,
       for any $f \in
	    H^{2}_{\cU}(\cF_{d})$ we have
	    $$
	     \| (\widehat \cO_{C,{\mathbf A}})^{*} S^{R}_{j} f
	     \|_{\cM} = \| S^{R}_{j}
	    (\widehat \cO_{C,{\mathbf A}})^{*} f \|_{H^{2}_{\cY}(\cF_{d})}
	    = \| (\widehat \cO_{C,{\mathbf A}})^{*} f
	    \|_{H^{2}_{\cU}(\cF_{d})}
	    $$
             for $j = 1, \dots, d$,	    since $S^{R}_{j}$ is isometric
             on  $H^{2}_{\cY}(\cF_{d})$.  Since, as was observed above,
	    $\widehat \cO_{C, {\mathbf A}}$ is isometric, it follows that
	    $$ \| P_{\operatorname{Ran} \widehat  \cO_{C, {\mathbf
             A}}} S_{j} f \| = \| f \| \quad\text{for all} \quad f \in
             \operatorname{Ran}\widehat  \cO_{C, {\mathbf A}}
	    $$
	    and hence $\operatorname{Ran} \widehat \cO_{C, {\mathbf A}}$ is
	    invariant under $S^{R}_{j}$ for $j = 1, \dots, d$.  As
	    $\operatorname{Ran} \widehat \cO_{C, {\mathbf A}}$ is also
	    invariant under $(S^{R}_{j})^{*}$ for each $j$ by
\eqref{intertwine1},
	    we conclude that $\operatorname{Ran} \widehat \cO_{C,{\mathbf
	    A}}$ is reducing for ${\mathbf S}^{R}$. Since
$\operatorname{Ran} C$ is
	    dense in $\cU$ by construction, we are now able to conclude that
	    $\operatorname{Ran} \widehat \cO_{C, {\mathbf A}}$ is all of
	    $H^{2}_{\cU}(\cF_{d})$ and hence $\widehat \cO_{C, \bA}
\colon \cM \to
	    H^{2}_{\cU}(\cF_{d})$ is actually unitary.  It then follows
	    finally that $\Theta = \iota \circ (\widehat \cO_{C, \bA})^{*}$ is
	    isometric and hence $\theta$ is inner as asserted. This
	    completes the proof of Theorem \ref{T:NC-BL}.
       \end{proof}

        A second application of these ideas is to operator model theory.
          For this application we are given only an operator-tuple
          ${\mathbf T} = (T_{1}, \dots, T_{d}) \in \cL(\cH)^{d}$ which is a row
          contraction, so $I - T_{1}T_{1}^{*} - \cdots - T_{d}T_{d}^{*}
          \ge 0$.  Set
          \begin{equation}  \label{T*defect}
          D_{{\mathbf T}^{*}}: = (I - T_{1} T_{1}^{*} - \cdots - T_{d}
          T_{d}^{*})^{1/2}\quad \text{and} \quad\cY : = \operatorname{Ran}
D_{{\mathbf
          T}^{*}}.
          \end{equation}
          We apply the ideas of the previous sections concerning the
          general pair $(C, \bA)$ to a pair of the special form
          $(D_{{\mathbf T}^{*}}, {\mathbf T}^{*})$.
          For simplicity we assume in addition that ${\mathbf T}^{*}$ is
          asymptotically stable, i.e.,
          $$
          \lim_{N \to \infty} \sum_{v \in \cF_{d} \colon |v| = N} \|
          {\mathbf T}^{*v}x \|^{2} = 0\quad  \text{for all}\quad  x\in\cH.
          $$
          Then we have the following dilation result.

          \begin{theorem} \label{T:NCop-model}  Suppose that ${\mathbf T}
	 = (T_{1}, \dots, T_{d})$ is a row contraction with ${\mathbf
	 T}^{*}$ asymptotically stable as above and define the defect
	 operator $D_{{\mathbf T}^{*}}$ and the coefficient space
	 $\cY$ as in \eqref{T*defect}.   Then there is a subspace
	 $\cM \subset H^{2}_{\cY}(\cF_{d})$ invariant for the
	 backward shift operator-tuple ${\mathbf S}^{R*}$ on
	 $H^{2}_{\cY}(\cF_{d})$ so that ${\mathbf T}$ is unitarily
	 equivalent to $P_{\cM}{\mathbf S}^{R}|_{\cM}$.  In particular,
	 ${\mathbf T}$ has a row-shift dilation unitarily equivalent
	 to ${\mathbf S}^{R}$ on $H^{2}_{\cY}(\cF_{d})$.
       \end{theorem}

       \begin{proof}  By the same arguments as in the proof of Theorem
           \ref{T:NC-BL}, we see that
           $$
           \widehat \cO_{D_{{\mathbf T}^{*}}, {\mathbf
           T}^{*}} \colon \cH \to H^{2}_{\cY}(\cF_{d})
           $$
           is isometric and satisfies the intertwining
           $$
           (S^{R}_{j})^{*} \widehat \cO_{D_{{\mathbf T}^{*}}, {\mathbf T}^{*}} =
           \widehat \cO_{D_{{\mathbf T}^{*}},{\mathbf T}^{*}} T_{j}^{*}.
           $$
           If we then set
           $$
             \cM = \operatorname{Ran} \widehat \cO_{D_{{\mathbf 
T}^{*}}, {\mathbf
	T}^{*}},
          $$
          then $\widehat \cO_{D_{{\mathbf T}^{*}}, {\mathbf T}*}$ implements the
          unitary equivalence between ${\mathbf T}$ and $P_{\cM} {\mathbf
          S}|_{\cM}$ as wanted.
           \end{proof}

           \begin{remark} \label{R:NC-vN}
	  {\em In the classical case $d=1$, the procedure for
	  constructing the unitary dilation of a contraction operator
	  via the observability operator as in the proof
	  of Theorem \ref{T:NCop-model} corresponds to the
	  construction of Douglas (see \cite{Douglas})
	  (see also \cite[Section I.10.1]{NF}) which is an
	  alternative to the more popular Sch\"affer-matrix
	  construction of the unitary dilation (see \cite[Section I.5]{NF}).
	  Popescu (see \cite{PopescuNF0, PopescuNF1}) used an analogue
	  of the Sch\"affer-matrix construction to construct the
	  row-unitary dilation of a row-contraction operator-tuple.
	  From the existence of this dilation, he went on to
	  verify a von Neumann inequality (see \cite{Popescu-vN}):
	  $$
	  \| p(T_{1}, \dots, T_{d}) \| \le \| p(S_{1}, \dots, S_{d})\|
	  $$
	  for any polynomial $p \in {\mathbb C}\langle z \rangle$ in
	  the noncommuting variable $z = (z_{1}, \dots, z_{d})$.  He
	  returned to this topic in \cite{popjfa} to give another
	  proof of the von Neumann inequality (actually a more
	  general version involving nonanalytic polynomials) based on
	  the Poisson transform: for ${\mathbf T}$ a {\em strict}
	  row-contraction (one can reduce the general case of a
	  row-contraction to the case of a strict row-contraction via
	  a limiting procedure), one defines the Poisson transform
	  $P({\mathbf T}) \colon \cL(H^{2}_{\cY}(\cF_{d}), \cH) \to
	  \cL(\cH)$ by
	  \begin{equation}  \label{Poissontrans}
	  P({\mathbf T})[X] = (\widehat \cO_{D_{{\mathbf T}^{*}}, {\mathbf
	  T}^{*}})^{*} X \widehat \cO_{D_{{\mathbf T}^{*}}, {\mathbf T}^{*}}.
	  \end{equation}
	  It is argued in \cite{popjfa} (as well as in
	  \cite{Chalendar} in the context of the classical case) that
	  this is an elementary (i.e., dilation-free) proof of the
	  von Neumann inequality.  Indeed, as argued in
	  \cite{Chalendar}, this proof of the von
	  Neumann inequality goes back to the paper of Heinz
	  \cite{Heinz}.  However, we would argue that the dilation is
	  very near the surface in this proof as well, since the
	  Poisson kernel, i.e., the observability operator
	  $\widehat \cO_{D_{{\mathbf T}^{*}}, {\mathbf T}^{*}}$, provides
	  the factorization of the Poisson transform
	  \eqref{Poissontrans} and is also the operator embedding the
	  state space $\cH$ into the dilation space
	  $H^{2}_{\cY}(\cF_{d})$ in the Douglas approach to dilation
	  theory.}
	  \end{remark}

      \section{The commutative-variable Arveson-space setting}
          \label{S:Arveson}

          \subsection{Output stability and Stein equations: the
          commutative-variable case}
          \label{S:C-Stein}

          To introduce the commutative multidimensional counterpart of the
          Hardy space $H^{2}({\mathbb D})$, we recall
          standard multivariable notations: for a multi-integer
          $$
          \bn =(n_{1},\ldots,n_{d}) \in \bbZ_+^d
          $$
          and a point
          $\blam=(\lambda_1,\ldots,\lambda_d)\in\C^d$,
          we set $|\bn| = n_{1}+n_{2}+\ldots +n_{d}$,
          $\bn!= n_{1}!n_{2}!\ldots n_{d}!$ and
\begin{equation}
          \blam^\bn = \lambda_{1}^{n_{1}}\lambda_{2}^{n_{2}}\ldots
\lambda_{d}^{n_{d}}.
\label{mnot}
\end{equation}
          The space can be derived from the full Fock space by
          ``letting the variables commute''.
          For this purpose we introduce the abelianization map
             $\ba \colon \, \cF_{d} \to{\mathbb Z}^{d}_{+}$
          given by
          $$
          \ba(i_{N} \cdots i_{1}) = (n_{1}, \dots, n_{d})\quad\text{where}
          \quad n_{k} = \#\{\ell\colon \; i_{\ell} = k\} \; \text{ for } \; k=1,
          \dots, d.
          $$
          A key combinatorial fact is that
          \begin{equation}
          \#\ba^{-1}(\bn) = \frac{|\bn|!}{\bn!}.
          \label{2.19a}
          \end{equation}
          We then consider the symmetric Fock space
          $\ell^{2}_{\cY}({\mathcal S}\cF_{d})$ equal to the subspace of
          $\ell^{2}_{\cY}(\cF_{d})$ spanned by the elements
          $\chi_{\bn} y$ ($n \in \bbZ^{d}_{+}$ and $y \in \cY$) where
          $\chi_{\bn}$ is given by
          $$
              \chi_{\bn} = \sum_{v \colon \ba(v) = \bn} \chi_{v}.
          $$
          Note that
          $$ \| \chi_{\bn}  \|^{2}_{\ell^{2}_{\bbC}(\cF_{d})} = \sum_{v \in
          \cF_{d} \colon \ba(v) = n} 1 = \frac{|\bn|!}{\bn !}
          $$
          and hence, if ${\mathcal B}$ is an orthonormal basis for ${\mathcal
          Y}$, then an orthonormal basis for $\ell^{2}_{\cY}({\mathcal S}
          \cF_{d})$ is the set
          $$
          \left\{ \sqrt{\frac{ \bn !}{|\bn|!}}\chi_{\bn} y \colon \bn \in
          \bbZ^{d}_{+}, \, y \in
          {\mathcal B}\right \}.
          $$
          It is then natural to identify $\ell^{2}_{\cY}({\mathcal S} \cF_{d})$
          with the weighted sequence space $\ell^{2}_{w, \cY}(\bbZ^{d}_{+})$
          consisting of all $\cY$-valued $\bbZ^{d}_{+}$-indexed sequences
          $ \{f_{\bn}\}_{\bn \in \bbZ^{d}_{+}} $ for which the norm given by
$$
              \| \{f_{\bn}\}_{\bn \in \bbZ^{d}_{+}} \|_{\ell^{2}_{w,
              \cY}(\bbZ^{d}_{+})} = \sum_{\bn \in \bbZ^{d}_{+}} w(\bn)
\|f_{\bn}\|^{2}
              \quad\text{where} \; \; w(\bn) = \frac{ \bn !}{|\bn| !},
$$
          is finite.  We abbreviate $\ell^{2}_{w, \bbC}(\bbZ^{d}_{+})$ to
          $\ell^{2}_{w}(\bbZ^{d}_{+})$ and observe that
\begin{equation}
\ell^{2}_{w,\cY}(\bbZ^{d}_{+}) = \ell^{2}_{w}(\bbZ^{d}_{+}) \otimes \cY.
\label{2.2aa}
\end{equation}
          The commutative $d$-variable $Z$-transform
          $$ \{f_{\bn}\}_{\bn \in \bbZ^{d}_{+}} \mapsto \widehat f^{\ba}(\blam)=
          \sum_{\bn \in \bbZ^{d}_{_{+}}} f_{\bn} \blam^{\bn}
          $$
          maps $\ell_{w}^{2}(\bbZ_+^d)$ unitarily onto the {\em Arveson space}
$$
          \cH(k_d):=\left\{f(\blam)=\sum_{\bn \in{\mathbb Z}^d_{+}}f_{\bn}
          \blam^\bn: \; \|f\|^{2}=\sum_{\bn \in {\mathbb Z}^d_{+}}
          \frac{\bn!}{|\bn|!}\cdot |f_{\bn}|^2<\infty\right\}
$$
          with inner product given by
          $$
            \langle f,g \rangle_{\cH(k_{d})}=\sum_{\bn \in \bbZ^{d}_{+}}
            \frac{\bn !}{|\bn|!} f_{n}\overline{g_{n}}
$$
if
$$
            f(\blam) = \sum_{\bn \in \bbZ^{d}_{+}} f_{\bn}
\blam^{\bn}\quad\text{and}\quad
            g(\blam) = \sum_{\bn \in \bbZ^{d}_{+}} g_{\bn} \blam^{\bn}.
          $$
          Then it follows that the set
          $\{ \sqrt{\frac{|\bn| !}{\bn !}} \blam^{\bn} \colon \bn
          \in \bbZ^{d}_{+}\}$ is an orthonormal basis for
          $\ell^{2}_{w}(\bbZ^{d}_{+})$. By general principles concerning
          reproducing kernel Hilbert spaces we see that
          $\cH(k_{d})$ is a reproducing kernel Hilbert space of functions
          analytic on the unit ball
$$
\B^{d} = \{ \blam = (\lambda_{1},
          \dots, \lambda_{d}) \in \bbC^{d} \colon \sum_{k=1}^{d}
          |\lambda_{k}|^{2} < 1 \}
$$
with reproducing kernel $k_{d}(\blam, \bzeta)$ given by
          $$
          k_{d}(\blam, \bzeta) = \sum_{\bn \in \bbZ^{d}_{+}}
          \frac{|\bn|!}{\bn!} \blam^{\bn} \overline{\bzeta}^{\bn} =
          \sum_{n=0}^{\infty} \left( \lambda_{1} \overline{\zeta_{1}} +
          \cdots + \lambda_{d} \overline {\zeta_{d}}\right)^{n}
          = \frac{1}{1 - \langle \blam, \bzeta \rangle}
          $$
          (see e.g. \cite{ArvesonIII}).
          This justifies the notation $\cH(k_d)$ for the space. In analogy to
          \eqref{2.2aa} we will use notation
          $\cH_{\cY}(k_d):=\cH(k_d)\otimes\cY$ for the tensor product
Hilbert space
          that is characterized by
          $$
          \cH_{\cY}(k_d)=\left\{f(\blam)=\sum_{\bn \in{\mathbb Z}^d_{+}}f_{\bn}
          \blam^\bn:\|f\|^{2}=\sum_{\bn \in {\mathbb Z}^d_{+}}
          \frac{\bn!}{|\bn|!}\cdot \|f_{\bn}\|_{\cY}^2<\infty\right\}.
          $$
          If we define the map $\Pi$ by
          \begin{equation}  \label{Pi}
             \Pi \colon \{f_{v} \}_{v \in \cF_{d}} \mapsto \left\{ 
\sum_{v \colon
             \ba(v) = \bn} f_{v} \right\}_{\bn \in \bbZ^{d}_{+}},
          \end{equation}
          then each basis vector $\chi_{v}\in\ell^{2}(\cF_{d})$
          ($v \in \cF_{d}$) is mapped via $\Pi$
          to its abelianization $\chi_{\bn} \in \ell^{2}_{w}(\bbZ^{d}_{+})$ and
          then $\Pi$ is extended to the whole space $\ell^{2}(\cF_{d})$ via
          linearity.   The norm on $\ell^{2}_{w}(\bbZ^{d}_{+})$ is 
arranged so as
          to make $\Pi$ a coisometry from $\ell^{2}_{\cY}(\cF_{d})$ onto
          $\ell^{2}_{w, \cY}(\bbZ^{d}_{+})$ with initial space equal to
          $\ell^{2}_{\cY}({\mathcal S} \cF_{d})$ and with kernel equal to the
            subspace $\ell^{2}_{\cY}({\mathcal S} \cF_{d})^{\perp}$ of
          $\ell^{2}_{\cY}(\cF_{d})$ given by
          $$
              \ell^{2}_{\cY}({\mathcal S} \cF_{d})^{\perp} =
              \left\{ \{ f_{v} \}_{v \in \cF_{d}} \colon \sum_{v \in
\cF_{d} \colon
              \ba(v) = \bn} f_{v} = 0 \; \text{ for each } \; \bn \in
\bbZ^{d}_{+}
          \right\}.
          $$
          If we introduce the $Z$-transformed version $\widehat \Pi \colon
          H^{2}(\cF_{d}) \to \cH(k_{d})$ via
          $$ \widehat \Pi \colon \sum_{v \in \cF_{d}} f_{v} z^{v} \mapsto
             \sum_{\bn \in \bbZ^{d}_{+}} \left[ \sum_{v \in \cF_{d} \colon
             \ba(v) = \bn} f_{v} \right] \blam^{\bn},
          $$
          then similarly $\widehat \Pi$ is a coisometry from $\cH^{2}(\cF_{d})$
          onto $\cH(k_{d})$ with initial space equal to the subspace
          $$
             H^{2}({\mathcal S} \cF_{d}) : = \left\{ \sum_{
            v \in \cF_{d}} f_{\ba(v)} z^{v}
             \colon \sum_{\bn \in \bbZ^{d}_{+}} |f_{\bn}|^{2}< \infty \right\}
          $$
          with kernel equal to
          $$
            H^{2}({\mathcal S} \cF_{d})^{\perp} = \left\{ \sum_{v \in
\cF_{d}} f_{v}
            z^{v} \in H^{2}(\cF_{d}) \colon \sum_{v \colon \ba(v) = \bn} f_{v} =
            0 \; \text{ for each } \; \bn \in \bbZ^{d}_{+} \right\}.
          $$
          This gives the natural link between the Fock-space norm on 
formal power
          series and the Arveson-space norm on analytic functions on the unit
          ball and is the basis for the application of noncommutative results
          to prove commutative results in \cite{Arias-Popescu, DP, Popescu-int}.

          By a commutative $d$-dimensional linear system we mean a linear system
          with evolution along the integer lattice $\bbZ^{d}_{+}$ rather than
          along the free semigroup $\cF_{d}$.  A particular type of such a
          system is a system of the Fornasini-Marchesini form given by
	 \eqref{2.14}.
          If we specify an initial condition $x(0) = x^{0} \in \cX$ along
          with an input
          sequence $\{u^{0}(\bn) \}_{\bn \in \bbZ^{d}_{+}}$ and impose
          the boundary
          conditions that $x(\bn) = 0$ whenever $\bn$ is outside the positive
          orthant $\bbZ^{d}_{+}$, then the system equations uniquely determine a
          full system trajectory $\{ u(\bn), x(\bn), y(\bn)\}$ consistent with
          $x(0) = x^{0}$ and $u(\bn) = u^{0}(\bn)$ for $\bn \in \bbZ^{d}_{+}$.

          If $\Pi$ is the projection map introduced in \eqref{Pi} formally
          extended to be defined on all $\cF_{d}$-indexed sequences to generate
          a $\bbZ^{d}_{+}$-indexed sequence
          \begin{align*}
	& \Pi \colon \{u(v) \}_{v \in \cF_{d}} \mapsto \left\{ \sum_{v \colon
          \ba(v) = \bn} u(v) \right\}_{\bn \in \bbZ^{d}_{+}},  \\
          & \Pi \colon \{x(v) \}_{v \in \cF_{d}} \mapsto \left\{ \sum_{v \colon
          \ba(v) = \bn} x(v) \right\}_{\bn \in \bbZ^{d}_{+}}, \\
          & \Pi \colon \{y(v) \}_{v \in \cF_{d}} \mapsto \left\{ \sum_{v \colon
          \ba(v) = \bn} y(v) \right\}_{\bn \in \bbZ^{d}_{+}},
          \end{align*}
          then one can check the claim: {\em  $\{ (\Pi u)(\bn), (\Pi x)(\bn),
          (\Pi y)(\bn) \}_{\bn \in \bbZ^{d}_{+}}$ satisfies the system
          equations \eqref{2.14} whenever $\{u(v), x(v), y(v) \}_{v \in
          \cF_{d}}$ satisfies the system equations \eqref{2.3}.}
          Indeed the first system equation in \eqref{2.3} can be rewritten in
          the form
          $$
             x(v) = \sum_{k=1}^{d} A_{k} x(k^{-1} v) + \sum_{k=1}^{d} B_{k}
	    u(k^{-1} v).
          $$
          Here we use the convention that
          $$  k^{-1} v = \begin{cases} v' & \text{if } v = k v' , \\
	     \text{undefined} & \text{otherwise}
	    \end{cases}
          $$
          for $v$ a word in $\cF_{d}$ and $k \in \{1, \dots, d \}$ a letter
          and that $x(k^{-1} v)$ is interpreted to be $0$ if $k^{-1} v$ is
          undefined.
          Summing over $v \in \cF_{d}$ with $\ba(v) = \bn$ then gives
          $$ (\Pi x)(\bn) = \sum_{k=1}^{d} A_{k} \sum_{v \colon \ba(v) = \bn}
          x( k^{-1} v) + \sum_{k=1}^{d} B_{k}   \sum_{v \colon \ba(v) = \bn}
	 u(k^{-1} v).
          $$
          Now observe that
          $$
          \{  k^{-1} v \colon \ba(v) = \bn \} = \{ v' \in \cF_{d} \colon
          \ba(v') = \bn - e_{k}\}
          $$
          and arrive at
          $$
             (\Pi x)(\bn) = \sum_{k=1}^{d} A_{k} (\Pi x)(\bn -e_{k}) +
             \sum_{k=1}^{d} B_{k} (\Pi u)(\bn - e_{k}).
          $$
          We see that $\{ (\Pi u)(\bn), (\Pi x)(\bn), (\Pi y)(\bn) \}$
          satisfies the first of the system equations \eqref{2.14}.  That
          $\{ (\Pi u)(\bn), (\Pi x)(\bn), (\Pi y)(\bn) \}$ satisfies the second
          system equation in \eqref{2.14} is a simple consequence of linearity.
          Conversely, given a trajectory $\{u(\bn), x(\bn), y(\bn) \}_{\bn \in
          \bbZ^{d}_{+}}$ of \eqref{2.14}, let $\{u_{\ell}(v)\}_{v \in
          \cF_{d}}$ be any
          $\cU$-valued $\cF_{d}$-indexed sequence such that $\Pi u_{\ell} = u$
          and set $x_{\ell}(\emptyset) = x(0)$.  Then the noncommutative
          system equations \eqref{2.3} recursively uniquely determine a
          full system trajectory $\{u_{\ell}(v), x_{\ell}(v), y_{\ell}(v)\}_{v
          \in \cF_{d}}$ of \eqref{2.3} with this preassigned input string and
          initial condition.  By the claim verified above, it follows that
          $(\Pi u_{\ell}, \Pi x_{\ell}, \Pi y_{\ell})$ is again a system
          trajectory.  By the uniqueness of solution of the initial value
          problem for the system \eqref{2.14}, it follows that
          $\{(\Pi u, \Pi x, \Pi y)\} = \{ ( u,x,y)\}$.
          Thus,  any trajectory $\{u(\bn), x(\bn), y(\bn)\}_{\bn \in
          \bbZ^{d}_{+}}$ can be lifted to a trajectory $\{u_{\ell}(v),
          x_{\ell}(v), y_{\ell}(v)\}_{v \in \cF_{d}}$ of \eqref{2.3}, i.e.,
          $\{ u_{\ell}(v), x_{\ell}(v), y_{\ell}(v) \}_{v \in \cF_{d}}$ is a
          trajectory of \eqref{2.3} such that
          $$
          \{ \Pi u_{\ell}, \Pi x_{\ell},\Pi y_{\ell}\} = \{ u,x,y\}.
         $$
        In this way we view the Fornasini-Marchesini
          commutative system \eqref{2.14} as the abelianization of the
          noncommutative Fornasini-Marchesini system \eqref{2.3}.

          Since the commutative Fornasini-Marchesini system \eqref{2.14} is
          just the abelianization of the noncommutative Fornasini-Marchesini
          system \eqref{2.3} and we have already derived the formula
          \eqref{NC-IO} for the solution of the noncommutative initial-value
          problem, we see that the solution of the initial-value problem for
          the commutative Fornasini-Marchesini system \eqref{2.14} is simply
          the abelianization of the corresponding formula for the
          noncommutative case:
          \begin{equation}  \label{IO}
	 (\widehat \Pi  \widehat y)(\blam) = C (I - Z(\blam) A)^{-1} x(0)
            + T_{\Sigma}(\blam) \cdot (\widehat \Pi \widehat u)(\blam)
             \end{equation}
             where the {\em transfer function} $T_{\Sigma}(\blam)$ for the
             commutative Fornasini-Marchesini system is given by
$$
	  T_{\Sigma}(\blam) = D + C (I - Z(\blam) A)^{-1} Z(\blam) B.
$$
        This gives a derivation of the transfer function relationship
       \eqref{IO}
           (via the connection with noncommutative systems) which is an
             alternative to the usual direct approach via commutative
             multivariable $Z$-transform (see e.g.~\cite{MammaBear}).

          The zero input string simplifies the system to
\begin{equation}
\left\{\begin{array}{rcl}
          x(\bn) & = & A_{1} x(\bn -e_{1}) + \cdots + A_{d} x(\bn - e_{d})\\
          y(\bn) & = & C x(\bn).
          \end{array} \right.
\label{comsyst}
\end{equation}
          Given a pair $(C, \bA)$, we have the option of considering $(C, \bA)$
          as coming from a noncommutative or a commutative system.  If we
          consider the associated noncommutative system, the output string
          associated with initial state $x(\emptyset) = x$ (and zero input
          string) is the $\cY$-valued function on  $\cF_{d}$ given by
$$
              \cO_{C,{\bf A}} x =  \{ C \bA^{v}x \}_{v \in \cF_{d}}
$$
          and $(C, \bA)$ is considered output stable if this output string
is
          in $\ell_{\cY}^{2}(\cF_{d})$ for all $x \in \cH$.
          We say that the commutative system \eqref{comsyst} is {\em
output stable}
          (and in this case we will say that the pair $(C,\bA)$ is
          {\em $\ba$-output stable}) if $\Pi (\cO_{C,{\bf A}} x) \in
          \ell^{2}_{w,\cY}(\bbZ^{d}_{+})$ for all $x \in \cH$, or equivalently,
          if $\widehat\Pi \widehat \cO_{C,{\bf A}} x)$ is in the Arveson space
          $\cH_{\cY}(k_{d})$ for all choices of initial state $x \in \cX$.  We
          note that
          $\widehat\Pi \widehat \cO_{C,{\bf A}} x$ can be computed explicitly as
          $$
          (\widehat\Pi \widehat \cO_{C,{\bf A}} x)(\blam) = C (I - Z(\blam)
          A)^{-1} x.
          $$
          Thus another equivalent formulation of $\ba$-output stability is:
\begin{definition}
A pair $(C, \bA)$ is $\ba$-output stable means that the function
$C(I - Z(\blam)A)^{-1} x$ belongs to  $\cH_{\cY}(k_{d})$
          for every $x \in \cH$, {\rm or equivalently (by the closed graph
theorem)}, the  operator $\widehat {\mathcal O}^{\ba}_{C,\bA}$ from
          $\cX$ to $\cH_\cY(k_d)$ defined by
          \begin{equation}
          \widehat {\mathcal O}^{\ba}_{C,\bA} = \widehat \Pi \widehat {\mathcal
          O}_{C, \bA} \colon  x \mapsto    C(I-Z(\blam) A)^{-1} x=
          \sum_{\bn \in {\mathbb Z}^{d}_{+}} \left(\sum_{v\in{\ba}^{-1}(\bn)}C
          \bA^vx\right)\, \blam^{\bn}
          \label{2.20}
          \end{equation}
          is bounded.
\label{D:3.10}
\end{definition}
          The inverse $Z$-transform sends the function
          $$\widehat y^{\ba}(\blam) = C (I - Z(\blam) A)^{-1} x
              = \sum_{\bn \in \bbZ^{d}_{+}} \left( \sum_{v \in \ba^{-1}(\bn)} C
              \bA^{v} x \right) \blam^{\bn}
          $$
             to the string
          $\{y(\bn)\}_{\bn\in\bbZ_+^d}$ with
          \begin{equation}
          y(\bn)=\sum_{v\in{\ba}^{-1}(\bn)}C\bA^vx,\quad\bn\in\bbZ^d_+
          \label{2.21}
          \end{equation}
          and $\widehat y^{\ba}$ belongs to $\cH_\cY(k_d)$ if and only if
          $\{y(\bn)\}_{\bn\in\bbZ_+^d}\in\ell_{w,\cY}^{2}(\bbZ^d_+)$.
          Thus, the operator $\widehat {\mathcal O}^{\ba}_{C,\bA}$ introduced in
          \eqref{2.20} is the $Z$-transformed version of the
observability operator
          \begin{equation}
          {\cO}^{\ba}_{C,\bA}\colon \; x\mapsto\left\{
          \sum_{v\in{\ba}^{-1}(\bn)}C\bA^vx\right\}_{\bn\in\bbZ^d_+}
          \label{2.22}
          \end{equation}
          and a pair $(C,\bA)$ is $\ba$-output stable if and only if
          ${\cO}^{\ba}_{C,\bA}$ is bounded as an operator from $\cX$ into
          $\ell_{w,\cY}^{2}(\bbZ^d_+)$.
          In this case it makes sense to introduce the {\em
observability gramian}
          $$
          {\mathcal G}^{\ba}_{C, \bA}:=({\cO}^{\ba}_{C, \bA})^{*}{\cO}^{\ba}_{C,
          \bA}=(\widehat{\cO}^{\ba}_{C, \bA})^{*}\widehat{\cO}^{\ba}_{C, \bA}
          $$
          and its representation in terms of strongly converging series
          \begin{equation}
          {\mathcal G}^{\ba}_{C,\bA}=\sum_{\bn\in{\mathbb
          Z}^{d}_{+}}\frac{\bn!}{|\bn|!}
          \left(\sum_{v,u\in\ba^{-1}(\bn)}\bA^{*v^\top}C^*C\bA^{u}\right)
          \label{2.23}
          \end{equation}
          follows immediately by definitions  \eqref{2.22}, \eqref{2.20}
          and the formulas for the inner products in
$\ell_{w,\cY}^{2}(\bbZ^d_+)$ and
          $\cH_\cY(k_d)$.

          \begin{definition}
          We say that the pair $(C,\bA)$ is {\em $\ba$-observable} if ${\mathcal
          G}^{\ba}_{C, \bA}$ is positive-definite and  {\em  exactly
          $\ba$-observable} if ${\mathcal G}^{\ba}_{C, \bA}$ is 
strictly positive
          definite.
          \label{D:2.2}
          \end{definition}

         By Theorem \ref{T:2-1.1} (2) we know that the observability
         gramian ${\mathcal G}_{C, {\mathbf A}}$ satisfies the Stein
         equation \eqref{3.4}.  It turns out that the abelianized
         observability gramian ${\mathcal G}^{\ba}_{C, {\mathbf A}}$
         satisfies a reverse Stein inequality (the reverse of \eqref{3.4a}).

         \begin{proposition} \label{P:reverseStein}
             Let $(C, {\mathbf A})$ be an $\ba$-output-stable pair and let
${\mathcal
             G}^{\ba}_{C, {\mathbf A}}$ be the abelianized observability
             gramian \eqref{2.23}.  Then ${\mathcal G}^{\ba}_{C, {\mathbf
             A}}$ satisfies the reverse Stein inequality
             \begin{equation}
             {\mathcal G}^{\ba}_{C,\bA} - A_{1}^{*}{\mathcal G}^{\ba}_{C,\bA}
             A_{1}- \cdots -A_{d}^{*}{\mathcal G}^{\ba}_{C,\bA} A_{d} \le C^*C.
             \label{3.7}
             \end{equation}
             Moreover, the following are equivalent:
             \begin{enumerate}
	  \item Equality holds in \eqref{3.7}.
	  \item  $\bA$ is $C$-abelian
             in the sense that
             \begin{equation}
             C{\mathbf A}^v=C\bA^u,\quad \mbox{whenever $v,u\in\cF_d$ and
             $\ba(v)=\ba(u)$}.
             \label{3.8}
             \end{equation}
            \item The observability gramian and the abelianized
            observability gramian are identical:
            $$
            {\mathcal G}^{\ba}_{C,{\mathbf A}} = {\mathcal G}_{C, {\mathbf
            A}}.
            $$
            \end{enumerate}
          \end{proposition}

          \begin{proof}
              It suffices to show that the operator $Q$ given by
              \begin{equation}
              Q:=C^*C-{\mathcal G}^{\ba}_{C,\bA}+\sum_{j=1}^d
              A_{j}^{*}{\mathcal G}^{\ba}_{C,\bA} A_{j}
              \label{3.9}
              \end{equation}
              is positive semidefinite.  To this end, plug \eqref{2.23} into
              \eqref{3.9} to get
              \begin{equation}
              Q=\sum_{N=1}^\infty Q_N
              \label{3.10}
              \end{equation}
           where $Q_{N}$ is given by
           \begin{align}
           Q_N &= \sum_{j=1}^d A_j^*\left[
           \sum_{\bm\in{\mathbb Z}^{d}_{+}\colon |\bm|=N-1}
           \frac{\bm!}{(N-1)!}\sum_{v,u\in\ba^{-1}(\bm)}
           \bA^{*v^\top}C^*C\bA^u\right]A_j  \notag \\
           & \qquad -\sum_{\bn\in{\mathbb Z}^{d}_{+}\colon
|\bn|=N}\frac{\bn!}{N!}
           \sum_{v,u\in\ba^{-1}(\bn)}\bA^{*v^\top}C^*C\bA^u.
           \label{3.11}
           \end{align}
           We introduce the notation
           \begin{equation}
           W(\bn)=\sum_{u\in\ba^{-1}(\bn)}\bA^u\quad (\bn\in{\mathbb Z}^{d}_{+})
           \label{3.12}
           \end{equation}
           and extend the notation to the all of ${\mathbb Z}^{d}$ by
           \begin{equation}
           W(\bn)= 0 \text{ if }  \bn \in {\mathbb Z}^{d}\setminus
           {\mathbb Z}^{d}_{+}.
           \label{3.13}
           \end{equation}
           With these definitions we have the equality
           \begin{equation}
           \label{3.14}
           W(\bn)=\sum_{i=1}^d W(\bn-e_i)A_i
           \qquad \mbox{for every $\bn=(n_1,\ldots, n_d) \in{\mathbb 
Z}^{d}_{+}$}
           \end{equation}
           where $e_1,\ldots,e_d\in\bbZ_+^d$ are defined in \eqref{2.16}.
           Write formula \eqref{3.11} in terms of \eqref{3.12} as
           \begin{eqnarray}
              Q_N&=&\sum_{j=1}^d A_j^*
           \left[\sum_{\bm\in{\mathbb Z}^{d}_{+}:\, |\bm|=N-1}
           \frac{\bm!}{(N-1)!}W(\bm)^*C^*CW(\bm)\right]A_j\nonumber\\
           &&-\sum_{\bn\in{\mathbb Z}^{d}_{+}:\, |\bn|=N}
           \frac{\bn!}{N!}\cdot W(\bn)^*C^*CW(\bn).
           \label{3.15}
           \end{eqnarray}
           Upon rearranging the terms in the first series in \eqref{3.15}
           and substituting formula \eqref{3.14} into the second series,
           we arrive at
           \begin{eqnarray}
           Q_N&=&\sum_{\bn\in{\mathbb Z}^{d}_{+}:\, |\bn|=N}
           \sum_{j=1}^d\frac{(\bn-e_j)!}{(N-1)!}
           A_j^*W(\bn-e_j)^*C^*CW(\bn-e_j)A_j \label{3.16}\\
           &&-\sum_{\{\bn\in{\mathbb Z}^{d}_{+}:\, |\bn|=N\}}
           \left(\sum_{i,j=1}^d\frac{\bn!}{N!}
           A_i^*W(\bn-e_i)^*C^*CW(\bn-e_j)A_j\right).\nonumber
           \end{eqnarray}
           We now consider the terms in \eqref{3.16} that correspond to a fixed
           $\bn=(n_1,\ldots,n_d)\in{\mathbb Z}^{d}_{+}$ (with $|\bn|=N$).
           Denoting the sum of these terms by $S_{\bn}$ we have
        \begin{align}
	S_{\bn} = & \sum_{j=1}^{d} \frac{(\bn - e_{j})!}{(N-1)!}
	A_{j}^{*}W(\bn - e_{j})^{*} C^{*}C W(\bn - e_{j}) A_{j} \notag \\
	& \qquad \qquad - \frac{\bn !}{N!} \sum_{i,j=1}^{d} A_{i}^{*}W(\bn -
	e_{i})^{*} C^{*}C W(\bn - e_{j}) A_{j} \notag \\
	= &\sum_{j=1}^{d} \frac{(\bn - e_{j})!}{N!} (N - n_{j})
	A_{j}^{*}W(\bn - e_{j})^{*} C^{*}C W(\bn - e_{j}) A_{j} \notag \\
	& \qquad \qquad - \frac{ \bn !}{N!} \sum_{i,j \in \{1, \dots,
d\} \colon i
	\ne j} A_{i}^{*} W(\bn - e_{i})^{*} C^{*}C W(\bn -e_{j}) A_{j}.
\label{3.17}
        \end{align}
Note that by convention \eqref{3.13}, the indices $i$ and $j$ in the
           latter summations vary on the set
           $$
           {\mathcal I}_{\bn}=\{\ell\in\{1,\ldots,d\}: \; n_\ell>0\}
           $$
           rather than $\{1,\ldots,d\}$. Furthermore, since
           $$
           N-n_j=|\bn|-n_j=\sum_{i\in{\mathcal I}_{\bn}: \, i\neq j}n_i
           $$
           and
           $$
           (\bn-e_j)!=(\bn-e_j-e_i)! \, n_i \quad (i\neq j),
           $$
           one can rewrite the first sum on the right hand side in
\eqref{3.17} as
           \begin{align*}
           & \sum_{i,j\in{\mathcal I}_{\bn} \colon i\neq j}
           \frac{(\bn-e_j)!}{N!} \, n_i A_j^*W(\bn-e_j)^*C^*CW(\bn-e_j)A_j \\
           & \qquad =\sum_{i,j\in {\mathcal I}_{\bn} \colon  i\neq j}
           \frac{(\bn-e_j-e_i)!}{N!} \, n_i^2A_j^*W(\bn-e_j)^*C^*CW(\bn-e_j)A_j.
           \end{align*}
           Plugging this into the right hand side in \eqref{3.17} leads us to
           \begin{eqnarray}
           S_{\bn}&=&\sum_{i,j\in{\mathcal I}_{\bn}: \, i\neq j}
           \frac{(\bn-e_j-e_i)!}{N!}
           \left[n_i^2A_j^*W(\bn-e_j)^*C^*CW(\bn-e_j)A_j\right.\nonumber\\
           &&\qquad\qquad\left.-n_in_j
A_i^*W(\bn-e_i)^*C^*CW(\bn-e_j)A_j\right]\nonumber\\
           &=&\sum_{i,j\in{\mathcal I}_{\bn}: \, i\neq j}
           \frac{1}{2} \frac{(\bn-e_j-e_i)!}{N!}
           R_{\bn,i,j}^*R_{\bn,i,j}\label{3.18}
           \end{eqnarray}
           where
           \begin{equation}
           R_{\bn,i,j}=C \left[n_iW(\bn-e_j)A_j-n_jW(\bn-e_i)A_i\right].
           \label{3.19}
           \end{equation}
           Representation \eqref{3.18} implies that $S_{\bn}$ is positive
           semidefinite and therefore $Q_N\ge 0$ for every $N\in{\mathbb N}$.
           By \eqref{3.10}, the operator $Q$ defined in \eqref{3.9} is positive
           semidefinite which completes the proof of \eqref{3.7}.

           We now show the equivalence of (1), (2) and (3) in the second
           part of Proposition  \ref{P:reverseStein}.

           \medskip
           {\em Proof of (1) $\Longrightarrow$ (2):}
           Assume condition (1), i.e., that the reverse Stein inequality
           \eqref{3.7} is satisfied with equality.
            Then representation \eqref{3.18} implies
           that $R_{\bn,i,j}=0$ for all $\bn\in{\mathbb Z}^{d}_{+}$.
           By \eqref{3.19}, this means that
           \begin{equation}
           n_i CW(\bn-e_j)A_j=n_jCW(\bn-e_i)A_i\qquad (\bn\in{\mathbb
Z}^{d}_{+}).
           \label{3.20}
           \end{equation}
           Now we shall prove \eqref{3.8} by induction (on the length of words
           $v,u\in\cF_d$). The basis of induction ( $|v|=|u|=0$) is
trivial. Assume
           that \eqref{3.8} holds true, whenever $|v|=|u|<N$. Then in 
particular,
           we have for every $\bm\in{\mathbb Z}^{d}_{+}$ with $|\bm|<N$:
           \begin{equation}
           CW(\bm)=\sum_{w\in\ba^{-1}(\bm)}C\bA^w=\frac{|\bm|!}{\bm!} C\bA^{w_0}
           \quad \mbox{for every} \; \; w_0\in\ba^{-1}(\bm).
           \label{3.21}
           \end{equation}
           Now take two words
           $v,u\in\cF_d$ of the length $N$ and let
           \begin{equation}
           \ba(v)=\ba(u)=:\bn=(n_1,\ldots,n_d).
           \label{3.22}
           \end{equation}
           If $v=\widetilde{v} i$ and $u=\widetilde{u}i$ for some
           $\widetilde{v},\widetilde{u}\in\cF_d$ and $i\in\{1,\ldots,d\}$, then
           we have $C\bA^{\widetilde{v}}=C\bA^{\widetilde{u}}$ by the induction
           hypothesis and therefore,
           $$
           C\bA^v=C\bA^{\widetilde{v}}A_i=C\bA^{\widetilde{u}}A_i=C\bA^u.
           $$
           Let $v=\widetilde{v} i$ and $u=\widetilde{u}j$ for some
           $i,j\in\{1,\ldots,d\}$ and $i\neq j$. By \eqref{3.22},
           $\ba(\widetilde{v})=\bn-e_i$ and $\ba(\widetilde{u})=\bn-e_j$.
           By \eqref{3.21}, we have
           \begin{eqnarray}
     CW(\bn-e_j)&=&\frac{(N-1)!}{(\bn-e_j)!}C\bA^{\widetilde{v}},\label{3.23}\\
CW(\bn-e_i)&=&\frac{(N-1)!}{(\bn-e_i)!}C\bA^{\widetilde{u}}.\label{3.24}
           \end{eqnarray}
           Multiplying \eqref{3.23} and \eqref{3.24} on the right by 
$n_iA_j$ and
           $n_jA_i$ respectively, we get
           \begin{eqnarray*}
           n_iCW(\bn-e_j)A_j=n_i\frac{(N-1)!}{(\bn-e_j)!}C\bA^{\widetilde{v}}A_j
           =n_i\frac{(N-1)!}{(\bn-e_j)!}C\bA^{\widetilde{v}j}=
           n_in_j\frac{(N-1)!}{\bn!}C\bA^v
           \end{eqnarray*}
           and
           \begin{eqnarray*}
           n_jCW(\bn-e_i)A_i=n_j\frac{(N-1)!}{(\bn-e_i)!}C\bA^{\widetilde{v}}A_i
           =n_j\frac{(N-1)!}{(\bn-e_i)!}C\bA^{\widetilde{u}i}=
           n_jn_i\frac{(N-1)!}{\bn!}Y\bA^u.
           \end{eqnarray*}
           By \eqref{3.20}, the left hand side expressions in the two latter
           equalities are equal. Upon comparing the right hand side
expressions we
           get $C\bA^v=C\bA^u$ , i.e., ${\mathbf A}$ is $C$-abelian as
           wanted.

           \medskip
           {\em Proof of (2) $\Longrightarrow$ (3):}
         Assume now that ${\mathbf A}$ is $C$-abelian, i.e., that
\eqref{3.8} holds.
         Then the  identify ${\mathcal G}^{\ba}_{C,{\mathbf A}} = {\mathcal
         G}_{C,{\mathbf A}}$ is an immediate consequence of the
         series representations \eqref{2.23} and \eqref{2.11} for
         ${\mathcal G}^{\ba}_{C,{\mathbf A}}$ and ${\mathcal G}_{C,{\mathbf A}}$
         respectively.

          \medskip
          {\em Proof of (3) $\Longrightarrow$ (1):}  We know from Theorem
          \ref{T:2-1.2} (2) that ${\mathcal G}_{C, {\mathbf A}}$ satisfies
          the Stein equation, i.e., ${\mathcal G}_{C, {\mathbf A}}$
          satisfies the Stein inequality \eqref{3.7} with equality.  Hence
          trivially ${\mathcal G}^{\ba}_{C, {\mathbf A}}$ satisfies \eqref{3.7}
          with equality whenever ${\mathcal G}^{\ba}_{C, {\mathbf A}} =
          {\mathcal G}_{C, {\mathbf A}}$.  This completes the proof of
          Proposition \ref{P:reverseStein}.
          \end{proof}

          \begin{example}  \label{E:reverseStein}
         {\rm  If $(C,\bA)$ is an output-stable pair, then
         by Theorem \ref{T:2-1.1} (2) ${\mathcal G}_{C,\bA}$ satisfies the
         Stein equation \eqref{3.4} and hence in particular
        $$
        {\mathcal G}_{C,\bA}-A_{1}^{*}{\mathcal G}_{C,\bA}
        A_{1}-\ldots-A_{d}^{*}{\mathcal G}_{C,\bA} A_{d}\ge 0,
        $$
        We now show that, for the abelianized case,
        the inequality in the reverse Stein inequality satisfied by the
        abelianized observability gramian ${\mathcal G}^{\ba}_{C, {\mathbf A}}$
        can be strict in the strong sense that the quantity
        $ {\mathcal G}^{\ba}_{C,\bA} - A_{1}^{*}{\mathcal G}^{\ba}_{C,\bA}
        A_{1}-\ldots-A_{d}^{*}{\mathcal G}^{\ba}_{C,\bA} A_{d}
        $ is not even positive semidefinite. As an
        example, let
$$
C=\begin{bmatrix}1 & 0 & 0 \end{bmatrix},\quad
A_1=\begin{bmatrix} 0 & \frac{1}{2} & 0 \\ 0 & 0 & 0 \\ -\frac{1}{2} & 0 &
0\end{bmatrix},\quad A_2=\begin{bmatrix} 0 & 0 & \frac{1}{2} \\
\frac{1}{2} & 0 & 0
\\ 0 & 0 &0\end{bmatrix}.
$$
A straightforward calculation shows that
$$
{\mathcal G}^{\ba}_{C,\bA} - A_{1}^{*}{\mathcal G}^{\ba}_{C,\bA}
        A_{1}-\ldots-A_{d}^{*}{\mathcal G}^{\ba}_{C,\bA} A_{d}=
        \begin{bmatrix} \frac{7}{8} & \frac{5}{8} & \frac{3}{8} \\
            \frac{5}{8} & 0 & \frac{1}{4} \\ \frac{3}{8} & \frac{1}{4} &
0\end{bmatrix},
$$
which is not positive semidefinite}.
        \end{example}

          Condition \eqref{3.8} is worth a formal definition.
          \begin{definition}
          Let $C\in{\mathcal L}(\cX,\cY)$. A $d$-tuple $\bA=(A_1,\ldots,A_d)$
          of bounded operators on $\cX$ will be called {\em $C$-abelian} if
          \eqref{3.8} holds.
          \label{D:2.6}
          \end{definition}

          One obvious way for a given operator $d$-tuple ${\mathbf A}$ to
          be $C$-abelian is for ${\mathbf A}$ itself to be commutative,
          i.e., for $A_{i}A_{j} = A_{j}A_{i}$ for all $1 \le i,j \le d$.
          We next show that, under an observability assumption, this is the
          only way.

          \begin{proposition}  \label{P:Ccom}
            Suppose that the output-stable pair $(C, {\mathbf A})$ is observable
            and that ${\mathbf A}$ is $C$-abelian.
            Then the $d$-tuple $\bA$ is commutative.
              \end{proposition}

              \begin{proof}
              Since $\bA$ is $C$-abelian, relations \eqref{3.8} hold.
              Fix $i,j\in\{1,\ldots,d\}$ and note that by \eqref{3.8},
              $$
              C\bA^v A_iA_j=C\bA^{vij}=C\bA^{vji}
              =C\bA^v A_jA_i \quad\mbox{for every $v\in\cF_d$},
              $$
              since $\ba(vij)=\ba(vji)$. Thus,
              $$
              C\bA^v(A_iA_j-A_jA_i)x=0
              $$
              for every $v\in\cF_d$ and $x\in\cX$. Since the pair $(C,\bA)$
is observable,
              we have by \eqref{2.13}
              $$
              (A_iA_j-A_jA_i)x=0
              $$
              holding for every $x\in\cX$, which proves the
commutativity relations
              $$
              A_iA_j=A_jA_i\quad\mbox{for} \; \; i,j=1,\ldots,d
              $$
and completes the proof.
          \end{proof}

          \begin{corollary}  \label{C:Ccom}
              Suppose that $(C, {\mathbf A})$ is an observable output-stable
              pair.  Then the abelianized observability gramian coincides
              with the observability gramian
              $$ {\mathcal G}^{\ba}_{C, {\mathbf A}} = {\mathcal G}_{C,
              {\mathbf A}}
              $$
              if and only if the operator $d$-tuple ${\mathbf A}$ is
              commutative.
          \end{corollary}

          \begin{proof}
              Combine (2) $\Longleftrightarrow$ (3) in Proposition
              \ref{P:reverseStein} with Proposition \ref{P:Ccom}.
          \end{proof}

          We next show that the observability gramian always dominates the
          abelianized observability gramian.

          \begin{proposition}  \label{P:Gineq}
          Let $(C,{\mathbf A})$ be an output-stable pair. Then:
          \begin{enumerate}
              \item $(C,{\mathbf A})$ is also $\ba$-output-stable with
           \begin{equation}  \label{Gineq}
	{\mathcal G}^{\ba}_{C,{\mathbf A}} \le {\mathcal G}_{C,{\mathbf A}}.
           \end{equation}

           \item Equality occurs in \eqref{Gineq} if and only if ${\mathbf
           A}$ is $C$-abelian:
           $$
           C {\mathbf A}^{v} = C {\mathbf A}^{u} \text{ whenever } u,v \in
           \cF_{d} \text{ with } \ba(u) = \ba(v).
          $$
          \end{enumerate}
          \end{proposition}

          \begin{proof}  Note that the second statement in Proposition
          \ref{P:Gineq} is just a restatement
          of (2) $\Longleftrightarrow$ (3) in Proposition \ref{P:reverseStein}.
          Thus it suffices only to prove the first statement.

          By definition, output-stability of $(C, {\mathbf A})$ simply means
          that ${\mathcal G}_{C, {\mathbf A}}$ is bounded, while
          $\ba$-output stability means that ${\mathcal G}^{\ba}_{C, {\mathbf
          A}}$ is bounded.  The fact that $\ba$-output stability
          follows from output-stability therefore follows immediately from
the general
          inequality \eqref{Gineq}.  Thus it suffices to prove \eqref{Gineq}.
          For this purpose, recall that
          $$
          \langle {\mathcal G}_{C, {\mathbf A}}x, x \rangle =
          \left\|{\cO}_{C, \bA}x\right\|^{2}_{\ell^2_{\cY}(\cF_d)}=
          \sum_{v\in\cF_d}\|C\bA^vx\|^2_{\cY},
          $$
          while
          $$ \langle {\mathcal G}^{\ba}_{C, {\mathbf A}} x, x \rangle =
          \left\|{\cO}^{\ba}_{C, \bA}x\right\|^{2}_{\ell^2_{\cY}(\bbZ_+^d)}=
          \sum_{\bn\in{\mathbb Z}^{d}_{+}}\frac{\bn!}{|\bn|!}
          \left\|\sum_{v\in\ba^{-1}({\mathbf n})}C\bA^vx\right\|^2_{\cY}.
          $$
          By the Cauchy-Schwarz inequality we have
          $$
          \left\|\sum_{v\in\ba^{-1}({\mathbf n})}C\bA^vx\right\|^2_{\cY}
          \le \left( \sum_{v \in
          \ba^{-1}(\bn)} \| C \bA^{v}x\|_{\cY} \right)^{2}
          \le \sum_{v \in \ba^{-1}(\bn)} \| C \bA^{v} x \|_{\cY}^{2} \cdot
          \frac{|\bn|!}{\bn !}.
          $$
          Therefore
         \begin{eqnarray*}
         \left\| \cO^{\ba}_{C, \bA}x\right\|^{2}_{\ell^{2}_{w,
         \cY}(\bbZ^{d}_{+})}&=&\sum_{\bn\in{\mathbb
Z}^{d}_{+}}\frac{\bn!}{|\bn|!}
         \left\|\sum_{v\in\ba^{-1}({\mathbf
n})}C\bA^vx\right\|^2_{\cY}\nonumber\\
         &\le & \sum_{\bn\in{\mathbb Z}^{d}_{+}}\sum_{v \in \ba^{-1}(\bn)} \| C
         \bA^{v} x \|_{\cY}^{2}\nonumber\\
         &=&\sum_{v\in\cF_d}\|C\bA^vx\|^2_{\cY}= \|\cO_{C, \bA}x \|^{2}
         \end{eqnarray*}
         and \eqref{Gineq} follows as wanted.
          \end{proof}

          \begin{example} \label{E:a-stable}
           {\rm The converse of Proposition \ref{P:Gineq} part (1)
          can fail, i.e., {\em there exists an output pair $(C, {\mathbf A})$
          which is $\ba$-output-stable but not output-stable.}
         For example take
$$
          C=\begin{bmatrix}1 & 0 & 0 \end{bmatrix},\quad
          A_1=\begin{bmatrix} 0 & 2 & 0 \\ 0 & 0 & 0 \\ -1 & 0 &
          0\end{bmatrix},\quad A_2=\begin{bmatrix} 0 & 0 & 2 \\ 1 & 0 & 0
          \\ 0 & 0 &0\end{bmatrix}.
$$
          Then  $C(I-\lambda_1A_1-\lambda_2A_2)^{-1}=\begin{bmatrix} 1 &
          2\lambda_1 & 2\lambda_2\end{bmatrix}$. Hence
          $$ \widehat \cO^{\ba}_{C,\bA} \colon x \to C(I -
          \lam_{1}A_{1} - \lam_{2}A_{2})^{-1} x
          $$
          maps $\cX = {\mathbb C}^{3}$ into $\cH(k_{2})$
         and thus $(C,\bA)$ is ${\ba}$-output stable.
          To show that $(C,\bA)$ is not output stable,
          note that
$$
(A_1A_2)^n= \begin{bmatrix} 2^n & 0 & 0 \\ 0 & 0 & 0 \\
0 & 0 & (-2)^n  \end{bmatrix}
$$
and therefore, $C(A_1A_2)^n=\left [\begin{array}{ccc}2^n & 0 &
0\end{array}\right]$, so that for $x=\left [\begin{array}{ccc}1 & 0 &
0\end{array}\right]^\top$,
$$
\sum_{v\in\cF_d}\|C\bA^v x\|_{\C}^2\ge \sum_{n\ge 0}2^n=\infty
$$
and therefore, the pair $(C,\bA)$ is not output-stable. We
conclude that $\ba$-output-stability has no obvious
characterization in terms of positive semidefiniteness of some solution
of a Stein inequality as in the noncommutative case
(see Theorem \ref{T:2-1.1} (2)).}
\end{example}

          As a corollary of the gramian inequality \eqref{Gineq}
          in Proposition \ref{P:Gineq}, we have the following.

          \begin{corollary} \label{C:Gineq}
              Let $(C, {\mathbf A})$ be an output-stable pair. Then:
          \begin{enumerate}
              \item $\operatorname{Ker} {\mathcal G}_{C, {\mathbf A}}
              \subset \operatorname{Ker} {\mathcal G}^{\ba}_{C, {\mathbf A}}$.
              Hence, if $(C,\bA)$ is $\ba$-observable (respectively, exactly
              $\ba$-observable, then $(C, \bA)$ is
              also observable (respectively, exactly observable).

              \item The subspace $\operatorname{Ker} {\mathcal G}_{C,
              {\mathbf A}} = \operatorname{Ker} {\mathcal O}_{C, {\mathbf
              A}}$ is invariant under the operator $A_{j}$ for each $j = 1,
              \dots, d$.

             \item  The subspace $\operatorname{Ker} {\mathcal
              G}^{\ba}_{C,{\mathbf A}}$ is invariant under  $A_{j}$ for each
              $j = 1, \dots, d$ if and only if $\operatorname{Ker} {\mathcal
              G}^{\ba}_{C,{\mathbf A}} =
              \operatorname{Ker}{\mathcal G}_{C,{\mathbf A}}$.
\end{enumerate}
          \end{corollary}

\begin{proof}
Statement (1) is an immediate consequence of the inequality  \eqref{Gineq}.
Statement (2) is easily checked from the definition of
${\mathcal O}_{C, {\mathbf A}}$.  Sufficiency in statement (3)
is then a consequence of statement (2).  It remains only to verify
necessity in statement (3).

Assume therefore that
$\operatorname{Ker} {\mathcal G}^{\ba}_{C,\bA}$ is invariant
under $A_j$ for each $j=1,\ldots,d$. Let $x$ be a vector in
$\operatorname{Ker}{\mathcal G}^{\ba}_{C,\bA}$. Then by the assumed invariance,
$\bA^ux\in{\rm Ker} \, {\mathcal G}^{\ba}_{C,\bA}$ for every $u\in\cF_d$.
Then we have
$$
\operatorname{Ker} \sum_{v\in\ba^{-1}(\bn)}C\bA^{vu} x=0\quad
\text{for every} \; \; u\in\cF_d.
$$
Then letting $\bn=0$ we get $C\bA^ux=0$ for every $u\in\cF_d$ and
therefore, $x\in \operatorname{Ker} {\mathcal G}_{C,\bA}$. Thus,
$\operatorname{Ker}{\mathcal G}^{\ba}_{C,\bA} \subset  {\rm Ker} \, {\mathcal
G}_{C,\bA}$ and since the reverse inclusion holds by the first statement,
equality follows.
\end{proof}

\begin{example} \label{E:a-obs}
{\rm We observed in part (1) of Corollary \ref{C:Gineq} that
$\ba$-observability for an output-stable pair $(C, \bA)$ implies
observability.
We now give an example to show that the converse can fail, i.e., there
exists an output-stable observable pair which is not $\ba$-observable.
For this purpose, let $d=2$, $\cX = {\mathbb C}^4$,
$\cY ={\mathbb C}$, $C=\begin{bmatrix}0
& 0 & 0 & 1\end{bmatrix}$ and ${\bf A}=(A_1,A_2)$, where
$$
A_1=\begin{bmatrix} -\frac{1}{16}&\frac{1}{16}&0&0\\
-\frac{1}{16}&\frac{1}{16}&-\frac{1}{16}&\frac{1}{16}\\ 0&0&0&0\\
0&0&-\frac{1}{16}&\frac{1}{16}\end{bmatrix},
\qquad A_2=\begin{bmatrix} \frac{1}{16}&0&0&-\frac{1}{16}\\
         -\frac{1}{16}&-\frac{1}{16}&-\frac{1}{16}&-\frac{1}{16}\\
         \frac{1}{16}&-\frac{1}{16}&\frac{1}{16}&-\frac{1}{16}\\
         -\frac{1}{16}&0&0&-\frac{1}{16} \end{bmatrix}.
$$
Then the pair $(C,\bA)$ is output stable. Now we show that $(C,{\bA})$ is
observable but not $\ba$-observable. Indeed, since
$$
CA_1A_2=\begin{bmatrix}-\frac{1}{128} & \frac{1}{256} & -\frac{1}{256} & 0
\end{bmatrix},
$$
we have
$$
CA_1=\frac{1}{16}\begin{bmatrix}0 & 0 & -1 & 1\end{bmatrix}, \quad
CA_2=-\frac{1}{16}\begin{bmatrix}1 & 0 & 0 & 1\end{bmatrix}
$$
and
$$
CA_1A_2=-\frac{1}{256}\begin{bmatrix}2 & -1 & 1 & 0\end{bmatrix}.
$$
Now it is clear that
$\operatorname{Ker} C\cap \operatorname{Ker} CA_1\cap
\operatorname{Ker}CA_2\cap
\operatorname{Ker}CA_1A_2=0$ which implies that
$\bigcap_{v\in\cF_d} \operatorname{Ker} C{\bf A}^v = 0$. Therefore, the pair
$(C,{\bf A})$ is observable. To show that $(C,{\bf A})$ is not
$\ba$-observable we first compute
$$
I-\lambda_1A_1-\lambda_2A_2=\begin{bmatrix}
1+\frac{\lambda_1}{16}-\frac{\lambda_2}{16} &-\frac{\lambda_1}{16}&0&
\frac{\lambda_2}{16}\\ \frac{\lambda_1}{16}+
\frac{\lambda_2}{16}&1-\frac{\lambda_1}{16}+\frac{\lambda_2}{16}&\frac{\lambda_1}{16}
+\frac{\lambda_2}{16}&-\frac{\lambda_1}{16}+\frac{\lambda_2}{16}\\
-\frac{\lambda_2}{16}&\frac{\lambda_2}{16}&1-\frac{\lambda_2}{16}&\frac{\lambda_2}{16}\\
\frac{\lambda_2}{16}&0&\frac{\lambda_1}{16}&1-\frac{\lambda_1}{16}+\frac{\lambda_2}{16}
\end{bmatrix}.
$$
A straightforward calculation gives
\begin{eqnarray*}
d(\lambda_1,\lambda_2)&:=&\det \, (I-\lambda_1A_1-\lambda_2A_2)\\ &=&
1-\frac{\lambda_1}{16}+\frac{\lambda_1\lambda_2}{128}-
\frac{\lambda_1^{2}\lambda_2}{2048}-\frac{\lambda_2^{2}}{64}+
\frac{\lambda_1\lambda_2^{2}}{2048}-\frac{\lambda_1\lambda_2^{3}}{16384}+
\frac{\lambda_2^{4}}{16384}.
\end{eqnarray*}
Note that
$$
\begin{bmatrix}y_1 & y_2 & y_3 & y_4\end{bmatrix}:=
C(I-\lambda_1A_1-\lambda_2A_2)^{-1}
$$
is the bottom row of the matrix $(I-\lambda_1A_1+\lambda_2A_2)^{-1}$
and we use the standard adjoint formula  for the inverse of a matrix  to get
$$
y_2=\frac{1}{d(\lambda_1,\lambda_2)}\cdot\left |
\begin{array}{ccc}1+\frac{\lambda_1}{16}-\frac{\lambda_2}{16}&-\frac{\lambda_1}{16}&0\\
-\frac{\lambda_2}{16}&\frac{\lambda_2}{16}&1-\frac{\lambda_2}{16}\\
\frac{\lambda_2}{16}&0&\frac{\lambda_1}{16}
            \end{array} \right|\equiv 0.
$$
Then it follows that the nonzero vector $x=\begin{bmatrix}0 & 1 & 0 &
0\end{bmatrix}^\top$ satisfies
$$
C(I-\lambda_1A_1-\lambda_2A_2)^{-1}x\equiv 0
$$
and therefore, the pair $(C,{\bf A})$ is not ${\bf a}$-observable.}
\end{example}

\subsection{Observability-operator range spaces and reproducing
kernel Hilbert spaces: the commutative-variable case} \label{S:C-Obs}

We seek the analogue of Theorem \ref{T:1.2} for the commuting
multivariable case. We extend multivariable power notation \eqref{mnot}
to any $d$-tuple ${\mathbf A} = (A_{1}, \dots, A_{d})$ of commuting
operators on a space $\cX$:
\begin{equation}
\bA^\bn:= A_{1}^{n_{1}}A_{2}^{n_{2}}\ldots A_{d}^{n_{d}}.
\label{2.2a}
\end{equation}
Note the connection between the commutative powers ${\mathbf
A}^{\bn}$ (with $\bn \in {\mathbb Z}^{d}_{+}$)
and the noncommutative powers ${\mathbf A}^{v}$ (with $v \in \cF_{d}$)
in case ${\mathbf A}$ is a commutative operator $d$-tuple:
$$
\bA^{v} = \bA^\bn \; \text{ where } \; \bn= \ba(v), \qquad
\sum_{v \in {\cF}_{d} \colon |v| = N}\bA^{v^{\top}}X\bA^{v} =
            \sum_{\bn \in {\mathbb Z}^{d}_{+} \colon |\bn| = N}
            \frac{N!}{\bn!} {\mathbf A}^{*{\mathbf n}} X {\bA}^{\bn}
$$
for any operator  $X$ on  $\cX$.
In case $(C, {\mathbf A})$ is an output stable pair with ${\mathbf A}$ a
commutative operator $d$-tuple,
the formulas \eqref{2.20}, \eqref{2.11}
and \eqref{2.23} for $\widehat{\cO}^{\ba}_{C, \bA}$,
${\mathcal G}_{C, \bA}$
and ${\mathcal G}^{\ba}_{C, \bA}$ collapse (in view of
\eqref{2.19a}) to
\begin{equation}
\widehat {\mathcal O}^{\ba}_{C,\bA} \colon \; x \mapsto
C(I-Z(\blam) A)^{-1} x=
\sum_{\bn \in {\mathbb Z}^{d}_{+}}\frac{|\bn|!}{\bn!}
\left(C\bA^{\bn}x\right)\, \blam^{\bn}
\label{5.5}
\end{equation}
and
\begin{equation}
{\mathcal G}_{C, \bA}={\mathcal G}^{\ba}_{C, \bA}
=\sum_{\bn\in\bbZ_+^d}\frac{|\bn|!}{\bn!}\bA^{*\bn}C^*C\bA^{\bn}.
\label{5.6}
\end{equation}

We next observe that
a natural commutative counterpart of operators $S_j$ introduced
in \eqref{4.5} are the operators $M_{\lambda_j}$ of
multiplication by the coordinate functions of $\C^d$ for $j=1,\ldots,d$
acting as contractions on the Arveson space $\cH_{\cY}(k_d)$. We will call
the commuting $d$-tuple
${\bf M}_\blam:=(M_{\lambda_1},\ldots,M_{\lambda_d})$
the shift of $\cH_{\cY}(k_d)$, whereas the commuting $d$-tuple
${\bf M}^*_\blam:=(M^*_{\lambda_1},\ldots,M^*_{\lambda_d})$ consisting
of the adjoints of $M_{\lambda_j}$'s (in the metric of $\cH_{\cY}(k_d)$)
will be referred to as to the {\em backward shift}. Recall that monomials
$\blam^{\bn}$ form an orthogonal basis for $\cH(k_d)$. As we have seen,
\begin{equation}
\langle \blam^{\bn}, \, \blam^{\bm}\rangle_{\cH(k_d)}
=\left\{\begin{array}{ccc}
{\displaystyle\frac{\bn!}{|\bn|!}} & \mbox{if} & \bn=\bm\\
0 & & \mbox{otherwise}.\end{array}\right.
\label{5.7}
\end{equation}
A simple calculation based on \eqref{5.7} gives
\begin{equation}
M_{\lambda_j}^* \blam^{\bf m}=\frac{m_j}{|{\bf m}|}
\blam^{{\bf m}-e_j} \; \; (m_j\ge 1)\quad\mbox{and}
\quad M_{\lambda_j}^* \blam^{\bf m}=0 \; \; (m_j=0)
\label{5.8}
\end{equation}
where ${\bf m}=(m_1,\ldots,m_d)$ and $e_j\in\bbZ_+^d$ is defined in
\eqref{2.16}. More generally,
\begin{equation}
\left({\bf M}_{\blam}^*\right)^{\bf n} \;  \blam^{\bf m}=
\left\{\begin{array}{ccc} {\displaystyle\frac{{\bf m}! |{\bf m}-{\bf
n}|!}{|{\bf m}|!({\bf m}-{\bf n})!}} \, \blam^{{\bf m}-{\bf n}}, &
\mbox{if} & m_j\ge n_j \; \; \mbox{for} \; j=1,\ldots,d,\\
0,&&\mbox{otherwise},\end{array}\right.
\label{5.9}
\end{equation}
where according to \eqref{2.2a}
$$
\left({\bf M}_{\blam}^*\right)^{\bf n}:=
\left(M_{\lambda_1}^*\right)^{n_1}\left(M_{\lambda_2}^*\right)^{n_2}
\cdots \left(M_{\lambda_d}^*\right)^{n_d}.
$$
The following proposition includes the analogue of Proposition
\ref{P:4.2} for the present commutative setting.

\begin{proposition}  \label{P:5.1}
Let ${\bf M}_{\blam}^*$ be the $d$-tuple of backward shifts on
$\cH_{\cY}(k_d)$ and let $G \colon \cH_{\cY}(k_{d}) \to \cY$ be the operator of
evaluation at $0 \in {\mathbb B}^{d}$
\begin{equation}
G\colon \; f(\blam)\to f(0).
\label{5.10}
\end{equation}
Then:
\begin{enumerate}
           \item
           For every $f \in \cH_{\cY}(k_{d})$ and every $\blam \in
           {\mathbb B}^{d}$ we have
           \begin{equation}
           f(\blam)-f(0)=\sum_{j=1}^d \lambda_j(M_{\lambda_j}^*f)(\blam).
           \label{5.12}
           \end{equation}

           \item The pair $(G,{\bf M}_{\blam}^*)$ is isometric:
\begin{equation}
I-M_{\lambda_1}M_{\lambda_1}^*-\ldots -M_{\lambda_d}M_{\lambda_d}^*=G^*G.
\label{5.11}
\end{equation}

\item
The abelianized observability operator associated with the pair $(G, {\mathbf
M}_{\blam}^{*})$ is the identity operator:
\begin{equation}  \label{G=I}
          \widehat{\cO}^{\ba}_{G, {\bf M}_{\blam}^*} =I_{\cH_{\cY}(k_d)}.
\end{equation}

\item
The $d$-tuple ${\bf M}_{\blam}^*$ is strongly stable,
that is,
\begin{equation}
\label{5.13a}
\lim_{N \to \infty} \sum_{v \in \cF_{d} \colon |v| = N}
\|({\bf M}_{\blam}^*)^v f
\|_{\cH_{\cY}(k_d)}^{2} =
\lim_{N \to \infty} \sum_{\bn \in\bbZ^{d}_{+}
\colon |\bn| = N} \frac{N!}{\bn!} \|({\bf M}_{\blam}^*)^\bn f
\|_{\cH_{\cY}(k_d)}^{2} =0
\end{equation}
for every  $f \in\cH_{\cY}(k_d)$.
\end{enumerate}
\end{proposition}

\begin{proof}[Proof of (1):] One can easily verify the identity \eqref{5.12} on
monomials $y\cdot\blam^{\bf m}$ (with $y\in\cY$ and ${\bf
m}\in\bbZ_+^d$) using \eqref{5.8}.  Then the result follows for all
$f \in \cH^{2}_{\cY}(k_{d})$ by linearity and continuity.

\medskip
\noindent
{\em Proof of (2):} Note that $G^{*} \colon \cY \to \cH_{\cY}(k_{d})$
is the identification of a vector $y \in \cY$ with the
constant function $y \in \cH_{\cY}(k_{d})$.  We then see that
\eqref{5.11} is simply the operator expression of \eqref{5.12}.

\medskip
\noindent
{\em Proof of (3):}
         From \eqref{5.9} and \eqref{5.10} we see that
$$
G\left({\bf M}_{\blam}^*\right)^{\bf n}f=\frac{\bn!}{|\bn|!}f_{\bn}\quad
\mbox{if} \quad f(\blam)=\sum_{\bm\in\bbZ_+^d}f_{\bm}\blam^{\bm} \; \;
\mbox{and} \; \; \bn\in\bbZ_+^d
$$
and therefore, according to definition \eqref{5.5},
$$
\widehat{\cO}^{\ba}_{G, {\bf M}_{\blam}^*}f:=
\sum_{\bn \in {\mathbb Z}^{d}_{+}}\frac{|\bn|!}{\bn!}
\left(G\left({\bf M}_{\blam}^*\right)^{\bf n}f\right)\, \blam^{\bn}
=\sum_{\bn \in {\mathbb Z}^{d}_{+}}f_{\bn}\blam^{\bn}=f(\blam).
$$
Since the latter equality holds for every
$f\in\cH_{\cY}(k_d)$,
\eqref{G=I} follows as asserted.

\medskip
\noindent
{\em Proof of (4):}
This can be derived directly from \eqref{5.9} or via Proposition
\ref{P:2-1.1'converse}
since $\widehat{\cO}^{\ba}_{G, {\bf M}_{\blam}^*}=I$ and therefore,
the pair $(G, {\bf M}_{\blam}^*)$ is exactly observable.
\end{proof}

\begin{remark}  \label{R:5.3}
{\rm Note that in contrast to the noncommutative case
(Proposition \ref{P:4.2}), the operator
$$
R=\begin{bmatrix}M_{\lambda_1}^*\\ \vdots \\ M_{\lambda_d}^* \\
G\end{bmatrix}:\cH_{\cY}(k_d)\to (\cH_{\cY}(k_d))^d\oplus \cY
$$
is not unitary (just isometric). A simple calculation shows that
$$
I-RR^*=\begin{bmatrix}P &0 \\ 0 & 0\end{bmatrix}\colon \;
\begin{bmatrix}(\cH_{\cY}(k_d))^d \\ \cY\end{bmatrix}\to
\begin{bmatrix}(\cH_{\cY}(k_d))^d \\ \cY\end{bmatrix}
$$
where $P$ is the orthogonal projection of $(\cH_{\cY}(k_d))^d$ onto
the subspace}
$$
\left\{h=\begin{bmatrix}h_1 \\ \vdots \\
h_d\end{bmatrix}\in(\cH_{\cY}(k_d))^d \colon \;
\sum_{j=1}^d\lambda_jh_j(\blam)\equiv 0\right\}.
$$
\end{remark}

If a pair $(C,\bA)$ is $\ba$-output stable, then the observability
operator $\widehat{\cO}^{\ba}_{C, \bA}\colon \; \cX\to \cH_{\cY}(k_d)$
is bounded and its range
\begin{equation}
\operatorname{Ran} \widehat{\cO}^{\ba}_{C, \bA}:=\{C(I-Z(\blam)A)^{-1} x: \;
x\in\cX\}
\label{5.1}
\end{equation}
is a linear manifold in $\cH_{\cY}(k_d)$.  We have the following
partial analogues of part (3) of Theorem \ref{T:2-1.2}.

\begin{theorem}  \label{T:3-1.2}
         Let $(C, {\mathbf A})$ be an $\ba$-output stable pair.  Then:
         \begin{enumerate}
             \item
             $\operatorname{Ran} \widehat{\mathcal O}^{\ba}_{C,
             {\mathbf A}}$ with the lifted norm
             \begin{equation}
             \left\|C(I-Z(\blam)A)^{-1}x\right\|_{\cH(K^{\ba}_{C,\bA})}=
             \|Q^{\ba}x\|_{\cX}
             \label{5.3}
             \end{equation}
             where $Q^{\ba}$ is the orthogonal projection of $\cX$ onto
             $(\operatorname{Ker}{\mathcal G}^{\ba}_{C,\bA})^\perp$
is isometrically equal to the reproducing kernel Hilbert space
$\cH(K^{\ba}_{C, {\mathbf A}})$ with reproducing kernel
$K^{\ba}_{C, {\mathbf A}}(\blam, \bzeta)$ given by
$$
K^{\ba}_{C,\bA}(\blam ,\bzeta)
=C(I-Z(\blam)A)^{-1}(I-A^*Z(\bzeta)^*)^{-1}C^*\quad(\blam,\bzeta\in\B^d).
$$
\item
$\operatorname{\overline{Ran}} \widehat{\mathcal O}^{\ba}_{C, {\mathbf
A}}$ with norm inherited form $\cH_{\cY}(k_{d})$ is a reproducing kernel
Hilbert space $\cH(K^{\ba,-}_{C,{\mathbf A}})$ with
reproducing kernel  $K^{\ba,-}_{C,{\mathbf A}}(\blam, \zeta)$ given by
$$
K^{\ba,-}(\blam, \bzeta) =
C(I-Z(\blam)A)^{-1} ({\mathcal G}^{\ba}_{C, {\mathbf A}})^{-1}
(I-A^*Z(\bzeta)^*)^{-1}C^*\quad(\blam,\bzeta\in\B^d).
$$
\end{enumerate}
\end{theorem}

We next discuss separately the case where ${\mathbf A}$ is
$C$-abelian and then the general case.

\subsubsection{$\cH(K^{\ba}_{C, {\mathbf A}})$ for the case where $\bA$
is $C$-abelian} \label{S:C-Obs-comA}

In case $(C,{\mathbf A})$ is an $\ba$-output-stable pair with
${\mathbf A}$ $C$-abelian, then we have the following commutative
analogue of Theorem \ref{T:2-1.2'}.

\begin{theorem}  \label{T:3-1.2c}
           Let $(C, {\mathbf A})$ be a contractive $\ba$-output-stable
pair such that
           operator $d$-tuple ${\mathbf A}$ is $C$-abelian.  Then:
         \begin{enumerate}
             \item
             The intertwining relations
             \begin{equation}
             M_{\blam_j}^*\widehat{\cO}^{\ba}_{C,
\bA}=\widehat{\cO}^{\ba}_{C, \bA}A_j
             \quad\text{for} \; \; j = 1, \dots, d
             \label{5.13}
             \end{equation}
             hold, and hence the linear submanifold $\operatorname{Ran}
             \widehat{\cO}^{\ba}_{C, \bA}$ of $\cH_{\cY}(k_d)$ is ${\bf
             M}_{\blam}^*$-invari\-ant.
             \item
             The operator $\widehat {\mathcal O}^{\ba}_{C, {\mathbf A}}$
             maps ${\mathcal X}$ contractively into $\cH_{\cY}(k_{d})$.
             This mapping is isometric if and only if $(C, {\mathbf A})$ is
             isometric and ${\mathbf A}$ is strongly stable.
             \item If ${\mathcal M}: =
             \operatorname{Ran}{\mathcal O}^{\ba}_{C, {\mathbf A}}$
             is given the lifted norm \eqref{5.3} (so ${\mathcal M}$ is
             isometrically equal to $\cH(K^{\ba}_{C, {\mathbf A}})$ by
             Theorem \ref{T:3-1.2} (1)), then the difference-quotient
             inequality
$$
             \sum_{j=1}^d \|M_{\lambda_j}^*f\|^2_{\cH(K^{\ba}_{C,\bA})}
             \le\|f\|^2_{\cH(K^{\ba}_{C,\bA})}-\|f(0)\|^2_{\cY}
$$
             holds for every $f\in\cH(K^{\ba}_{C,\bA})$. Moreover, the
             difference-quotient identity
$$
             \sum_{j=1}^d \|M_{\lambda_j}^*f\|^2_{\cH(K^{\ba}_{C,\bA})}
             =\|f\|^2_{\cH(K^{\ba}_{C,\bA})}-\|f(0)\|^2_{\cY}
$$
             holds for every $f\in\cH(K^{\ba}_{C,\bA})$ if and only if
the subspace
             $(\operatorname{Ker} {\mathcal G}_{C,\bA})^\perp$
             is $\bA$-invariant and the restriction $(C^{0},\bA^{0})$ 
(defined in
             \eqref{matrix-decom}) of $(C,\bA)$ to the subspace
             $(\operatorname{Ker} {\mathcal G}_{C,\bA})^\perp$ is isometric.
        \end{enumerate}
        \end{theorem}

        \begin{proof}
            By \eqref{5.5} and \eqref{5.8}, we have for every $x\in\cX$,
           \begin{eqnarray*}
           (M_{\lambda_{j}})^{*}\widehat{\cO}^{\ba}_{C, \bA} x &=&
           (M_{\lambda_{j}})^{*} \left( \sum_{\bn \in {\bbZ}^{d}_{+}}
           \frac{|\bn|!}{\bn!} C\bA^{\bn}x \cdot \blam^{\bn} \right) \\
           & =& \sum_{\bn \in \bbZ^{d}_{+}}
           \frac{|\bn|!}{\bn!} \frac{ n_{j}}{|\bn|} C\bA^{\bn}x
           \cdot\blam^{\bn-e_{j}}  \\
           &=& \sum_{\bn \in\bbZ^{d}_{+}} \frac{|\bn-e_{j}|!}{(\bn - e_{j})!}
           C\bA^{\bn-e_{j}}\cdot \blam^{\bn - e_{j}} A_{j}x \\
           & =& \left(\sum_{\bn \in\bbZ^{d}_{+}}\frac{|\bn|!}{\bn!} C\bA^{\bn}
           \blam^{\bn} \right) \cdot A_{j}x \\
           & =&\widehat{\cO}^{\ba}_{C, \bA} \cdot A_{j}x
           \end{eqnarray*}
           and \eqref{5.13} follows. This completes the proof of statement
           (1) in the theorem.

           Since the pair $(C,\bA)$ is contractive and $\bA$ is
           $C$-abelian, we have
           $$
           {\mathcal G}^{\ba}_{C,\bA}={\mathcal G}_{C,\bA}\le Q=Q^\ba\le I.
           $$
           Therefore,
           \begin{equation}
           \|\widehat{\cO}^{\ba}_{C,\bA}x\|_{\cH_{\cY}(k_d)}=
           \langle {\mathcal G}^{\ba}_{C,\bA}x, \, x\rangle_{\cX}^{\frac{1}{2}}
           \le \|Q^\ba x\|_{\cX}= 
\|\widehat{\cO}^{\ba}_{C,\bA}x\|_{\cH(K^{\ba}_{C,\bA})}\le\|x\|_{\cX}.
           \label{5.16}
           \end{equation}
           Now the arguments used in the proof of Theorem \ref{T:2-1.2'} can
           be used to prove the remaining statements in the theorem.
        \end{proof}

        For the converse direction we have the following result.

        \begin{theorem}  \label{T:3c-converse}
        Let ${\mathcal M}$ be a Hilbert space of $\cY$-valued functions included
        into $\cH_{\cY}(k_d)$ and let us assume that ${\mathcal M}$ is
        ${\bf M}_{\blam}^*$-invariant.
        \begin{enumerate}
        \item If the inequality
        \begin{equation}
        \sum_{j=1}^d \|M_{\lambda_j}^*f\|^2_{{\mathcal M}}
        \le\|f\|^2_{{\mathcal M}}-\|f(0)\|^2_{\cY}
        \label{5.17}
        \end{equation}
        holds for every $f\in{\mathcal M}$, then
        ${\mathcal M}=\operatorname{Ran}(\widehat{\cO}^{\ba}_{C,\bA})$ for
        a contractive
        and exactly observable (with respect to ${\mathcal M}$) pair $(C,\bA)$
        with the commutative $d$-tuple $\bA=(A_1,\ldots,A_d)$.
        In particular,
        ${\mathcal M}$ is contractively included in $\cH_{\cY}(k_d)$.
        \item If the equality
        \begin{equation}
        \sum_{j=1}^d \|M_{\lambda_j}^*f\|^2_{{\mathcal M}}
        =\|f\|^2_{{\mathcal M}}-\|f(0)\|^2_{\cY}
        \label{5.18}
        \end{equation}
        holds for every $f\in{\mathcal M}$, then ${\mathcal M}=
        \operatorname{Ran}
        (\widehat{\cO}^{\ba}_{C,\bA})$ for an isometric and  exactly
        observable (with respect to ${\mathcal M}$) pair $(C,\bA)$
        with the commutative $d$-tuple $\bA$. By part $(1)$, ${\mathcal M}$ is
        contractively included in $\cH_{\cY}(k_d)$. Moreover, it is 
isometrically
        included in $\cH_{\cY}(k_d)$ if and only if the restriction of the
        backward shift ${\bf M}_{\blam}^*$ to ${\mathcal M}$ is strongly
stable, i.e.,
        \begin{equation}
        \label{5.19}
        \lim_{N \to \infty} \sum_{\bn \in\bbZ^{d}_{+}
        \colon |\bn| = N} \frac{N!}{\bn!} \|({\bf M}_{\blam}^*)^\bn f
        \|_{{\mathcal M}}^{2} =0\quad\mbox{for every} \; \; f \in{\mathcal M}.
        \end{equation}
        \end{enumerate}
        \label{T:4.16}
        \end{theorem}

        \begin{proof} Define operators
$C\colon \, {\mathcal M}\to \cY$ and $A_1,\ldots,A_d\colon \, {\mathcal
M}\to {\mathcal M}$ by
\begin{equation}
C=G\vert_{{\mathcal M}}\colon \; f\to f(0)\quad\mbox{and}\quad
A_j=M_{\lambda_j}^*\vert_{{\mathcal M}} \; \; (j=1,\ldots,d).
\label{5.14}
\end{equation}
        Thus, the $d$-tuple $\bA$ is the
        restriction of the backward-shift tuple ${\bf M}_{\blam}^*$ to 
${\mathcal
        M}$. By part (3) of Proposition \ref{P:5.1}, it follows that
        $\widehat{\cO}^{\ba}_{C, \bA}=I_{\mathcal M}$ and thus, the pair
        $(C,\bA)$ is exactly observable (with respect to ${\mathcal M}$) and
        the range of the associated observability operator
        $\widehat{\cO}^{\ba}_{C, \bA}$ coincides (algebraically) with ${\mathcal
        M}$. Now we write \eqref{5.17} in terms of operators \eqref{5.14}
        as
        $$
        \sum_{j=1}^d \|A_jf\|^2_{{\mathcal M}}
        \le\|f\|^2_{{\mathcal M}}-\|Cf\|^2_{\cY}\qquad (f\in{\mathcal M})
        $$
        and conclude that the pair $(C,\bA)$ is contractive. Similarly,
        assumption \eqref{5.18} means that the chosen pair $(C,\bA)$ is
isometric.
        Furthermore, if ${\mathcal M}$ is included in $\cH_{\cY}(k_d)$
        isometrically, relation \eqref{5.19} holds since ${\bf M}_{\blam}^*$
        is strongly stable (see part (4) of Proposition \ref{P:5.1}).
        Conversely, if \eqref{5.19}
        holds, that is, if the commutative $d$-tuple $\bA=(A_1,\ldots,A_d)$
        is strongly stable on ${\mathcal M}$, the Stein equation \eqref{3.4}
        has a unique positive semidefinite solution. Since the pair $(C,\bA)$ is
        isometric (recall that we are proving isometrical inclusion under
        assumption \eqref{5.18}), this unique solution is the identity operator.
        On the other hand the observability gramian  ${\mathcal G}^{\ba}_{C,\bA}
        ={\mathcal G}_{C,\bA}$ defined by the convergent series \eqref{5.6}
        satisfies the same Stein equation (as observed in part (2) of
        Theorem \ref{T:2-1.1}). Thus,
        ${\mathcal G}^{\ba}_{C,\bA}={\mathcal G}_{C,\bA}=I$.
        Note that the inequality \eqref{3.6} holds with $H = {\mathcal G}_{C,
        \bA} = I$, i.e.,
        $$ \sum_{ v \in \cF_{d} \colon |v| < N} \bA^{* v^{\top}} C^{*}C
        A^{v} \le I - \sum_{v \in \cF_{d} \colon |v| = N+1} \bA^{*
        v^{\top}} \bA^{v}.
        $$
       Taking strong limits as $N \to \infty$ and noting that $I =
       {\mathcal G}_{C, \bA} =
        \sum_{v \in \cF_{d}} \bA^{*v ^{\top}} C^{*} C \bA^{v}$ then
        gives
        $$ I \le I - \operatorname{s-lim}_{N \to \infty} \sum_{v \in
\cF_{d} \colon |v|
        = N} \bA^{*v^{\top}} \bA^{v}
        $$
        from which the strong-stability of $\bA$ follows.
        Then ${\mathcal M}=\operatorname{Ran}(\widehat{\cO}^{\ba}_{C,
        \bA})$ is isometrically included in $\cH_{\cY}(k_d)$ by statement (2) in
        Theorem \ref{T:3-1.2c}.
\end{proof}

We have the following analogue of Theorem \ref{T:unique} for the
present commutative situation.

\begin{theorem}  \label{T:c-unique}
       Suppose that  $(C, \bA)$
       and $(\widetilde C,  \widetilde \bA)$ are two observable
       output-stable pairs with both $\bA$ and $\widetilde \bA$ commutative
       such that $K^{\ba}_{C, \bA}(\blam, \bzeta) =
       K^{\ba}_{\widetilde C, \widetilde \bA}(\blam, \bzeta)$ for all
       $\blam, \bzeta \in {\mathbb B}^{d}$. Then there is a unitary
       operator $U \colon {\mathcal X} \to \widetilde {\mathcal X}$ such that
       \begin{equation}  \label{Uintertwine}
         C = \widetilde C U \quad\text{and}\quad A_{j} = U^{-1}
\widetilde A_{j} U\quad
         \text{for} \; \; j = 1, \dots, d.
       \end{equation}
       \end{theorem}

       \begin{proof}
           Suppose that $(C, \bA)$ and $(\widetilde C, \widetilde \bA)$ are
             as in the hypothesis of the theorem.  The identity of the
             kernels $K^{\ba}_{C, \bA}$ and $K^{\ba}_{\widetilde C, \widetilde
             \bA}$ implies equality of the respective coefficients of
$\blam^{\bn}
             \bzeta^{{\mathbf m}}$ for each $\bn, {\mathbf m} \in {\mathbb
             Z}^{d}_{+}$:
             $$ \frac{|\bn|!}{\bn !} \frac{ |{\mathbf m}|!}{ {\mathbf m}!}
C \bA^{\bn}
             \bA^{*{\mathbf m}} C^{*} =
             \frac{|\bn|!}{\bn !} \frac{ |{\mathbf m}|!}{ {\mathbf m}!}
\widetilde C
             \widetilde \bA^{\bn} \widetilde \bA^{*{\mathbf m}}
             \widetilde C^{*}.
             $$
             If we define a mapping $U$ by
	 \begin{equation} \label{defineU}
	  U \colon \bA^{* \bn } C^{*} y \mapsto \widetilde \bA^{* \bn }
	 \widetilde C^{*} y
	\end{equation}
             it follows that $U$ extends by linearity to an isometry from
             $$ {\mathcal D}: = \operatorname{span} \{\bA^{*{\mathbf m}} C^{*} y
             \colon {\mathbf m} \in {\mathbb Z}^{d}_{+} \text{ and } y \in \cY\}
             $$
             onto
             $$ {\mathcal R}: = \operatorname{span} \{\widetilde
\bA^{*{\mathbf m}}
             \widetilde C^{*} y
             \colon {\mathbf m} \in {\mathbb Z}^{d}_{+} \text{ and } y \in \cY\}
             $$
             Since both $(C, \bA)$ and $(\widetilde C, \widetilde \bA)$ are
             observable, we see that ${\mathcal D}$ is dense in
${\mathcal X}$ and
             that ${\mathcal R}$ is dense in $\widetilde {\mathcal X}$.
Hence $U$
             extends to a unitary operator from ${\mathcal X}$ onto $\widetilde
             {\mathcal X}$ by continuity.  From the defining equations
             \eqref{defineU} for $U$ we see that
             $$ U C^{*} = \widetilde C^{*}\quad \text{and}\quad U A_{j}^{*} =
             \widetilde A_{j}^{*} U.
             $$
             By taking adjoints and using that $U$ is unitary, we arrive at the
             intertwining equations \eqref{Uintertwine} as wanted.
        \end{proof}

Theorem \ref{T:c-unique} can be adapted to give the following result
concerning containment between two backward-shift-invariant subspaces
rather than equality; the finite-dimensional case appears as
Proposition 1.2 in \cite{BR1}.

\begin{theorem}  \label{T:containment}
       Let ${\mathcal M}$ and $\widetilde{\mathcal M}$ be two
       backward-shift-invariant subspaces of the Arveson space
       $\cH_{\cY}(k_d)$ with realizations
          \begin{equation}
          {\mathcal M}=\operatorname{Ran} \widehat{\mathcal O}_{C,{\bf A}}\quad
          \mbox{and}\quad \widetilde{\mathcal M}=\operatorname{Ran}
          \widehat{\mathcal
          O}_{\widetilde{C},\widetilde{\bf A}}
          \label{7.1}
          \end{equation}
          where the $d$-tuples ${\bf A}=(A_1,\ldots,A_d)\in\cX^d$ and
          $\widetilde{\bf A}=
          (\widetilde{A}_1,\ldots,\widetilde{A}_d)\in\widetilde{\cX}^d$
          are commutative and strongly stable and the pairs
          $(C,{\bf A})$ and $(\widetilde{C},\widetilde{\bf A})$ are isometric.
          Then ${\mathcal M}\subseteq\widetilde{\mathcal M}$ if and 
only if there
          exists an isometry $V: \, \cX\to \widetilde{\cX}$ such that
          \begin{equation}
          C=V\widetilde{C}\quad\mbox{and}\quad VA_j=\widetilde{A}_jV \; \;
          (j=1,\ldots,d).
          \label{7.2}
          \end{equation}
\end{theorem}

\begin{proof}
          The necessity part is clear. For the sufficiency part, assume that
          ${\mathcal M}\subseteq\widetilde{\mathcal M}$. By Theorem
\ref{T:3-1.2c},
          there exist unitary operators $U: \, {\mathcal M}\to \cX$
          and $\widetilde{U}: \, \widetilde{\mathcal M}\to\widetilde{\cX}$
          such that
          $$
          U^*A_jU=M_{\lambda_j}^*\vert_{\mathcal M},\quad
          \widetilde{U}^*\widetilde{A}_j\widetilde{U}=
          M_{\lambda_j}^*\vert_{\widetilde{\mathcal M}}\quad (j=1,\ldots,d).
          $$
          and
          $$
          CU=G\vert_{\mathcal M},\quad
          \widetilde{C}\widetilde{U}=G\vert_{\widetilde{\mathcal M}}
          $$
          where the operator $G: \, \cH_{\cY}\to\cY$ is defined in \eqref{5.10}.
          Let ${\mathcal I} \colon  {\mathcal M}\to \widetilde{\mathcal M}$
          be the inclusion operator. Clearly ${\mathcal I}$ is
isometric. Then the
          operator $V=U^*{\mathcal I}\widetilde{U} \colon \cX\to 
\widetilde{\cX}$
          is isometric and satisfies \eqref{7.2}.
        \end{proof}

        \subsubsection{$\cH(K^{\ba}_{C, {\mathbf A}})$: The general case}
        \label{S:C-Obs-noncomA}

        In case the $\ba$-output-stable pair $(C, {\mathbf A})$ is such that
        ${\mathbf A}$ is not $C$-abelian, it can happen that the
        associated reproducing kernel Hilbert space is not invariant under
        the backward-shift tuple ${\mathbf M}_{\blam}^*$, as the following
        example shows.

        \begin{example} \label{E:not-shift*-inv}
        {\rm Let
$$C=\begin{bmatrix}\frac{\sqrt{3}}{2} & 0 \\ 0 & \frac{\sqrt{3}}{2}
\end{bmatrix},\quad
A_1=\begin{bmatrix} 0 & 0 \\  \frac{1}{2} & 0 \end{bmatrix},\quad
A_2=\begin{bmatrix} 0 & \frac{1}{2}\\  0 & 0\end{bmatrix}.
$$
Then a straightforward calculation gives
\begin{eqnarray}
K^{\ba}_{C, {\mathbf
A}}(\blam,\bzeta)&=&C(I-Z(\blam)A)^{-1}(I-A^*Z(\bzeta)^*)^{-1}C^*\nonumber\\
&=&\frac{3}{(4-\lambda_1\lambda_2)(4-\overline{\zeta}_1\overline{\zeta}_2)}
\left[\begin{array}{cc}
2 & \lambda_2 \\ \lambda_1 & 2\end{array}\right]\left[\begin{array}{cc}
2 & \overline{\zeta}_1 \\ \overline{\zeta}_2 & 2\end{array}\right].\nonumber
\end{eqnarray}
Thus $K^{\ba}_{C, {\mathbf A}}(\lambda,w)$ is positive definite on
$\B^2\times\B^2$ and the space $\cH(K^{\ba}_{C, {\mathbf A}})$ is
spanned by the two rational functions
$$
f_1(\lambda)=\frac{4}{4-\lambda_1\lambda_2}\left[\begin{array}{c}2\\
\lambda_1\end{array}\right]
\quad\mbox{and}\quad
f_2(\lambda)=\frac{4}{4-\lambda_1\lambda_2}\left[\begin{array}{c}\lambda_2\\
2\end{array}\right].
$$
Furthermore, since
$$
M_{\lambda_1}^*(\lambda_1^{n_1}\lambda_2^{n_2})=\frac{n_1}{n_1+n_2}
\lambda_1^{n_1-1}\lambda_2^{n_2},
$$
and since
$$
\frac{4\lambda_1}{4-\lambda_1\lambda_2}=\sum_{j=0}^\infty
\frac{\lambda_1^{j+1}\lambda_2^j}{4^j},
$$
it holds that
$$
M_{\lambda_1}^*\left(\frac{4\lambda_1}{4-\lambda_1\lambda_2}\right)=
\sum_{j=0}^\infty\frac{j+1}{2j+1}\left(\frac{\lambda_1\lambda_2}{4}\right)^{j}.
$$
The latter function is rational if and only if the single-variable
function $F(z) = \sum_{j=0}^{\infty} \frac{j+1}{2j+1} z^{j}$ is
rational.  By the well-known Kronecker theorem, $F$ in turn is
rational if and only if the associated infinite Hankel matrix
$$ {\mathbb H} = [s_{i+j}]_{i,j=0}^{\infty} \quad\text{where}\quad
     s_{k} = \frac{k+1}{2k+1}
$$
has finite rank.
However one can check that the finite Hankel matrices ${\mathbb H}_{n}
= [s_{i+j}]_{i,j=0}^{n}$ have full rank for all $n=0,1,2,\dots$ and
hence $F(z)$ is not rational.
Therefore, $M_{\lambda_1}^*f_1$ does not belong to $\cH(K^{\ba}_{C,
{\mathbf A}})$ and
hence $\cH(K^{\ba}_{C, {\mathbf A}})$ is not invariant  under
$M_{\lambda_1}^*$.}
       \end{example}

        For the general case, there is a simple replacement for
        ${\mathbf M}_{\blam}^{*}|_{\cH(K^{\ba}_{C, {\mathbf
        A}})}$.  Specifically, given an $\ba$-output-stable pair $(C,
        {\mathbf A})$, we define an operator-tuple ${\mathbf T} = (T_{1},
        \dots, T_{d})$ on $\operatorname{Ran} \widehat{\cO}^{\ba}_{C,{\mathbf
A}}$ by
        \begin{equation} \label{6.10a}
        T_j \widehat{\cO}^{\ba}_{C, \bA}x=\widehat{\cO}^{\ba}_{C,
        \bA}A_jx\quad
        \text{for} \; \; x \in(\operatorname{Ker}{\mathcal
G}^{\ba}_{C,\bA})^\perp
        \; \text{ and } \; j=1,\dots,d.
        \end{equation}
        We then have
        \begin{eqnarray}
        f(\blam)-f(0)&=&C(I-Z(\blam)A)^{-1}x-Cx\nonumber\\
        &=&C(I-Z(\blam)A)^{-1}Z(\blam)Ax\nonumber\\
        &=&\sum_{j=1}^d \lambda_j C(I-Z(\blam)A)^{-1}A_jx\nonumber\\
        &=&\sum_{j=1}^d  \lambda_j \cdot(\widehat{\cO}^{\ba}_{C,
        \bA}A_jx)(\blam)\nonumber\\
        &=&\sum_{j=1}^d  \lambda_j\cdot
        (T_j\widehat{\cO}^{\ba}_{C, \bA}x)(\blam)
        =\sum_{j=1}^d  \lambda_j\cdot (T_j f)(\blam).
        \label{Gleason}
        \end{eqnarray}
        We next give the following analogue of Theorem \ref{T:3-1.2c} for
        the general case.

        \begin{theorem}  \label{T:3-1.2nc}
            Let $(C,\bA)$ be a contractive pair with $C\in{\mathcal
L}(\cX,\cY)$ and
            $\bA=(A_1,\ldots,A_d)\in{\mathcal L}(\cX)^d$. Then:
            \begin{enumerate}
            \item The $Z$-transformed observability operator
            $\widehat{\cO}^{\ba}_{C,\bA}$
            is a contraction of $\cX$ into the reproducing kernel Hilbert
            space $\cH(K^\ba_{C,\bA})$.
            It is an isometry if and only if the the pair $(C,\bA)$ is
            $\ba$-observable.

            \item The space $\cH(K^{\ba}_{C,\bA})$ is contractively
included in the
            Arveson space $\cH_{\cY}(k_d)$; it is isometrically included in
            $\cH_{\cY}(k_d)$ if and only if $\widehat{\cO}^{\ba}_{C,\bA}$ (as an
            operator from $\cX$ into $\cH_{\cY}(k_d)$) is a partial isometry.
            \item For every function $f\in\cH(K^{\ba}_{C,\bA})$ it holds that
            \begin{equation}
            f(\blam)-f(0)=\sum_{j=1}^d \lambda_j(T_jf)(\blam)\quad(\blam
            \in\B^d)
            \label{6.16a}
            \end{equation}
            and
            \begin{equation}
            \sum_{j=1}^d \|T_jf\|^2_{\cH(K^\ba_{C,\bA})}
            \le\|f\|^2_{\cH(K^\ba_{C,\bA})}-\|f(0)\|^2_{\cY}
            \label{6.17a}
            \end{equation}
            where $T_1,\ldots,T_d\in{\mathcal L}(\cH(K^\ba_{C,\bA}))$ are the
            operators defined in \eqref{6.10a}.

            \item Equality holds in \eqref{6.17a} for every
            $f\in\cH(K^\ba_{C,\bA})$   if and only if the subspace
            $({\rm Ker} \, {\mathcal G}_{C,\bA})^\perp$ is
            $\bA$-invariant and the restriction $(C^{0},\bA^{0})$ (defined in
            \eqref{matrix-decom}) of $(C,\bA)$ to the subspace
            $(\operatorname{Ker} {\mathcal G}_{C,\bA})^\perp$ is isometric.

            \item If $\cH(K^{\ba}_{C, {\mathbf A}})$ is isometrically 
included in
            $\cH_{\cY}(k_d)$, then $T_j=M^*_{\lambda_j}$ for $j=1,\ldots,d$
            and therefore, $\cH_{\cY}(k_d)$ is ${\bf M}^*_{\blam}$-invariant.
            \end{enumerate}
           \end{theorem}

           \begin{proof} Since the pair $(C,\bA)$ is contractive, the
           identity operator $H = I_{{\mathcal X}}$ solves the Stein
           inequality \eqref{3.4a}.  Then
           ${\mathcal G}^{\ba}_{C,\bA}\le {\mathcal G}_{C,\bA}\le I_{\cX}$
           (by part (1) of Proposition \ref{P:Gineq} and  part (2) of Theorem
           \ref{T:2-1.1}). Thus,
           $$
           {\mathcal G}^{\ba}_{C,\bA}\le Q^{\ba}\le I_{\cX}
           $$
           where $Q^{\ba}$ is the orthogonal projection of $\cX$ onto
           $({\rm Ker} \,{\mathcal G}^{\ba}_{C,\bA})^\perp$. Therefore
           it holds for every $x\in\cX$ that
           \begin{equation}
           \|\widehat{\cO}^{\ba}_{C,\bA}x\|_{\cH_{\cY}(k_d)}=
           \langle {\mathcal G}^{\ba}_{C,\bA}x, \, x\rangle_{\cX}^{\frac{1}{2}}
           \le \|Q^{\ba} x\|_{\cX}=
 \|\widehat{\cO}^{\ba}_{C,\bA}x\|_{\cH(K^{\ba}_{C,\bA})}\le\|x\|_{\cX}.
           \label{5.100}
           \end{equation}
           We have the equality instead of the first inequality in \eqref{5.100}
           if and only if  ${\mathcal G}^{\ba}_{C,\bA}= Q^{\ba}$, that 
is, if and
           only if $\widehat{\cO}^{\ba}_{C,\bA}$ is a partial isometry.
Furthermore,
           the second inequality in \eqref{5.100} can be replaced by equality
           if and only if $Q^{\ba}=I_{\cX}$, i.e., if and only if
           the pair $(C,\bA)$ is $\ba$-observable. This completes the proof
           of the two first assertions in the theorem. The multivariable
           difference-quotient relation \eqref{6.16a} follows by the calculation
           \eqref{Gleason}. Furthermore, for every $x\in({\rm Ker} \,{\mathcal
           G}^{\ba}_{C,\bA})^\perp$,
           $$
           \|\widehat{\cO}^{\ba}_{C,\bA}x\|_{\cH(K^{\ba}_{C,\bA})}=
           \|Q^{\ba}x\|_{\cX}=\|x\|_{\cX}
           $$
           and thus, $\widehat{\cO}^{\ba}_{C,\bA}$ maps unitarily
           $(\operatorname{Ker}{\mathcal G}^{\ba}_{C,\bA})^\perp$ onto
           $\cH(K^{\ba}_{C,\bA})$.
           Therefore, by \eqref{6.10a}, $T_j$ is unitarily equivalent to the
           compression of $A_j$ to $({\rm Ker}
           \,{\mathcal G}^{\ba}_{C,\bA})^\perp$ and hence
           $$
           \|T_j\|\le \|A_j\|\quad\mbox{for $j=1,\ldots,d$}.
           $$
           In particular, $T_j\in{\mathcal L}(\cH(K^{\ba}_{C,\bA}))$. For
an element
           $f=\widehat{\cO}^{\ba}_{C,\bA}x\in\cH(K^{\ba}_{C,\bA})$, we have
           \begin{align*}
           \sum_{j=1}^{d} \| T_{j}f\|^{2}_{\cH(K^{\ba}_{C,\bA})}&=
           \sum_{j=1}^{d} \| T_{j} \widehat{\mathcal O}^{\ba}_{C,\bA}x
           \|^{2}_{\cH(K^{\ba}_{C,\bA})}\\
             &= \sum_{j=1}^{d}  \|\widehat {\mathcal
           O}^{\ba}_{C, \bA} A_{j} x \|^{2}_{\cH(K^{\ba}_{C,\bA})} \\
           &=\sum_{j=1}^{d}  \|Q^{\ba} A_{j} x \|^{2}_{\cH(K^{\ba}_{C,\bA})}
\\
           &  \le \sum_{j=1}^{d} \| A_{j}x \|^{2}_{\cX} \\
              & \le \|x\|^2_{\cX} - \| C x \|^{2}_{\cY}
           =\|f\|^2_{\cH(K^{\ba}_{C,\bA})}-\|f(0)\|^2_{\cY}
           \end{align*}
           where the first inequality holds since $Q^{\ba}\le I$ and the
second since $(C, \bA)$ is a contractive pair. This proves inequality
\eqref{6.17a} and it is readily seen that equalities hold throughout in
the last calculation for every $x\in(\operatorname{Ker}{\mathcal
G}^{\ba}_{C,\bA})^\perp$ if and only
           the subspace $(\operatorname{Ker} {\mathcal
G}^{\ba}_{C,\bA})^\perp$ is
           $\bA$-invariant and the restriction $(C^{0},\bA^{0})$ (defined in
           \eqref{matrix-decom}) of $(C,\bA)$ to the subspace
           $(\operatorname{Ker} {\mathcal G}^{\ba}_{C,\bA})^\perp$ is
isometric.

           Finally suppose that $\cH(K^{\ba}_{C, {\mathbf A}})$ is included
           isometrically in $\cH_{\cY}(k_{d})$. Then the assumption
           \eqref{6.17a} becomes
           \begin{equation}
           \sum_{j=1}^d \|T_jf\|^2_{\cH_{\cY}(k_d)}\le\|f\|^2_{\cH_{\cY}(k_d)}
           -\|f(0)\|^2_{\cY} \text{ for every }  f\in{\mathcal M}.
           \label{6.12a}
           \end{equation}
           Then we take the inner product
           of both parts in equality \eqref{Gleason}
           with $f$:
           $$
           \langle f-f(0), \, 
f\rangle_{\cH_{\cY}(k_d)}=\|f\|_{\cH_{\cY}(k_d)}^2-
           \|f(0)\|_{\cY}^2
           $$
           and
           $$
           \sum_{j=1}^d \langle M_{\lambda_j}T_jf, \, f\rangle_{\cH_{\cY}(k_d)}
           =\sum_{j=1}^d \langle T_jf, \,
M_{\lambda_j}^*f\rangle_{\cH_{\cY}(k_d)}.
           $$
           Thus,
           \begin{equation}
           \|f\|_{\cH_{\cY}(k_d)}^2-\|f(0)\|_{\cY}^2=\sum_{j=1}^d
\langle T_jf, \,
            M_{\lambda_j}^*f\rangle_{\cH_{\cY}(k_d)}.
           \label{6.13b}
           \end{equation}
           For any $f$ in $\cH_{\cY}(k_{d})$, applying the identity
           \eqref{5.11} to $f$ and then taking the inner product with $f$
           gives us
           \begin{equation}
           \|f\|^2_{\cH_{\cY}(k_d)}-\|f(0)\|^2_{\cY}=
           \sum_{j=1}^d \langle M_{\lambda_j}M_{\lambda_j}^*f, \,
           f\rangle_{\cH_{\cY}(k_d)}=\sum_{j=1}^d
           \|M_{\lambda_j}^*f\|^2_{\cH_{\cY}(k_d)}.
           \label{5.12a}
           \end{equation}
           Now we conclude from \eqref{5.12a}, \eqref{6.13b} and \eqref{6.12a}
           that
           \begin{eqnarray*}
           \|f\|_{\cH_{\cY}(k_d)}^2-\|f(0)\|_{\cY}^2&=&
           \sum_{j=1}^d \|M_{\lambda_j}^*f\|^2_{\cH_{\cY}(k_d)}\\
           &=&\sum_{j=1}^d \langle T_jf, \,
           M_{\lambda_j}^*f\rangle_{\cH_{\cY}(k_d)}
           \ge \sum_{j=1}^d \|T_jf\|^2_{\cH_{\cY}(k_d)},
           \end{eqnarray*}
from which we get
   \begin{eqnarray*}
0&=&\sum_{j=1}^d \|M_{\lambda_j}^*f\|^2_{\cH_{\cY}(k_d)}-
\sum_{j=1}^d \langle T_jf, \,
           M_{\lambda_j}^*f\rangle_{\cH_{\cY}(k_d)},\\
0&\ge& -\sum_{j=1}^d \langle T_jf, \,
           M_{\lambda_j}^*f\rangle_{\cH_{\cY}(k_d)}
           +\sum_{j=1}^d \|T_jf\|^2_{\cH_{\cY}(k_d)}.
\end{eqnarray*}
Adding these inequalities and using that ${\displaystyle\sum_{j=1}^d}
\langle T_jf, \, M_{\lambda_j}^*f\rangle$ is real
then gives
           \begin{eqnarray*}
           0 &\ge& \sum_{j=1}^d \left(\|M_{\lambda_j}^*f\|^2-\langle T_jf, \,
           M_{\lambda_j}^*f\rangle-\langle M_{\lambda_j}^*f, \, T_jf\rangle
           +\|T_jf\|^2\right)\\
           &=&\sum_{j=1}^d \|M_{\lambda_j}^*f-T_jf\|_{\cH_{\cY}(k_d)}^2.
           \end{eqnarray*}
           Therefore, $M_{\lambda_j}^*f=T_jf$ for $j=1,\ldots,d$ and for every
           $f\in{\mathcal M}$ as asserted.
           This completes the proof of Theorem \ref{T:3-1.2nc}.
           \end{proof}

        \subsection{The Gleason problem: a uniqueness result}
\label{S:gleason}

Let ${\mathcal M}$ be a Hilbert space of $\cY$-valued functions. A tuple
${\bf T}=(T_1,\ldots,T_d)$ of operators $T_j\in{\mathcal M}$ is called
a solution of the  {\em Gleason problem}  (see \cite{gleason, Henkin})
if relation \eqref{6.16a} holds
for every $f\in{\mathcal M}$. Let us say that ${\bf T}$ {\em is a
contractive solution of the Gleason problem} if in addition
\begin{equation}
\sum_{j=1}^d \|T_jf\|^2_{\cM}\le\|f\|^2_{\cM}-\|f(0)\|^2_{\cY}\quad
\mbox{for every}\quad f\in\cM
\label{6.17aa}
\end{equation}
or, equivalently, if the pair $({\bf T}, G)$ is contractive where
$G: \, \cM\to\cY$ is defined by
\begin{equation}
G\colon \; f(\blam)\to f(0).
\label{5.10a}
\end{equation}
We have the following analogue of Theorem
\ref{T:4.16} characterizing contractively included subspaces ${\mathcal
M}$ of $\cH_{\cY}(k_d)$ of the form ${\mathcal M}=\cH(K^{\ba}_{C, \bA})$;
for the general case where ${\mathcal M}$ is not ${\bf
M}_{\blam}^*$-invariant, one simply replaces ${\bf
M}_{\blam}^*$ with some contractive solution ${\bf T}$ of the Gleason
problem on ${\mathcal M}$.

           \begin{theorem}  \label{T:3nc-converse}
           Let ${\mathcal M}$ be a Hilbert space of $\cY$-valued functions
           and let us assume that there exists a contractive solution
	  ${\bf T}=(T_1,\ldots,T_d)$ of the Gleason problem (i.e.,
	  $T_{j} \in {\mathcal L}({\mathcal M})$ such that
	  \eqref{6.16a} and \eqref{6.17aa}
           hold for every $f\in{\mathcal M}$). Then ${\mathcal M}$ is
          isometrically
           equal to a reproducing kernel Hilbert space $\cH(K^{\ba}_{C, \bA})$
           for a contractive pair $(C, \bA)$.
           Therefore, ${\mathcal M}$ is contractively included in the
           Arveson space
           $\cH_{\cY}(k_d)$.
           \end{theorem}
           \begin{proof}
	Take $C = G|_{{\mathcal M}}$ where $G$ is given by \eqref{5.10a},
	$\bA = {\mathbf T}$ on
           ${\mathcal M}$.  Then \eqref{6.17aa} says that $(C, \bA)$ is
contractive.
           Iteration of \eqref{6.16a} says that, for each $f \in {\mathcal M}$,
	\begin{align*}
	 f(\blam) = & \sum_{j_{1}=1}^{d} \lambda_{j_{1}}\left[ (T_{j_{1}}f)(0)
	+ \sum_{j_{2}=1}^{d} \lambda_{j_{2}} \left[(T_{j_{2}}
	T_{j_{1}}f)(0) + \sum_{j_{3}=1}^{d} \lambda_{j_{3}}\left[(T_{j_{3}}
	T_{j_{2}}  T_{j_{1}}f(0) + \cdots \right.\right.\right.\\
	& \qquad \left.\left.\cdots + \sum_{j_{k}=1}^{d}
           \lambda_{j_{k}}\left[(T_{j_{k}}
           \cdots  T_{j_{2}}T_{j_{1}}f)(0) + \cdots \right] \cdots
           ]\right]\right].
	\end{align*}
	This unravels to the tautology
	$$ f(\blam) = C (I - Z(\blam)A)^{-1} f
	$$
	so we recover ${\mathcal M}$ as ${\mathcal M} = \operatorname{Ran}
	\widehat{\mathcal O}^{\ba}_{C, \bA}$ with $\| C (I - Z(\cdot)A)^{-1}
          f\|_{\cH(K^{\ba}_{C, \bA})} = \| f \|_{{\mathcal M}}$, i.e.,
          ${\mathcal M} = \cH(K^{\ba}_{C, \bA})$ isometrically.  From the
          fact that $(C, \bA)$ is contractive, we have seen that ${\mathcal
          G}^{\ba}_{C, \bA}\le {\mathcal G}_{C, \bA} \le I_{{\mathcal M}}$.
          Then
	$$
           \| f \|^{2}_{\cH^{2}_{\cY}(k_{d})} =\| C (I - Z(\cdot)A)^{-1} f
	  \|^{2}_{\cH^{2}_{\cY}(k_{d})} = \langle {\mathcal G}^{\ba}_{C, \bA}
           f, f\rangle_{{\mathcal M}} \le \|f \|^{2}_{{\mathcal M}}
	$$
	and we also have the contractive inclusion property.
\end{proof}

Combining Theorems \ref{T:3-1.2nc} and \ref{T:3nc-converse} gives the
following uniqueness result for contractive solutions of the Gleason
problem on a subspace ${\mathcal M}$ contained in $\cH_{\cY}(k_d)$
isometrically.

\begin{theorem}\label{P:uniquegleason}
Suppose that ${\mathcal M}$ is a subspace of $\cY$-valued functions
contained in  $\cH_{\cY}(k_d)$ isometrically and that ${\bf
T}=(T_1,\ldots,T_d)$ is a  contractive solution of the Gleason
problem on ${\mathcal M}$. Then ${\mathcal M}$ is ${\bf
M}^*_{\blam}$-invariant and ${\bf T}={\bf M}^*_{\blam}$.
\end{theorem}

\begin{proof} By Theorem \ref{T:3nc-converse}, there is a contractive pair
$(C, \, {\bf A})$ so that ${\mathcal M}=\cH(K^{\ba}_{C, \bA})$
isometrically. As ${\mathcal M}$ is contained in $\cH_{\cY}(k_d)$
isometrically, we conclude that $\cH(K^{\ba}_{C, \bA})$
is contained in $\cH_{\cY}(k_d)$ isometrically. Part $(5)$ in Theorem
\ref{T:3-1.2nc} then asserts that the subspace  ${\mathcal
M}=\cH(K^{\ba}_{C, \bA})$ is ${\bf
M}^*_{\blam}$-invariant and that $T_j=M_{\lambda_j}^*$ for $j=1,\ldots,d$.
\end{proof}

We note that the proof of Theorem \ref{T:unique} is like the proof
        of the State-Space-Isomorphism Theorem for structured noncommutative
        multidimensional linear systems in \cite{BGM2}.  It is known that
        the State-Space-Isomorphism Theorem (and related Kalman reduction
        procedure) fails in general for
        commutative multidimensional linear systems---see
        e.g.~\cite{Galkowski} for a recent account of the situation.
        The fact that uniqueness does hold in the
        special commutative situation in Theorem \ref{T:c-unique} shows
        that the technique in the proof of the State-Space-Isomorphism
        Theorem is salvageable in special commutative situations.

        A uniqueness result for solutions of the Gleason problem
        somewhat different from that in Theorem
        \ref{P:uniquegleason} was obtained in \cite{AlpayDubi}; rather
        than assuming that ${\mathbf T}$ is a contractive solution of
        the Gleason problem on $\cM = \cH(K_{C, \bA})$ contained
        isometrically in $\cH_{\cY}(k_{d})$ as in Theorem
        \ref{P:uniquegleason}, Alpay and Dubi in \cite{AlpayDubi}
        assume instead that ${\mathbf T}$ is a commutative solution of
        the Gleason problem and are then able to conclude that necessarily
        ${\mathbf T} = {\mathbf M}_{\blam}^{*}|_{\cM}$.  This latter
        result can be seen as an immediate consequence of our Theorem
        \ref{T:c-unique} above since,
        by the construction in the proof of Theorem \ref{T:3nc-converse},
        solutions   $(C,\bA)$ of
        $K^{\ba}_{C, \bA} = K$ are in one-to-one correspondence with
        solutions ${\mathbf T}$ of the Gleason problem.
We illustrate the preceding analysis by two examples.

\begin{example}\label{E:3.3.1}
{\rm Consider the subspace
$\cM={\rm span}\{1,\lambda_1,\lambda_2\}\subset\cH(k_2)$ and define
the operators $T_{a,1}$ and $T_{a,2}$ on $\cM$  by
\begin{equation}
T_{a,1}: \; f\mapsto \beta+a\alpha\lambda_2,\quad
T_{a,2}: \; f\mapsto \gamma-a\alpha \lambda_1
\label{july17}
\end{equation}
where $f(\blam)=\alpha+\beta \lambda_1+\gamma\lambda_2$ is the generic
element in $\cM$ and where $a$ is a fixed complex number. It is
readily checked that
$$
f(\blam)-f(0)=\beta \lambda_1+\gamma\lambda_2=\lambda_1(T_{a,1}f)(\blam)+
\lambda_2(T_{a,2}f)(\blam),
$$
so the tuple $(T_{a,1},T_{a,2})$ solves the Gleason problem on $\cM$. Let
$A_{a,1}$ and $A_{a,2}$ be the matrices of $T_{a,1}$ and $T_{a,2}$ with
respect to the basis $\{1,\lambda_1,\lambda_2\}$ of $\cM$ and let $C$ be
the matrix of  the operator $G:\cM\to\cY$ defined in \eqref{5.10a}:
\begin{equation}\label{july16}
C=\begin{bmatrix}1 & 0 & 0 \end{bmatrix},  \quad
A_{a,1}=\begin{bmatrix} 0 & 1 & 0 \\ 0 & 0 & 0 \\ a &0 &0
   \end{bmatrix},  \quad
A_{a,2}=\begin{bmatrix} 0 & 0 & 1 \\ -a & 0 & 0
     \\ 0 & 0 & 0\end{bmatrix}.
\end{equation}
A straightforward calculation shows that
$$
C(I-\lambda_1A_{a,1}-\lambda_2A_{a,2})^{-1}=\begin{bmatrix}1 &
\lambda_1& \lambda_2\end{bmatrix}
$$
which realizes $\cM$ as the range of the observability operator of a pair
$(C,\bA_a)$. Different choices of $a$ in \eqref{july16} lead to
non-equivalent realizations of $\cM$. Note that $A_{a,1}$ and $A_{a,2}$ do
not commute unless $a=0$, in which case the operators $T_{0,1}$ and
$T_{0,2}$ are equal to  backward shifts $M_{\lambda_1}^*$ and
$M_{\lambda_2}^*$,  respectively; in other words, the matrices
\begin{equation}\label{july15}
  C=\begin{bmatrix}1 & 0 & 0 \end{bmatrix},  \quad
     A_{0,1}=\begin{bmatrix} 0 & 1 & 0 \\ 0 & 0 & 0 \\ 0 &0 &
     0\end{bmatrix},  \quad
     A_{0,2}=\begin{bmatrix} 0 & 0 & 1 \\ 0 & 0 & 0
     \\ 0 & 0 &0\end{bmatrix}
\end{equation}
provide a commutative realization of  $\cM$ which is  unique (up to
unitary equivalence) by Theorem \ref{T:c-unique}. Note also
that the tuple $(T_{a,1},T_{a,2})$ defined in \eqref{july17} is never a
contractive solution of the Gleason problem unless $a=0$.}
\end{example}
\begin{example}\label{E:3.3.2}
{\rm  Consider the subspace
$$
\cM={\rm span}\left\{\frac{4}{4-\lambda_1\lambda_2},
\frac{\lambda_1}{4-\lambda_1\lambda_2},
\frac{\lambda_2}{4-\lambda_1\lambda_2}\right\}\subset\cH(k_2)
$$
and define the operators $T_{a,1}$ and $T_{a,2}$ on $\cM$  by
$$
T_{a,1}: \; f\mapsto
\frac{\beta+a\alpha\lambda_2}{4-\lambda_1\lambda_2},\quad
T_{a,2}: \; f\mapsto
\frac{\gamma+(1-a)\alpha\lambda_1}{4-\lambda_1\lambda_2}
$$
where $a$ s a fixed complex number and where
$$
f(\blam)=\frac{4\alpha+\beta\lambda_1+\gamma\lambda_2}{4-\lambda_1\lambda_2}
$$
is the generic element in $\cM$. Thus,
$f(0)=\alpha$ and it is readily checked that
$$
f(\blam)-f(0)=\frac{\alpha\lambda_1\lambda_2 +\beta
\lambda_1+\gamma\lambda_2}{4-\lambda_1\lambda_2}=\lambda_1(T_{a,1}f)(\blam)+
\lambda_2(T_{a,2}f)(\blam),
$$
so the tuple $(T_{a,1},T_{a,2})$ solves the Gleason problem on $\cM$.
As in the previous example, take the matrices
$$
C=\begin{bmatrix}1 & 0 & 0 \end{bmatrix},  \quad
A_{a,1}=\begin{bmatrix} 0 & \frac{1}{4} & 0 \\ 0 & 0 & 0 \\ a &0 &0
   \end{bmatrix},  \quad
A_{a,2}=\begin{bmatrix} 0 & 0 & \frac{1}{4} \\ 1-a & 0 & 0
     \\ 0 & 0 & 0\end{bmatrix}
$$
where  $A_{a,1}$ and $A_{a,2}$ are the matrices of $T_{a,1}$ and
$T_{a,2}$ with respect to the basis
$\left\{\frac{4}{4-\lambda_1\lambda_2},
\frac{\lambda_1}{4-\lambda_1\lambda_2},
\frac{\lambda_2}{4-\lambda_1\lambda_2}\right\}$ of $\cM$ and $C$ is
the matrix of  the operator $G:\cM\to\cY$ defined in \eqref{5.10a}.
For every choice of $a$,
$$
C(I-\lambda_1A_{a,1}-\lambda_2A_{a,2})^{-1}=\begin{bmatrix}
\frac{4}{4-\lambda_1\lambda_2} & \frac{\lambda_1}{4-\lambda_1\lambda_2} &
\frac{\lambda_2}{4-\lambda_1\lambda_2}\end{bmatrix}
$$
which realizes $\cM$ as the range of the observability operator of a pair
$(C,\bA_a)$. Different choices of $a$ in \eqref{july16} lead to
non-equivalent realizations of $\cM$. Note that $A_{a,1}$ and $A_{a,2}$
never commute which is not surprising since $\cM$ is not backward-shift
invariant as has been established in Example \ref{E:not-shift*-inv}.}
\end{example}

        \subsection{Applications of observability operators: the
           commutative setting}  \label{S:C-appl}

           In this subsection we discuss applications of observability
           operators for the commutative setting.
            This subsection parallels
           Subsection \ref{S:NC-appl}.

           For subspaces of $\cH_{\cY}(k_{d})$ invariant under the forward
           shift operator-tuple $M_{\blam}$, we have the following analogue
           of the Beurling-Lax-Halmos-de Branges theorem due originally
to Arveson
           \cite{ArvesonIII} and McCullough-Trent \cite{MCT} (for the
	 case of isometric inclusion); in fact, one can check that our proof,
	 namely, the commutative adaptation of the proof of Theorem
	 \ref{T:NC-BL},  follows that of \cite{Arveson} if one makes the
           substitution $L = (\widehat \cO^{\ba}_{D_{{\mathbf T}^{*}}, {\mathbf
           T}^{*}})^{*}$ (where $L$ is the key operator appearing in
           \cite{Arveson}).   In general, an operator
           $\Theta$ between two Arveson spaces $\cH_{\cU}(k_{d})$ and
           $\cH_{\cY}(k_{d})$ is said to be {\em multiplier} if $\Theta$
           intertwines the respective coordinate-function multipliers:
           $$ \theta M_{\lam_{j}} f = M_{\lam_{j}} \theta f \quad \text{for all}
           \quad f \in \cH_{\cU}(k_{d}).
           $$
           It is straightforward to see that a multiplier $\Theta$
           necessarily has the form
           $$ \Theta = M_{\theta} \colon f(z) \to \theta(z) \cdot f(z)
           $$
           where $\theta(z) = \sum_{\bn \in {\mathbb Z}^{d}_{+}} \theta_{\bn}
           z^{\bn}$ is a bounded, holomorphic $\cL(\cU, \cY)$-valued
           function on ${\mathbb B}^{d}$, but not all bounded, holomorphic,
           operator-valued functions on ${\mathbb B}^{d}$ are multipliers
           (see e.g. \cite{AglerMcCarthy}). In case the multiplication
	 operator has operator norm at most $1$, we say that $\theta$
	 is a {\em contractive multiplier} and belongs to the
	 (commutative) multivariable Schur-class $\cS_{d}(\cU, \cY)$.
           Unlike the convention in the classical case, such a multiplier
           $\theta$ is
           said to be {\em inner} if in addition $M_{\theta}$ is a partial
           isometry.

           \begin{theorem} \label{T:C-BL}
	    \begin{enumerate}
		\item A Hilbert space $\cM$ is such that
		\begin{enumerate}
		    \item $\cM$ is contractively contained in
		    $\cH_{\cY}(k_{d})$,
		    \item $\cM$ is invariant under the Arveson-shift
	$d$-tuple ${\mathbf M}_{\blam}$,
	\item the $d$-tuple
	$$ {\mathbf M}_{\cM \blam} = (M_{\cM,\lam_{1}}, \dots,
	M_{\cM, \lam_{d}}) \text{ where } M_{\cM,\lam_{j}}: =
	M_{\lam_{j}}|_{\cM} \text{ for } j = 1, \dots, d
	$$
	is a row contraction:
	$$
	 M_{\cM, \lam_{1}} (M_{\cM, \lam_{1}})^{*} + \cdots + M_{\cM,
	 \lam_{d}}(M_{\cM, \lam_{d}})^{*} \le I_{\cM},
	$$
	and
	\item $({\mathbf M}_{\cM \blam})^{*}$ is strongly stable,
	i.e.,
	$$
	  \sum_{\bn \in {\mathbb Z}^{d}_{+} \colon |\bn| = N}
	  \frac{|\bn|!}{\bn !} \left\|{\mathbf M}_{\cM,
	  \blam} )^{* \bn} f \right\|^{2}_{\cM} \to 0 \text{ as } n
	  \to \infty \text{ for all } f \in \cM
	$$
	\end{enumerate}
	if and only if there is a coefficient Hilbert space $\cU$ and
	a contractive multiplier $\theta \in \cS_{d}(\cU, \cY)$ so that
	$\cM = \theta \cdot \cH_{\cU}(k_{d})$ with lifted norm
	$$
	 \| \theta \cdot f \|_{\cM} = \| Q f \|_{\cH_{\cU}(k_{d})}
	 $$
	 where $Q$ is the orthogonal projection onto
	 $(\operatorname{Ker} M_{\theta})^{\perp} \subset
	 \cH_{\cU}(k_{d})$.

	 \item The subspace $\cM$ in part (1) above is isometrically
contained in
	 $\cH_{\cY}(k_{d})$ if and only if the corresponding
	 contractive multiplier $\theta \in \cS_{d}(\cU, \cY)$ can be
	 taken to be inner.
	\end{enumerate}

	\end{theorem}

	\begin{proof}  The proof is a straightforward commutative
	    adaptation of the proof of Theorem \ref{T:NC-BL} and
	    hence will be left to the reader.
	    We remark that, for the case where $\cM$ is contained
	    isometrically in $\cH_{\cY}(k_{d})$, we are unable to obtain
	    a representer $\theta$ for which $M_{\theta}$ is
	    isometric but rather only a representer with $M_{\theta}$
	    partially isometric.  Indeed, one can check that the
	   argument in the proof of Theorem \ref{T:NC-BL} breaks down
	   because, for the case here, $M_{\lam_{j}}$ is only
	   contractive rather than isometric.
         \end{proof}

           \begin{remark}  \label{R:C-BL} {\em As observed in
\cite{Arveson}, from the
	   function-theory point of view Theorem \ref{T:C-BL} is not a
	   true analogue of the classical Beurling-Lax theorem since the
	   characterization of $\theta$ is purely operator-theoretic with
	   no information on the boundary behavior of the associated
	   multiplier $\theta(z)$.  This deficiency has now been remedied
	   in the paper of Greene-Richter-Sundberg \cite{GRS}.}
           \end{remark}

           The following is the analogue of Theorem \ref{T:NCop-model}; we
           omit the proof as it exactly parallels the proof of Theorem
           \ref{T:NCop-model}.  The result goes back to Drury \cite{Dr}.

           \begin{theorem} \label{T:Cop-model}  Suppose that ${\mathbf T}
	       = (T_{1}, \dots, T_{d})$ is a commutative row-contractive
	       operator-tuple
	       with ${\mathbf T}^{*}$ asymptotically stable and
define the defect
	       operator $D_{{\mathbf T}^{*}}$ and the coefficient space
	       $\cY$ as in \eqref{T*defect}.   Then there is a subspace
	       $\cM \subset \cH_{\cY}(k_{d})$ invariant for the
	       backward shift operator-tuple $M_{\blam}^{*}$ on
	       $\cH_{\cY}(k_{d})$ so that ${\mathbf T}$ is unitarily
	       equivalent to $P_{\cM}M_{\blam}^{*}|_{\cM}$.  In particular,
	       ${\mathbf T}$ has a Arveson-shift dilation unitarily equivalent
	       to $M_{\blam}$ on $\cH_{\cY}(k_{d})$.
	 \end{theorem}

         As a corollary of this result one can arrive at the von Neumann
         inequality
         $$
          \| p(T_{1}, \dots, T_{d}) \|  \le \| p(M_{\lam_{1}},
           \ldots,M_{\lam_{d}})\|
          $$
          of Drury \cite{Dr} and Arveson \cite{ArvesonIII} (see Remark
          \ref{R:NC-vN} for the noncommutative case).

	 \begin{remark} \label{R:Cop-model}
       {\em The result in Theorem \ref{T:Cop-model} is tied to the unit ball
         with associated multivariable resolvent operator
         $(I - \lam_{1}T_{1}^{*} - \cdots - \lam_{d}T_{d}^{*})^{-1}$,
         associated defect operator $D_{{\mathbf T}^{*}} = (I -
         T_{1}T_{1}^{*} - \cdots - T_{d}T_{d}^{*})^{1/2}$, associated
         observability operator of the form $\widehat
\cO^{\ba}_{D_{{\mathbf T}^{*}},
         {\mathbf T}^{*}} = D_{{\mathbf T}^{*}} (I - \lam_{1}T_{1}^{*} -
\cdots - \lam_{d}
          T_{d}^{*})^{-1}$ and associated ambient kernel function
$k(\blam, \bzeta) =
          1/(1-\lam_{1}\overline{\zeta_{1}} - \cdots -
          \lam_{d} \overline{\zeta_{d}})$.
          We mention that there has been a lot of work centering around other
         types of kernels and giving a model theory for other classes of
         operator-tuples by using appropriately modified observability-like
         operators.  Specifically,
         M\"uller-Vasilescu \cite{MV} for the commutative ball case with
         $k(\blam, \bzeta) = 1/(1-\lam_{1}\overline{\zeta_{1}} -
       \cdots - \lam_{d}\overline{\zeta_{d}})^{m}$,  Curto-Vasilescu
       \cite{CV1, CV2} for the commutative polydisk case with
       $k(\blam, \bzeta) = (1/(1-\lam_{1}\overline{\zeta_{1}}) \cdots
       (1-\lam_{d}\overline{\zeta_{d}}))^{m}$, and Pott \cite{Pott} and
       Bhattacharyya-Sarkar \cite{BS}  for the commutative case with
       $k(\blam,\bzeta) = 1/(1-P(\lam_{1}\overline{\zeta_{1}}, \dots,
       \lam_{d}\overline{\zeta_{d}}))$ with $P$
       equal to a ``positively regular polynomial''. The most general
       form of results along this line is due to
       Ambrozie-Engli\v s-M\"uller \cite{AEM} and Arazy-Engli\v s
       \cite{ArazyEnglis}:}
       given a positive-definite kernel $k(\blam, \bzeta)$ on a domain
       ${\mathcal D} \subset {\mathbb C}^{d}$ and a $d$-tuple of
       operators ${\mathbf T} = (T_{1}, \dots, T_{d})$ with Taylor
       spectrum contained in $\overline{\mathcal D}$ for which one can
       make sense of the defect operator $D_{{\mathbf T}^{*}} :=
       \frac{1}{k}(T,T)$ and of the observability operator
       $$
       \cO_{D_{{\mathbf T}^{*}}, {\mathbf T}^{*}} \colon x \mapsto
       D_{{\mathbf T}^{*}} k(\blam, T)
       $$
       {\em (for example, if $k(\blam, \bzeta)$ has no zeros in
       ${\mathcal D} \times {\mathcal D}$ and ${\mathbf T}$ has Taylor spectrum
       contained in ${\mathcal D}$)}, then, under the assumption that
       $D_{{\mathbf T}^{*}} \ge 0$ and that an additional stability
       condition on ${\mathbf T}^{*}$ holds, ${\mathcal O}_{D_{{\mathbf
       T}^{*}}, {\mathbf T}^{*}}$ implements a unitary equivalence
       between ${\mathbf T}$ and $P_{\cM} M_{\blam} |_{\cM}$, {\em where
       $$
       \cM = \operatorname{Ran} \cO_{D_{{\mathbf T}^{*}}, {\mathbf T}^{*}}
       \subset \cH(k) \otimes \cY \text{ with }
       \cY : = \overline{\operatorname{Ran}} D_{{\mathbf T}^{*}},
       $$
       where $M_{\blam} = (M_{\lam_{1}}, \dots, M_{\lam_{d}})$ is the
       operator-tuple of multiplication by the coordinate functions on
       $\cH(k) \otimes \cY$, and where $\cM$ is invariant under each of
       $M_{\lam_{1}}^{*}, \dots, M_{\lam_{d}}^{*}$.
       The noncommutative case is not as well developed at this writing,
       but there is the paper of Popescu \cite{popjfa} which handles the
       case of a Cartesian product of noncommutative balls (and therefore
       including a noncommutative polydisk).
       We expect that many of the ideas of the present paper, including the
       interplay between the noncommutative and commutative settings and
       the connections with system theory, have some parallels in these
       other situations.}
       \end{remark}

\end{document}